\documentclass[smallextended,numbook,runningheads]{svjour3}     % onecolumn (second format)
\smartqed  % flush right qed marks, e.g. at end of proof\usepackage{amssymb}
\usepackage{xspace}
\usepackage{subfigure}
\usepackage{amsfonts}
\usepackage{amsmath}
\usepackage{enumerate}
\usepackage{framed}
\usepackage{graphicx}
\usepackage{lscape}
\usepackage{bm}
\usepackage{color}
\usepackage{multicol}
\usepackage{cite}
\newcommand{\secref}[1]{{Section~\ref{#1}}}
\renewcommand{\eqref}[1]{{(\ref{#1})}}

\newcommand{\R}{\mathbb{R}}
\newcommand{\N}{\mathbb{N}}
\newtheorem{thm}[theorem]{Theorem}
\newtheorem{prp}[theorem]{Proposition}
\newtheorem{dfn}[theorem]{Definition}
\newtheorem{ass}[theorem]{Assumption}
\newtheorem{lem}[theorem]{Lemma}
\newcommand{\thmref}[1]{{Theorem~\ref{#1}}}
\newcommand{\lemref}[1]{{Lemma~\ref{#1}}}
\newcommand{\assref}[1]{{Assumption~\ref{#1}}}
\newcommand{\propref}[1]{{Proposition~\ref{#1}}}
\usepackage{graphicx}
\begin{document}
	
	\title{Strong convergence of some  Magnus-type schemes  for the finite element discretization of non-autonomous parabolic SPDEs driven by additive fractional Brownian motion and  Poisson random measure}
	
	%placed here. 
	% General acknowledgments should be placed at the end of the article.
	%}
	
	
	\titlerunning{ Magnus-type schemes for non-autonomous SPDEs driven by additive noises}        % if too long for running head
	
	\author{Aurelien Junior Noupelah \and Jean Daniel Mukam        \and
		Antoine Tambue %etc.
	}
	

	\institute{A. J. Noupelah \at The African Institute for Mathematical Sciences (AIMS) of Cameroon, P.O. Box 608, Crystal Gardens, Limbe, Cameroon.\\
	Equipe MASS, Laboratoire de Math\'ematiques,  UFD MIA, EDOSFA, Universit\'e de Douala, P.O. BOX 24157, Douala, Cameroun.\\
	\email{aurelien.noupelah@aims-cameroon.org, noupsjunior@yahoo.fr}           %  \\
	\and
	J. D. Mukam \at
	Department of Mathematics, Bielefeld University, 33615 Bielefeld, Germany.\\
	% Fax: +123-45-678910\\
	\email{jmukam@math.uni-bielefeld.de, jean.d.mukam@aims-senegal.org}           %  \\
	\and
	A. Tambue \at
	Department of Computing Mathematics and Physics,  Western Norway University of Applied Sciences, Inndalsveien 28, 5063 Bergen.\\
	%Center for Research in Computational and Applied Mechanics (CERECAM), and Department of Mathematics and Applied Mathematics, University of Cape Town,
	%7701 Rondebosch, South Africa. \\
	%The African Institute for Mathematical Sciences(AIMS), 6-8 Melrose Road, Muizenberg 7945, South Africa\\
	\email{antonio@aims.ac.za}
	}
	
	\date{Received: date / Accepted: date}

	\maketitle
	
\begin{abstract}
The aim of this work is to provide the  strong convergence results of numerical approximations of a general second order non-autonomous semilinear stochastic partial differential equation (SPDE) driven simultaneously by an additive fractional Brownian motion (fBm) with Hurst parameter $H\in(\frac 12, 1)$  and a Poisson random measure, more realistic in modelling real world phenomena.  Approximations in space are  performed by the  standard finite element method and in time by the stochastic Magnus-type integrator  or the  linear semi-implicit Euler method.  We investigate the  mean-square errors estimates of our fully discrete schemes and   the results show how  the convergence orders depend on the regularity of the initial data and the driven processes.  To the best of our knowledge, these two schemes are the first numerical methods  to 
approximate the  non-autonomous  semilinear stochastic partial differential equation (SPDE) driven simultaneously by an additive fractional Brownian motion with Hurst parameter $H$  and a Poisson random measure.
%These theoretical findings are accompanied by several numerical examples.
\keywords{Non-autonomous stochastic parabolic partial differential equations \and  Fractional Brownian motion  \and Poisson random measure \and Finite element method \and Strong convergence \and  Magnus-type integrator  \and semi-implicit Euler method.}
%SPDE, Fractional Brownian motion, Additive noise, Finite element method, Linear implicit scheme, Exponential integrator-Euler scheme, Exponential Rosenbrock-Euler scheme, Strong convergence
\subclass{65C30 \and 65J08  \and 65M60  \and 65M12 \and 65M15}
\end{abstract}
%\begin{linenumbers}
\section{Introduction}
\label{intro}
We analyse the strong numerical approximations of the second order semilinear stochastic partial differential equation which can be written  in the following non-autonomous SPDE 
%defined in $\Lambda\subset\R^d$, $d \in \{1,2,3\}$
%Dirichlet, Neumann, Robin boundary conditions or mixed Dirichlet and Neumann). 
  \begin{eqnarray}
\label{SPDE}
  \left\{
    \begin{array}{ll}
       dX(t)+A(t)X(t)dt=F(X(t))dt+\phi(t) dB^H(t)+\int_{\chi}\psi(z)\widetilde{N}(dz,dt),\quad t\in[0,T] \\
      X(0)=X_0
    \end{array}
  \right.
  \end{eqnarray}
  in the Hilbert space $\mathcal{H}=L^2(\Lambda,\R)$, where $\Lambda \subset\R^d$, $d \in \{1,2,3\}$ is bounded and has smooth boundary. $T>0$ is the final time.
  The mark set $\chi$ is defined by $\chi:=\mathcal{H}\setminus\{0\}$.  For any topological set $\Gamma$,  let $\mathcal{B}(\Gamma)$ be the smallest $\sigma$-algebra containing all open sets of $\Gamma$. Let $(\chi, \mathcal{B}(\chi), \nu)$ be a $\sigma$-finite measurable space and $\nu$ ( with $\nu\not\equiv 0$) a L\'{e}vy measure on $\mathcal{B}(\chi)$ such that
\begin{eqnarray*}
\nu(\{0\})=0\quad \text{and}\quad \int_{\chi}\min(\Vert z\Vert^2, 1)\nu(dz)<\infty.
\end{eqnarray*}
 Let $N(dz, dt)$ be the $\mathcal{H}$-valued  Poisson  distributed $\sigma$-finite measure on the product $\sigma$-algebra $\mathcal{B}(\chi)$ and $\mathcal{B}(\mathbb{R}_{+})$ with intensity $\nu(dz)dt$, where $dt$ is the Lebesgue measure on $\mathcal{B}(\mathbb{R}_{+})$. In  \eqref{SPDE}, $\widetilde{N}(dz, dt)$ stands for the compensated Poisson random measure defined by
 \begin{eqnarray*}
 \widetilde{N}(dz, dt):=N(dz, dt)-\nu(dz)dt.
\end{eqnarray*} 
%on the Hilbert space $\mathcal{H}=L^2(\Lambda)$.
 $F:[0,T]\times\mathcal{H}\rightarrow \mathcal{H}$, $\phi:[0,T]\rightarrow L^0_2$ and $\psi:\chi\rightarrow \mathcal{H}$ are deterministic mappings that will be specified  later, $X_0$ is the initial data which is random, and for all $s\in[0,T]$, $-A(s)$ is a second order linear operator, not necessarily self-adjoint, unbounded  and generator of an analytic semigroup $S_s(t):=e^{-tA(s)}$, $t\geq 0$ (see Section 2.2 for more details). The term  $B^H(t)$ in \eqref{SPDE} is an $\mathcal{H}$-valued $Q$-cylindrical fractional Brownian motion with Hurst parameter $H\in(\frac 12,1)$ in a filtered probability space $(\Omega, \mathcal{F}, \mathbb{P},\{\mathcal{F}_t\}_{t\geq 0})$. The  covariance operator $Q: \mathcal{H}\rightarrow \mathcal{H}$ of $B^H(t)$ is assumed to be positive definite and  self-adjoint. 
The filtered probability space $(\Omega, \mathcal{F}, \mathbb{P},\{\mathcal{F}_t\}_{t\geq 0})$ is assumed to fulfil the usual conditions
(see \cite[Definition 2.2.11]{Pre}). It is well known \cite{Car} that the noise can be represented as
\begin{eqnarray*}
%\label{def_fBm}
B^H(t)=\sum_{i\in\N^d}\beta^H_i(t)Q^{\frac 12}e_i=\sum_{i\in\N^d}\sqrt{q_i}\beta^H_i(t)e_i,\,\,\,t\in[0,T].
\end{eqnarray*}
where $q_i,\,e_i,\,i\in\N^d$ are respectively the eigenvalues and eigenfunctions of the covariance operator $Q$, and $\beta_i^H$ are mutually independent and identically distributed fractional Brownian motions (fBm). The cylindrical fBm $B^H$ and the compensated Poisson random measure $\widetilde{N}$ are assumed to be independent.
 
The self-similar and long-range dependence properties of the fBm make this process a suitable candidate to model  many phenomena  like financial markets (see e.g., \cite{Che,Hu, Man}) and traffic networks (see e.g., \cite{Lel,Wil}). In many real life situations such as in finance, the unpredictable nature of many events 
 such as market crashes, announcements made by the central banks,  changing    credit risk, insurance in a changing risk, changing face of operational risk \cite{Con,Pla} might have sudden and significant impacts on the stock price. When modelling such phenomena,  it is recommended to  incorporate a non-Gaussian noise such as L\'evy process or Poisson random measure in order to obtain a more realistic model.
%  Deterministic counter part of \eqref{SPDE} is used to model many real world phenomena in many fields such as quantum fields theory, electromagnetism, nuclear physics and transport in porous media.
 In most cases,  analytical solutions for SPDEs of type \eqref{SPDE} are rarely known explicitly and therefore it is important to construct discrete approximations. Due to the complexity of the linear operator $A(s)$, its semigroup $S_s(t)$ and the resolvent operator $(I+tA(s))^{-1}$, $t,s\in[0,T]$, numerical approximations
 of non-autonomous SPDEs are not well understood in the literature. To the best of our knowledge, the first numerical method  to approximate the  deterministic non-autonomous PDE was presented in\cite{Mag},  then over time some modified versions of that Magnus scheme have been proposed  in \cite{Bla1,Bla2,Bla3,Gon,Hip,Hoc}. It is just very recently that few  works has been done on  numerical methods for non-autonomous SPDEs driven by cylindrical Brownian motion (see e.g., \cite{Mukmaa,Mukmab,Tamma,Tammab}). 
 
  It is important to  mention that  if $H\neq \frac 12$ the process $B^H$ is not a semi-martingale  and the standard stochastic calculus techniques are therefore obsolete while studying SPDEs of type \eqref{SPDE}.  Alternative approaches to the standard It\^o calculus are therefore required in order to build a stochastic calculus framework for such fBm.   In recent years, there have been various developments of stochastic calculus and stochastic differential equations with respect to 
  the fBm especially for $H\in(\frac 12, 1)$ (see, for example \cite{Alo,Car,Mis})
%More precisely , to build a mathematical setting for such fBm,
and  theory of SPDEs driven by fractional Brownian motion  has been also studied. For example, linear and semilinear stochastic equations in a Hilbert space with an infinite dimensional fractional Brownian motion are considered in \cite{Duna,Dunb}. In contrast to standard Brownian ($H=1/2$), where there are numerous literature on numerical algorithms for autonomous SPDEs, few works have been done  for numerical methods for fBm for autonomous SPDEs of type \eqref{SPDE}.
Indeed, standard explicit, linear implicit and stochastic exponential schemes have been investigated in the literature for autonomous SPDEs of type \eqref{SPDE} (see \cite{Kam,Kim,Wanc} when $\psi(z)=0$ and \cite{Nou} when $\psi(z)=z_0\neq 0$). 
The work in \cite{Nou} deals with autonomous SPDEs with constant linear operator $A(t) = A$.  This is very restrictive as modelling real world phenomena with
time dependent linear operator is more realistic than modelling with time independent linear operator (see e.g. \cite{Bla3} and references therein). In contrast to autonomous SPDE driven by fBm  where at least few numerical schemes exist, numerical schemes  for such non-autonomous SPDE of type \eqref{SPDE} driven by fBm and Poisson measure  have been lacked in scientific literature, to the best of our knowledge.
Our goal in this work is to build  the first  two numerical methods  to approximate the  non-autonomous  semilinear stochastic partial differential equation (SPDE)  \eqref{SPDE} driven simultaneously by an additive fractional Brownian motion with Hurst parameter $H\in(\frac 12, 1)$  and a Poisson random measure by proposing stochastic Magnus-type  scheme and stochastic  implicit  scheme  for fBm and analyzing their strong convergences.
% updating  works in  \cite{Mukmab,Tammab}.
%These theoretical findings are accompanied by several numerical examples.
%Our goal in this work  is to extend keys time stepping methods, which have been built in \cite{Mukmab,Tammab} for standard Brownian motion ($H=1/2$). 
%To the best of our knowledge, the numerical schemes presented in this paper  are the first numerical methods for non-autonomous SPDEs driven simultaneously by fBm and Poisson random measure.
 Strong convergence analysis of  such schemes toward the mild solution are not straightforward since  the  the process $B^H$ is not a semi-martingale. 
Our results will be based on many novel intermediate lemmas. Indeed, our  schemes here are based on finite element method for spatial discretization so that we gain the flexibility of 
these methods to deal with complex boundary conditions and we can apply well-developed techniques such as upwinding to deal with advection.
%The main idea of this article is to consider the linear operator not necessarily self-adjoint, to use the finite element method for space discretization and make a comparative study of schemas i
 %For time discretization,  we first provide the strong convergence of the stochastic Magnus-type integrator scheme \cite{Mukmaa,Tamma}.
% We also provide the  strong convergence  of the Magnus implicit scheme \cite{Mukmab}.
  Our results on strong convergence demonstrate how the convergence orders depend on the regularity of the noise and the initial data. Specifically, the semi-implicit Euler method achieves suboptimal convergence rates of $\mathcal{O}(h^{2-\epsilon})$ in space and $\mathcal{O}(\Delta t^{1-\epsilon})$ in time, where $\epsilon$ is a small positive number. In contrast, the stochastic Magnus-type integrator scheme achieves an optimal rate of $\mathcal{O}(\Delta t)$.
 % , .
% 
% Recently, stochastic exponential integrators [14, 22, 37] wereappeared as non standard explicit methods efficient for SPDE (1). All stochas-tic exponential integrators analyzed in the literature for SPDEs [14, 22, 37] arebounded on the nonlinear problem as in (1) where the linear partAand thenonlinear functionFare explicitly known a priori. Such approach is justifiedin situations where the nonlinear functionFis small. Indeed whenFis smallthe linear operatorAdrives the SPDE (1) and the good stability of the ex-ponential integrators and semi-implicit method are ensured. In fact, in morerealistic applications the nonlinear functionFcan be stronger. Typical exam-ples  are  semilinear  advection  diffusion  reaction  equations  with  stiff  reactionterm. In such cases, the SPDE (1) is driven by the nonlinear operatorFandboth exponential integrators [14,22,37] and semi-implicit Euler [26] will behaveas explicit Euler-Maruyama scheme, therefore their good stability propertiesare lost. To overcome this issue we propose in this work a new scheme calledStochastic Exponential Rosenbrock Scheme (SERS

The rest of the paper is structured as follows. In Section \ref{preli}, Mathematical settings for cylindrical fBm and random Poisson measure are presented, along with the well posedness and regularities results of the mild solution of SPDE \eqref{SPDE}.  In Section \ref{numschemes}, numerical schemes based on stochastic Magnus-integrator  and   implicit scheme for SPDE \eqref{SPDE} are presented. We give some regularity estimates of the semi-discrete problem in \secref{conv}, along with the space error.  In Section \ref{convSMTI}, the strong convergence proof of the stochastic Magnus-type integrator for fBm is presented. We end the paper in \secref{convimpl} by presenting the convergence analysis of the stochastic implicit scheme for fBm.
%Space strong convergence proof of schemes is presented in Subsection \ref{convspa}, time strong convergence proof of SMTI scheme is presented in Subsection \ref{convSMTI} and time strong convergence proof of Magnus implicit scheme is presented in Subsection \ref{convimpl}. 
%We end the paper in Section \ref{numerik} with numerical experiments to illustrate our theoretical findings.

\section{Mathematical setting and well posedness}
\label{preli}
Let us start by reviewing briefly some standard results on fractional calculus and presenting notations, main function spaces and norms that will be used in this paper.
\begin{dfn}\textbf{[Fractional Brownian motion]}\cite{Dio,Mas,Mis,Nou,Wanb}
\label{fBm}An $\mathcal{H}$-valued Gaussian process $\{B^H(t),t\in[0,T]\}$ on $(\Omega, \mathcal{F}, \mathbb{P},\{\mathcal{F}_t\}_{t\geq 0})$ is called a fractional Brownian motion with Hurst parameter $H\in(0,1)$ if 
\label{definfBm}
\begin{itemize}
\item $\mathbb{E}[B^H(t)]=0$ for all $t\in\R$,
\item $\text{Cov}(B^H(t),B^H(s))=\frac 12\left(|t|^{2H}+|s|^{2H}-|t-s|^{2H}\right)$ for all  $t,s\in\R$,
\item $\{B^H(t),t\in[0,T]\}$ has continuous sample paths  $\mathbb{P}$ a.s.,
\end{itemize} 
where $\text{Cov}(X,Y)$ denotes the covariance operator for the Gaussian random variables $X$ and $Y$ and $\mathbb{E}$ stands for the mathematical expectation on $(\Omega, \mathcal{F}, \mathbb{P},\{\mathcal{F}_t\}_{t\geq 0})$. 
\end{dfn}  
Notice that if $H=\frac 12$, the fractional Brownian motion coincides with the standard Brownian motion. Throughout this paper the Hurst parameter $H$ is assumed to belong to $(\frac 12,1)$.  

Let $\left(K,\langle .,.\rangle_K,\|.\|\right)$ be a separable Hilbert space. For $p\geq 2$ and for a Banach space U, we denote by $L^p(\Omega,U)$ the Banach space of $p$-integrable $U$-valued random variables. 
We denote by $L(U,K)$ the space of bounded linear mapping from $U$ to $K$ endowed with the usual operator norm $\|.\|_{L(U,K)}$ and $\mathcal{L}_2(U,K)=HS(U,K)$ the space of Hilbert-Schmidt operators from $U$ to $K$ equipped  with the following norm
\begin{eqnarray}
\label{normL2UK}
\left\|l\right\|_{\mathcal{L}_2(U,K)}:=\left(\sum_{i\in\N^d}\|l\psi_i\|^2\right)^{\frac 12},\,\,\,\,\,l\in \mathcal{L}_2(U,K),
\end{eqnarray}
where $(\psi_i)_{i\in\N^d}$ is an orthonormal basis on $U$. The sum in $\eqref{normL2UK}$ is independent of the choice of the orthonormal basis of $U$. We use the notation $L(U,U)=:L(U)$ and $\mathcal{L}_2(U,U)=:\mathcal{L}_2(U)$. It is well known that for all $l\in L(U,K)$ and $l_1\in \mathcal{L}_2(U)$, $ll_1\in \mathcal{L}_2(U,K)$ and 
\begin{eqnarray*}
%\label{inL2UK}
\|ll_1\|_{\mathcal{L}_2(U,K)}\leq \left\|l\right\|_{L(U,K)}\|l_1\|_{\mathcal{L}_2(U)}.
\end{eqnarray*} 
We denote by $L^0_2:=HS(Q^{\frac 12}(\mathcal{H}),\mathcal{H})$ the space of Hilbert-Schmidt operators from $Q^{\frac 12}(\mathcal{H})$ to $\mathcal{H}$ with corresponding norm $\|.\|_{L^0_2}$ defined by 
\begin{eqnarray}
\label{normL02}
\|l\|_{L_2^0}:=\left\|lQ^{\frac 12}\right\|_{HS}=\left(\sum_{i\in\N^d}\|lQ^{\frac 12}e_i\|^2\right)^{\frac 12},\,\,\,\,\,l\in L_2^0,
\end{eqnarray}
where $(e_i)_{i\in\N^d}$ is an orthonormal basis of $\mathcal{H}$. The sum in \eqref{normL02} is also independent of the choice of the orthonormal basis of $\mathcal{H}$. Let $L^2_{\nu}(\chi\times[0,T];\mathcal{H})$ be the space of all mappings $\theta: \chi\times[0,T]\times \Omega\rightarrow\mathcal{H}$ such that $\theta$ is jointly measurable and $\mathcal{F}_t$-adapted for all $z\in\chi$, $0\leq s\leq T$ satisfying
\begin{eqnarray*}
%\label{wellposjump}
\int_0^T\int_{\chi}\left\|\theta(z,s)\right\|^2\nu(dz)ds<\infty.
\end{eqnarray*}
\begin{lem}
\label{int}
\begin{enumerate}
\item [(i)] Let $\phi\in L^2([0,T];L^0_2)$, then the following holds
\begin{eqnarray}
\mathbb{E}\left[\left\|\int_0^T\phi(s)dB^H(s)\right\|^2\right]&\leq& C(H)\sum_{i\in\N^d}\left(\int_0^T\left\|\phi(s)Q^{\frac 12}e_i\right\|^{\frac 1H}ds\right)^{2H}\label{fracint1}\\
&\leq&C(H,T)\int_0^T\left\|\phi(s)\right\|^2_{L^0_2}ds.\label{fracint2}
\end{eqnarray}
\item [(ii)] Let $\theta\in L^2_{\nu}(\chi\times[0,T];\mathcal{H})$, then the following holds
\begin{eqnarray}
\label{jumpint}
\mathbb{E}\left[\left\|\int_0^T\int_{\chi}\theta(z,s)\tilde{N}(dz,ds)\right\|^2\right]=\int_0^T\int_{\chi}\left\|\theta(z,s)\right\|^2\nu(dz)ds.
\end{eqnarray}
\end{enumerate}
\end{lem}
\textit{Proof:} The equality \eqref{jumpint} is the It\^o isometry property in the jump case and its proof can be found in \cite[(3.56)]{Man}. The proof of (i) can be found in \cite[(9)-(11)]{Nou}. Note that in the case $H=\frac 12$ the constants $C(H)$ and $C(H,T)$ in \eqref{fracint1} and \eqref{fracint2} are 1 and the inequalities become the equalities. In that case they are also known as  It\^o isometry, see \cite[(4.30)]{Pra}.

More details on the definition of stochastic integral with respect to $Q$-fractional Brownian motion,  the compensated Poisson random measure and their properties are given in e.g \cite{Alo,Duna,Dunb,Mis,Man}.  To ensure the existence and the uniqueness of solution of \eqref{SPDE} and for the purpose of the convergence analysis, we make the following assumptions.  
\begin{ass}
\label{linear}
\begin{enumerate}
\item [(i)] We assume that $\mathcal{D}(A(t))=\mathcal{D}(A(s))=:\mathcal{D}$, $0\leq s, t\leq T$ and  the family of linear operators $A(t):\mathcal{D}\subset\mathcal{H}\rightarrow\mathcal{H}$ is uniformly sectorial on $0\leq t\leq T$, i.e. there exist constants $C>0$ and $\theta\in\left(\frac 12\pi;\pi\right)$ such that 
\begin{eqnarray*}
%\label{linear1}
\left\|\left(\lambda I-A(t)\right)^{-1}\right\|_{L(L^2(\Lambda))} \leq \frac C{{\left|\lambda\right|}},\quad\lambda\in S_{\theta},
\end{eqnarray*}
where $S_{\theta}:=\{\lambda\in\mathbb{C},\lambda=\rho e^{i\phi},\rho>0,0\leq|\phi|\leq \theta\}$. As in \cite{Hip}, by a standard scaling argument, we assume $-A(t)$ to be invertibe with bounded inverse.
\item[(ii)]  There exists a positive constant $K_1$ such that 
\begin{eqnarray*}
\left\|\left(A(t)-A(s)\right)\left(A(0)\right)^{-1} \right\|_{L(\mathcal{H})}&\leq& K_1|t-s|,\quad c,t\in[0;T],\label{linear2}\\
\left\|\left(A(0)\right)^{-1}\left(A(t)-A(s)\right) \right\|_{L(\mathcal{D},\mathcal{H})}&\leq& K_1|t-s|,\quad c,t\in[0;T].\label{linear3} 
\end{eqnarray*}  
\item [(iii)] Since we are dealing with non smooth data, we also assume that 
\begin{eqnarray}
\label{linear4}
\mathcal{D}\left((A(t))^{\alpha}\right)=\mathcal{D}\left((A(0))^{\alpha}\right),\quad \alpha\in[0,1],
\end{eqnarray}
and there exists a positive constant $K_2$ such that the following estimate holds
\begin{eqnarray}
\label{linear5}
K_2^{-1}\left\|(A(0))^{\alpha}u\right\|\leq\left\|(A(t))^{\alpha}u\right\|\leq K_2\left\|(A(0))^{\alpha}u\right\|,\quad u\in\mathcal{D}\left((A(0))^{\alpha}\right),\quad t\in[0;T]. 
\end{eqnarray}
\end{enumerate}
\end{ass}
\begin{ass}
\label{initdata}The initial data $X_0:\Omega\rightarrow \mathcal{H}$ is assumed to be $\mathcal{F}_0/\mathcal{B}(\mathcal{H})$ measurable and $X_0\in L^2\left(\Omega,\mathcal{D}\left((A(0))^{\frac{2H+\beta-1}2}\right)\right)$ with $0<\beta\leq 1$.
\end{ass}
\begin{ass}
\label{noiseandjump} We assume that for some constant $\delta\in\left[\frac{2H+\beta-1}2;1\right]$, the deterministic mapping $\phi:[0;T]\times\Omega\rightarrow L^0_2$ and the jump function $\psi(z)$ satisfy the following estimates
\begin{eqnarray}
\left\|(A(0))^{\frac{\beta-1}2}\phi(t)\right\|_{L^0_2}\leq C<\infty,\label{noise1}\quad
\int_{\chi}\left\|(A(0))^{\frac{2H+\beta-2}2}\psi(z)\right\|^2\nu(dz)\leq C<\infty,\label{jump}\\
\left\|(A(0))^{\frac{\beta-1}2}(\phi(t_2)-\phi(t_1))\right\|_{L^0_2}\leq C(t_2-t_1)^{\delta},\quad 0\leq t_2\leq t_1\leq T.\label{noise2}
\end{eqnarray}
\end{ass}
Next as in \cite{Sei,Mukmab}, we make the following assumptions on operators $F$ and $A(t)$.
\begin{ass}
\label{nonlinear}
The nonlinear operator $F:[0,T]\times\mathcal{H}\rightarrow\mathcal{H}$ is $\frac{2H+\beta-1}2$-H\"{o}lder continuous with respect to the first variable and Lipschitz continuous with respect to the second variable,  i.e. there exists a constant $L$ such that
\begin{eqnarray*}
%\label{nonlin}
\left\|F(s,0)\right\|\leq L, \; \left\|F(t,u)-F(s,v)\right\|\leq L\left(\left|t-s\right|^{\frac{2H+\beta-1}2}+\left\|u-v\right\|\right), \; t,s\in[0,T],\; u,v\in\mathcal{H}.
\end{eqnarray*}
\end{ass} 

Next we introduce the space related to the fractional powers of the linear operators $A(t)$. For $\alpha\in\R$, define $\mathcal{\dot H}_{t}^{\alpha}=\mathcal{D}\left((A(t))^{\frac{\alpha}2}\right)$ with the norm $\left\|\cdot\right\|_{\alpha,t}=\left\|(A(t))^{\frac{\alpha}2}\cdot\right\|$.
Due to \eqref{linear4}, \eqref{linear5} and in means of simplicity, we  write $\mathcal{\dot H}^{\alpha}:=\mathcal{\dot H}_t^{\alpha}$ and $\left\|\cdot\right\|_{\alpha}:=\left\|\cdot\right\|_{\alpha,t}$. 
\begin{prp}(\cite[Remark 2.1]{Mukmab})
\label{semigroup} Let $\alpha\geq 0$ and $0\leq \gamma \leq 1$, then there exist a constant $C_1>0$ such that
\begin{eqnarray*}
\left\|(A(t))^{\alpha}e^{-sA(t)}\right\|_{L(\mathcal{H})}\leq C_1 s^{-\alpha},\quad \label{semigroup1}\;
\left\|(A(t))^{-\gamma}\left(I-e^{-rA(t)}\right)\right\|_{L(\mathcal{H})}\leq C_1 r^{\gamma},\quad s>0,\quad r\geq0.\label{semigroup2} 
\end{eqnarray*}
\end{prp}
\begin{prp}(\cite[Remark 2.2]{Mukmab})
\label{evolution} Let $\Delta(T):=\{(t,s), 0\leq s\leq t\leq T\}$. Under \assref{linear} there exists a unique evolution system \cite[Definition 5.3, Chapter 5]{Pazy} $U:\Delta(T)\rightarrow L(\mathcal{H})$ satisfying:  
\begin{enumerate}
\item [(i)] There exists a positive constant $K_0$ such that 
\begin{eqnarray*}
%\label{evolution1}
\left\|U(t,s)\right\|_{L(\mathcal{H})}\leq K_0,\quad 0\leq s\leq t\leq T.
\end{eqnarray*}
\item [(ii)] $U(\cdot,s)\in\mathcal{C}^1( ]s,T];L(\mathcal{H}))$, $0\leq s\leq T$,
\begin{eqnarray*}
%\label{evolution2}
\frac{\partial U}{\partial t}(t,s)=A(t)U(t,s),\quad \left\|A(t)U(t,s)\right\|_{L(\mathcal{H})}\leq \frac{K_0}{t-s},\quad 0\leq s<t\leq T. 
\end{eqnarray*}
\item [(iii)] $U(t,\cdot)v\in\mathcal{C}^1([0,t[;\mathcal{H})$, $0<t\leq T$, $v\in\mathcal{D}(A(0))$ and 
\begin{eqnarray*}
\label{evolution3}
\frac{\partial U}{\partial s}(t,s)v=-U(t,s)A(s)v,\quad \left\|A(t)U(t,s)(A(s))^{-1}\right\|_{L(\mathcal{H})}\leq K_0,\quad 0\leq s\leq t\leq T. 
\end{eqnarray*}
\end{enumerate}
\end{prp}
The following assumption will be needed in our convergence estimate to achieve optimal convergence order in time without any logarithmic reduction when dealing with Magnus-integrator scheme.
\begin{ass}
\label{splitassumption}
Let $A(t)=A^s(t)+A^{ns}(t)$, where $A^s(t)$ and $A^{ns}(t)$ are respectively the self-adjoint and the non self-adjoint parts of $A(t)$.
We assume that  the family $(\lambda_n(t))_{n\in \mathbb{N}}$ 
of positive eigenvalues of $-A^s(t)$ corresponding to the eigenvectors $(e_n(t))_{n\in \mathbb{N}}$ are such that for $v\in H$  
\begin{eqnarray}
 \label{supcond}
 \underset{0\leq t\leq T}{\sup} \lambda_n(t) < C(n), \,\quad  \underset{0\leq t\leq T}{\sup} (e_n(t),v) < C_1(v,n).
 \end{eqnarray}
 where  $C (n)$ and $C_1(v, n)$  are two positive constants.
 \end{ass}

Note that if equation \eqref{SPDE} is only driven by a fractional Brownian motion, then the above assumptions are sufficient to achieve a convergence order greater than $\frac 12$ (see \cite{Nou}). However since our problem \eqref{SPDE} is also driven by a random Poisson measure we require additional following assumption as in \cite{Muka,Mukmaa,Mukmab,Nou}.
\begin{ass}
\label{derivate} The nonlinear function $F:[0,T]\times\mathcal{H}\rightarrow\mathcal{H}$ is differentiable with respect to the second variable, with bounded derivative, ie. there exists $C_2\geq 0$ such that
\begin{eqnarray*}
%\label{derivate1}
\left\|F'(t,u)v\right\|\leq C_2\left\|v\right\|,\quad t\in[0,T],\quad u,v\in\mathcal{H},  
\end{eqnarray*}
where $F'(t,u)=\frac{\partial F}{\partial u}(t,u)$ for $t\in[0;T]$ and $u\in \mathcal{H}$.

Moreover, the  derivative operator $F'$ satisfies the following estimates
\begin{eqnarray*}
\Vert F'(s,0)\Vert_{L(\mathcal{H})} \leq K_1, \; \Vert (-A(t))^{-\eta}\left(F'(t, u)-F'(t,v)\right)\Vert_{L(\mathcal{H})}\leq K_1\Vert u-v\Vert,
\end{eqnarray*}
for all $s, t\in[0,T]$ and $u, v\in H$ and for some $\eta\in\left(\frac{3}{4},1\right)$, where $F'(t,u):=\frac{\partial F}{\partial u}(t, u)$, $u \in \mathcal{H}$, $t\in[0, T]$.
\end{ass}   
The following definition is an extension of \cite[Definition 1.1]{Sei} in a Jump case.
\begin{dfn}
\label{mild}
A $\mathcal{H}$-valued $\{\mathcal{F}_t\}$-adapted measurable stochastic process $\{X(t), 0\leq t\leq T\}$ is called mild solution to the problem \eqref{SPDE} on $[0;T]$, if
\begin{enumerate}
\item [(i)] $\{X(t), 0\leq t\leq T\}$ is joint measurable and $\mathbb{E}\left[\int_0^T\left\|X(s)\right\|^2 ds\right]<\infty$.
\item [(ii)] For arbitrary $t\in[0,T]$,
\begin{eqnarray}
\label{solpb}
X(t)&=&U(t,0)X_0+\int_0^tU(t,s)F(s,X(s))ds+\int_0^tU(t,s)\phi(s)dB^H(s)\nonumber\\
&+&\int_0^t\int_{\chi}U(t,s)\psi(z)\tilde{N}(dz,ds)
\end{eqnarray}
holds almost surely, where $U(t,s)$ is the evolution system defined in \propref{evolution}.
\end{enumerate}
\end{dfn}
The well posedness result is given in the following theorem.
\begin{thm}
\label{wellposedness}
Assume that Assumptions \ref{noiseandjump}, \ref{nonlinear}, \ref{linear} (i), (ii) are satisfied. Then the non-autonomous problem \eqref{SPDE} admit a unique mild solution $X(t)\in(\mathcal{D}[0,T],\mathcal{H})$ given by \eqref{solpb} and satisfying 
\begin{eqnarray*}
%\label{boundsol}
\mathbb{E}\left[\sup_{0\leq t\leq T}\left\|X(t)\right\|^2\right]<\infty,
\end{eqnarray*} 
where $(\mathcal{D}[0,T],\mathcal{H})$ denotes the space of all c\`adl\`ag processes  defined on $[0;T]$ with values in $\mathcal{H}$.
\end{thm}
\textit{Proof}
\cite[Theorem 6.1]{Dio} gives the result of existence and uniqueness of the mild solution $X$ for SPDEs of type \eqref{SPDE} only driven by fractional Brownian motion. As in the work done in \cite{Alb}, we can extend it to \eqref{SPDE} by combining with \cite{Dio}.

\section{Numerical schemes and Main results}
\label{numschemes}
For the seek of simplicity, we assume the family of linear operators $A(t)$\footnote{ Indeed the operators $A(t)$ are identified to their $L^2$ realizations
given in \eqref{family} (see \cite{Fujita}).} to be of second order and has the following form
\begin{eqnarray}
\label{family}
A(t)u=\sum_{i,j=1}^d\frac{\partial}{\partial x_i}\left(q_{ij}(x,t)\frac{\partial u}{\partial x_j}\right)-\sum_{j=1}^dq_j(x,t)\frac{\partial u}{\partial x_j}.
\end{eqnarray}
We require the coefficients $q_{i,j}$ and $q_j$ to be smooth functions of the variable $x\in\overline{\Lambda}$ and H\"{o}lder-continuous with respect to $t\in[0,T]$. 
We further assume that there exists a positive constant $c$ such that the following  ellipticity condition holds
\begin{eqnarray}
\label{ellip}
\sum_{i,j=1}^dq_{ij}(x,t)\xi_i\xi_j\geq c\vert \xi\vert^2, \quad (x,t)\in\overline{\Lambda}\times [0,T].
\end{eqnarray}
It is well known that the family of linear operators defined by \eqref{family} satisfy \assref{linear}, see e.g. \cite[Section 7.6]{Pazy}, \cite[Chapter III, Section 12]{Fujita} or \cite{Mukmab}.

We consider the discretization of the spatial domain by a finite element triangulation \cite{Mukb,Wanb}.
Let $\mathcal{T}_h$ be a set of disjoint intervals of $\Lambda$ (for $d=1$), a triangulation of $\Lambda$ (for $d=2$) or a set of tetrahedra (for $d=3$) 
with maximal length $h$ satisfying the usual regularity assumptions.\\
Let $V_h\subset \mathcal{H}$ denote the space of continuous functions that are piecewise linear over triangulation $\mathcal{T}_h$. Next, we introduce the projection $P_h$ from $L^2(\Lambda)$ to $V_h$ define  by
\begin{eqnarray*}
%\label{proj}
\left\langle P_hu,\mathcal{E}\right\rangle =\left\langle u,\mathcal{E}\right\rangle,\quad u\in L^2(\Lambda),\quad \mathcal{E}\in V_h.
\end{eqnarray*} 
The discrete operator $A_h(t):V_h\rightarrow V_h$ is defined by
\begin{eqnarray*}
%\label{discrete_op}
\left\langle A_h(t)\varphi,\mathcal{X}\right\rangle =-a(t)\left( \varphi, \mathcal{X}\right) ,\quad \varphi,\mathcal{X}\in V_h,
\end{eqnarray*}
where $a(t)$ is the corresponding bilinear form of $-A(t)$.

Like the operator $A(t)$, the discrete operator $A_h(t)$ is also the generator of an analytic semigroup on $L^2(\Lambda)$ denoted by  $S_t^h(s)=:e^{-sA_h(t)}$  and  \propref{semigroup} holds for $A_h(t)$ uniformly on $h>0$ and $t\in[0,T]$,  i.e. for all $\alpha\geq 0$ and $0\leq \gamma\leq 1$, there exists a positive constant $C'_1$ such that the following estimates hold uniformly on $h>0$ and $t\in[0,T]$, see \cite{Mukmab}.  
\begin{eqnarray*}
\left\|(A_h(t))^{\alpha}e^{-sA_h(t)}\right\|_{L(\mathcal{H})}\leq C'_1 s^{-\alpha},\; s>0,\label{semigroup3}\;
\left\|(A_h(t))^{-\gamma}\left(I-e^{-rA_h(t)}\right)\right\|_{L(\mathcal{H})}\leq C'_1 r^{\gamma},\quad r\geq0. %\label{semigroup4} 
\end{eqnarray*}
The semi-discrete space version of problem \eqref{SPDE} is to find $X^h(t)=X^h(\cdot,t)$, $t\in[0,T]$, such that
\begin{eqnarray}
\label{disSPDE}
\left\{
  \begin{array}{ll}
    dX^h(t)+A_h(t)X^h(t)dt=P_hF(t,X^h(t))dt+P_h\phi(t) dB^H(t)+\int_{\chi}P_h\psi(z)\tilde{N}(dz,dt), \\
    X^h(0)=P_hX_0.
  \end{array}
\right.
\end{eqnarray}
Now applying the stochastic exponential integrator scheme
\cite{Lor,Muka,Nou} to \eqref{disSPDE}, we obtain the following scheme called Stochastic Magnus-type Integrator \cite{Mukmaa,Tammab} (\textbf{SMTI})
\begin{eqnarray}
\label{SMTI}
\left\{
  \begin{array}{lll}
    X^h_0=P_hX_0, \\
    X^h_{m+1}=e^{-\Delta t A_{h,m}}X^h_m+\Delta t \varphi_1(\Delta t A_{h,m})P_hF(t_m,X^h_m)+e^{-\Delta t A_{h,m}}P_h\phi(t_m) \Delta B^H_m\\
    \quad\quad \quad+\int_{\chi}e^{-\Delta t A_{h,m}}P_h\psi(z)\tilde{N}(dz,\Delta t),\quad m\in\{0,1,\cdot\cdot\cdot,M-1\}.
  \end{array}
\right.
\end{eqnarray}
Furthermore applying the linear implicit Euler method \cite{Kam,Wanb,Mukb,Yan,Nou} to \eqref{disSPDE}, yields the following fully discrete scheme (see \cite{Mukmab})
\begin{eqnarray}
\label{impl}
\left\{
  \begin{array}{lll}
    Y^h_0=P_hX_0, \\
    Y^h_{m+1}=S^m_{h,\Delta t}Y^h_m+\Delta t S^m_{h,\Delta t}P_hF(t_m,Y^h_m)+S^m_{h,\Delta t}P_h\phi(t_m) \Delta B^H_m\\
    \quad\quad\quad+\int_{\chi}S^m_{h,\Delta t}P_h\psi(z)\tilde{N}(dz,\Delta t),\quad m\in\{0,1,\cdot\cdot\cdot,M-1\},
  \end{array}
\right.
\end{eqnarray} 
where $t_m=m\Delta t$, with $\Delta t= \frac TM$ for a given $M\in\N$, 
$$\Delta B^H_m:=B^H_{m+1}-B^H_m=\sum_{i\in\N^d}\sqrt{q_i}\left(\beta_i^H(t_{m+1})-\beta_i^H(t_m)\right)e_i,$$
$$ A_{h,m}:=A_h(t_m),\; S^m_{h,\Delta t}:=(I+\Delta t A_{h,m})^{-1}$$  and the linear operator $\varphi(\Delta t A_{h,m})$ is given by 
\begin{eqnarray*}
\varphi(\Delta t A_{h,m}):=\frac 1{\Delta t}\int_0^{\Delta t}e^{-(\Delta t-s)A_{h,m}}ds. \nonumber
\end{eqnarray*}
With the above two schemes \eqref{SMTI} and \eqref{impl} in hand, we state our main results in the following theorem.
\begin{thm}
\label{mainthm}
Let $X(t_m)$ be the mild solution of \eqref{SPDE} at time $t_m=m\Delta t$, given by \eqref{solpb}. Let $X^h_m$ be the numerical approximations through \eqref{SMTI} and $Y^h_m$ be the numerical approximations through \eqref{impl}. 
Under Assumptions \ref{linear}-\ref{derivate}, the following  estimations holds
\begin{eqnarray*}
%\label{mainthmSMTI}
\left(\mathbb{E}\|X(t_m)-X^h_m\|^2\right)^{\frac 12}\leq C\left(h^{2H+\beta-1-\epsilon}+\Delta t^{\frac{2H+\beta-1}2}\right),
\end{eqnarray*}
and 
\begin{eqnarray*}
%\label{mainthmimpl}
\left(\mathbb{E}\|X(t_m)-Y^h_m\|^2\right)^{\frac 12}\leq C \left(h^{2H+\beta-1-\epsilon}+\Delta t^{\frac{2H+\beta-1}2-\epsilon}\right),
\end{eqnarray*}
where $\epsilon$ is a positive number small enough.  
\end{thm}
To prove  the main results, we will split  the full error in two terms, namely the spatial error and the temporal error and estimate both errors separately. 
We  will start with the spatial error. 
\begin{remark}
Note that for SPDE driven only with fractional Brownian motion ( $\psi=0$ in \eqref{SPDE}), the same orders of convergence are achieved in  \thmref{mainthm} without Assumption \ref{derivate} 
\footnote{That means without information on the Frechet derivative $F'$}.
\end{remark}
\section{Regularity estimates and space error}
\label{conv}
Throughout the next of this paper, we take $C$ as a generic constant that may change from one place to another.
\subsection{Preliminaries and regularity estimates of the semi-discrete problem}
\label{convspa}
To prove our main result, we need some preparatory results.
\begin{lem}(\cite{Tamma})
\label{equi} Let \assref{linear} be fulfilled
\begin{enumerate}
\item [(i)] For any $\gamma\in[0,1]$, the following equivalence norms holds
\begin{eqnarray}
\label{equi1}
C^{-1}\left\|(A_h(0))^{-\gamma}v\right\|\leq\left\|(A_h(t))^{-\gamma}v\right\|\leq C\left\|(A_h(0))^{-\gamma}v\right\|,\quad v\in V_h,\quad t\in[0;T].
\end{eqnarray}
\item [(ii)] For any $\gamma\in[0,1]$, the following equivalence norms holds
\begin{eqnarray}
\label{equi2}
C^{-1}\left\|(A_h(0))^{\gamma}v\right\|\leq\left\|(A_h(t))^{\gamma}v\right\|\leq C\left\|(A_h(0))^{\gamma}v\right\|,\quad v\in V_h,\quad t\in[0;T].
\end{eqnarray}
\item [(iii)] For any $\alpha\in[-\frac 12,\frac 12]$, it holds that (see \cite[Lemma 2.2]{Mukmaa}) 
\begin{eqnarray}
\label{equi3}
\left\|(A_{h,k})^{\alpha}P_hv\right\|\leq C\left\|(A_{h,l})^{\alpha}v\right\|\leq C\left\|(A(0))^{\alpha}v\right\|,\quad v\in V_h,\quad 0\leq k,l\leq M-1.
\end{eqnarray}
\item [(iv)] The following estimates hold
\begin{eqnarray*}
&&\left\|\left(A_h(t)-A_h(s)\right)\left(A_h(r)\right)^{-1} u^h\right\|\leq C|t-s|\left\|u^h\right\|,\quad r,s,t\in[0;T],\quad u^h\in V_h.\nonumber\\
&&\left\|\left(A_h(r)\right)^{-1}\left(A_h(s)-A_h(t)\right)u^h \right\|\leq C|s-t|\left\|u^h\right\|,\quad r,s,t\in[0;T], \quad u^h\in V_h\cap\mathcal{D}.\nonumber
\end{eqnarray*}
\end{enumerate}
\end{lem}
\begin{prp}(\cite[Remark 3.1]{Mukmab})
\label{disevol} There exists a unique evolution system $U_h:\Delta (T)\rightarrow L(\mathcal{H})$, satisfying
\begin{eqnarray}
\label{disevol1}
U_h(t,s)=S^h_s(t,s)+\int_s^tS^h_{\tau}(t-\tau)R^h(t,s)d\tau,
\end{eqnarray}
\end{prp}
where $R^h(t,s):=\sum_{m=1}^{\infty}R^h_m(t,s)$ with $R^h_m(t,s)$ given by \cite[(6.3),P149]{Pazy}
\begin{eqnarray*}
R^h_1(t,s):=\left(A_h(s)-A_h(t)\right)S^h_s(t-s),\quad R^h_{m+1}(t,s):=\int_{s}^tR^h_1(t,s)R^h_m(\tau,s)d\tau,\quad m\geq 1.
\end{eqnarray*}
Note also that from \cite[(6.6)]{Pazy}, the following identity holds
\begin{eqnarray*}
\label{disevol2}
R^h(t,s)=R^h_1(t,s)+\int_{s}^tR^h_1(t,\tau)R^h(\tau,s)d\tau.
\end{eqnarray*}
The mild solution of the semi-discrete problem \eqref{disSPDE} can therefore be represented as follows
\begin{eqnarray}
\label{soldispb}
X^h(t)&=&U_h(t,0)P_hX_0+\int_0^tU_h(t,s)P_hF(s,X^h(s))ds+\int_0^tU_h(t,s)P_h\phi(s)dB^H(s)\nonumber\\
&+&\int_0^t\int_{\chi}U_h(t,s)P_h\psi(z)\tilde{N}(dz,ds).
\end{eqnarray}
\begin{prp}(\cite[Lemma 3.2]{Mukmab})
\label{disevol'} Under \assref{linear}, the evolution system $U_h:\Delta(T)\rightarrow L(\mathcal{H})$ satisfies:  
\begin{enumerate}
\item [(i)] $U_h(\cdot,s)\in\mathcal{C}^1( ]s,T];L(\mathcal{H}))$, $0\leq s\leq T$ and
\begin{eqnarray*}
\frac{\partial U_h}{\partial t}(t,s)=-A_h(t)U_h(t,s),\quad\left\|A_h(t)U_h(t,s)\right\|_{L(\mathcal{H})}\leq \frac{C}{t-s},\quad 0\leq s<t\leq T.\nonumber
\end{eqnarray*}
\item [(iii)] $U_h(t,\cdot)v\in\mathcal{C}^1([0,t[;\mathcal{H})$, $0<t\leq T$, $v\in\mathcal{D}(A_h(0))$ and 
\begin{eqnarray*}
\frac{\partial U_h}{\partial s}(t,s)v=U_h(t,s)A_h(s)v,\quad \left\|A_h(t)U_h(t,s)(A_h(s))^{-1}\right\|_{L(\mathcal{H})}\leq C\quad 0\leq s\leq t\leq T.\nonumber
\end{eqnarray*}
\item [(iii)] For any $(t,r),(r,s)\in\Delta(T)$ it holds that
\begin{eqnarray*}
U_h(s,s)=I\quad\text{and}\quad U_h(t,r)U_h(r,s)=U_h(t,s).
\end{eqnarray*}
\end{enumerate}
\end{prp}
\begin{lem}
\label{propdisevol} Let \assref{linear} be fulfilled
\begin{enumerate}
\item [(i)] The following estimates hold
\begin{eqnarray}
\label{propdisevol1}
\left\|R^h(t,s)\right\|_{L(\mathcal{H})}\leq C,\quad \left\|U_h(t,s)\right\|_{L(\mathcal{H})}\leq C,\quad 0\leq s\leq t\leq T.
\end{eqnarray}
\item [(ii)] For any $0\leq\alpha,\gamma\leq 1$ and  $0\leq s\leq t\leq T$ the following estimates hold
\begin{eqnarray}
\left\|(A_h(r))^{\alpha}U_h(t,s)\right\|_{L(\mathcal{H})}\leq C(t-s)^{-\alpha},\quad r\in[0,T].\label{propdisevol2}\\
\left\|U_h(t,s)(A_h(s))^{\alpha}\right\|_{L(\mathcal{H})}\leq C(t-s)^{-\alpha},\quad r\in[0,T].\label{propdisevol3}\\
\left\|(A_h(r))^{\alpha}U_h(t,s)(A_h(s))^{-\gamma}\right\|_{L(\mathcal{H})}\leq C(t-s)^{\gamma-\alpha},\quad r\in[0,T].\label{propdisevol4}\\
\left\|(A_h(r))^{\alpha}U_h(t,s)(A_h(s))^{\gamma}\right\|_{L(\mathcal{H})}\leq C(t-s)^{-(\gamma+\alpha)},\quad r\in[0,T].\label{propdisevol5}
\end{eqnarray}
\item [(iii)] For any $0\leq\gamma\leq 1$ and $0\leq s\leq t\leq T$, the following estimates hold
\begin{eqnarray*}
\left\|(U_h(t,s)-I)(A_h(s))^{-\gamma}\right\|_{L(\mathcal{H})}\leq C(t-s)^{\gamma},\label{propdisevol6}\\
\left\|(A_h(r))^{-\gamma}(U_h(t,s)-I)\right\|_{L(\mathcal{H})}\leq C(t-s)^{\gamma},\quad r\in[0,T]. %\label{propdisevol7}
\end{eqnarray*}
\end{enumerate}
\end{lem}
\textit{Proof:} The proof of all results in this lemma except \eqref{propdisevol5} can be found in \cite[Lemma 5]{Tamma}. For the proof of \eqref{propdisevol5},  we just have to decompose $U_h(t,s)=U_h(t,\frac{t+s}2)U_h(\frac{t+s}2,s)$, then applying \eqref{propdisevol2} and \eqref{propdisevol3}, we get our result.$\hfill\square$

The following lemma is an important result which plays a crucial role to obtain optimal regularity results,  useful in this work.
\begin{lem}
\label{inte}
\begin{enumerate}
\item [(i)] For any $0\leq \rho,\gamma,\eta\leq 1$, $0\leq \tau_1\leq\tau_2\leq T$, $\tau\in[0,T]$, then the following estimates hold
\begin{eqnarray}
&&\int_{\tau_1}^{\tau_2}\left\|U_h(\tau_2,r)(A_h(r))^{\frac{\rho}2}\right\|^2_{L(\mathcal{H})}dr\leq C(\tau_2-\tau_1)^{1-\rho}.\label{inte1}\\
&&\int_{\tau_1}^{\tau_2}\left\|(A_h(\tau))^{\frac{\rho}2}U_h(\tau_2,r)(A_h(r))^{-\frac{\eta}2}\right\|^2_{L(\mathcal{H})}dr\leq C(\tau_2-\tau_1)^{1-\rho+\eta},\; 0\leq \rho-\eta\leq 1.\label{inte2}\\
&&\int_{\tau_1}^{\tau_2}\left\|(A_h(\tau))^{\frac{\rho}2}U_h(\tau_2,r)(A_h(r))^{\frac{\gamma}2}\right\|^2_{L(\mathcal{H})}dr\leq C(\tau_2-\tau_1)^{1-(\rho+\gamma)},\; 0\leq \rho+\gamma\leq 1.\label{inte3}
\end{eqnarray}
\item [(ii)] For any $0\leq \alpha,\kappa\leq 1$  such that $0\leq \alpha+\kappa\leq 2H$, with $H\in(\frac 12,1)$,  the following estimate holds
\begin{eqnarray}
\label{inte4}
\int_{\tau_1}^{\tau_2}\left\|(A_h(\tau))^{\frac{\alpha}2}U_h(\tau_2,r)(A_h(r))^{\frac{\kappa}2}\right\|^{\frac 1H}_{L(\mathcal{H})}dr\leq C(\tau_2-\tau_1)^{\frac{2H-(\alpha+\kappa)}{2H}},\quad \tau\in[0,T].
\end{eqnarray}
\end{enumerate}
\end{lem}
\textit{Proof:} See \cite[Lemma 2.7]{Mukmaa} for the proof of \eqref{inte1}. Concerning the proof of \eqref{inte2}, using \eqref{disevol1}, triangle inequality, the estimate $(a+b)^2\leq 2a^2+2b^2$ we have
\begin{eqnarray}
\label{inte2'}
&&\int_{\tau_1}^{\tau_2}\left\|(A_h(\tau))^{\frac{\rho}2}U_h(\tau_2,r)(A_h(r))^{-\frac{\eta}2}\right\|^2_{L(\mathcal{H})}dr\nonumber\\
&\leq& 2\int_{\tau_1}^{\tau_2}\left\|(A_h(\tau))^{\frac{\rho}2}S^h_r(\tau_2-r)(A_h(r))^{-\frac{\eta}2}\right\|^2_{L(\mathcal{H})}dr\nonumber\\
&+&2\int_{\tau_1}^{\tau_2}\left\|\int_{r}^{\tau_2}(A_h(\tau))^{\frac{\rho}2}S^h_{\sigma}(\tau_2-\sigma)R^h(\sigma,r)(A_h(r))^{-\frac{\eta}2}d\sigma\right\|^2_{L(\mathcal{H})}dr\nonumber\\
&=:&2J_1+2J_2.
\end{eqnarray}
Thanks to \eqref{equi2}, stability properties of the discrete semigroup $S^h_r(t)$ and \cite[(71)]{Mukmaa}. This yields
\begin{eqnarray}
\label{inte21}
J_1&=&\int_{\tau_1}^{\tau_2}\left\|(A_h(\tau))^{\frac{\rho}2}S^h_r(\tau_2-r)(A_h(r))^{-\frac{\eta}2}\right\|^2_{L(\mathcal{H})}dr\nonumber\\
&\leq& C\int_{\tau_1}^{\tau_2}\left\|(A_h(r))^{\frac{\rho}2}S^h_r(\tau_2-r)(A_h(r))^{-\frac{\eta}2}\right\|^2_{L(\mathcal{H})}dr\nonumber\\
&\leq& C\int_{\tau_1}^{\tau_2}\left\|(A_h(\tau))^{\frac{\rho-\eta}2}S^h_r(\tau_2-r)\right\|^2_{L(\mathcal{H})}dr\\
&\leq& C(\tau_2-\tau_1)^{1-\rho+\eta}.
\end{eqnarray}
Triangle inequality, H\"older inequality, \eqref{propdisevol1}, stability properties of discrete semigroup $S^h_r(t)$, \cite[(71)]{Mukmaa} yield
\begin{eqnarray}
\label{inte22}
J_2&=&\int_{\tau_1}^{\tau_2}\left\|\int_{r}^{\tau_2}(A_h(\tau))^{\frac{\rho}2}S^h_{\sigma}(\tau_2-\sigma)R^h(\sigma,r)(A_h(r))^{-\frac{\eta}2}d\sigma\right\|^2_{L(\mathcal{H})}dr\nonumber\\
&\leq&\int_{\tau_1}^{\tau_2}(\tau_2-\sigma)\int_{r}^{\tau_2}\left\|(A_h(\tau))^{\frac{\rho}2}S^h_{\sigma}(\tau_2-r)\right\|^2_{L(\mathcal{H})}\left\|R^h(\sigma,r)(A_h(r))^{-\frac{\eta}2}\right\|^2_{L(\mathcal{H})}d\sigma dr\nonumber\\
&\leq&C\int_{\tau_1}^{\tau_2}(\tau_2-r)\left(\int_{r}^{\tau_2}\left\|(A_h(\tau))^{\frac{\rho}2}S^h_{\sigma}(\tau_2-\sigma))\right\|^2_{L(\mathcal{H})}d\sigma\right) dr\nonumber\\
&\leq&C\int_{\tau_1}^{\tau_2}(\tau_2-r)(\tau_2-r)^{1-\rho}dr\nonumber\\
&\leq&C(\tau_2-\tau_1)^{3-\rho}.
\end{eqnarray}
Hence substituting \eqref{inte21} and \eqref{inte22} in \eqref{inte2'} gives
\begin{eqnarray*}
\int_{\tau_1}^{\tau_2}\left\|(A_h(\tau))^{\frac{\rho}2}U_h(\tau_2,r)(A_h(r))^{-\frac{\eta}2}\right\|^2_{L(\mathcal{H})}dr&\leq& C(\tau_2-\tau_1)^{1-\rho+\eta}+C(\tau_2-\tau_1)^{3-\rho}\nonumber\\
&\leq& C(\tau_2-\tau_1)^{1-\rho+\eta}.\nonumber
\end{eqnarray*}
The proof of \eqref{inte2} is thus completed. For the proof of \eqref{inte3}, as in \eqref{inte2'} we also split it in two following terms
\begin{eqnarray}
\label{inte3'}
&&\int_{\tau_1}^{\tau_2}\left\|(A_h(\tau))^{\frac{\gamma}2}U_h(\tau_2,r)(A_h(r))^{\frac{\rho}2}\right\|^2_{L(\mathcal{H})}dr\nonumber\\
&\leq& 2\int_{\tau_1}^{\tau_2}\left\|(A_h(\tau))^{\frac{\gamma}2}S^h_r(\tau_2-r)(A_h(r))^{\frac{\rho}2}\right\|^2_{L(\mathcal{H})}dr\nonumber\\
&+&2\int_{\tau_1}^{\tau_2}\left\|\int_{r}^{\tau_2}(A_h(\tau))^{\frac{\gamma}2}S^h_{\sigma}(\tau_2-\sigma)R^h(\sigma,r)(A_h(r))^{\frac{\rho}2}d\sigma\right\|^2_{L(\mathcal{H})}dr\nonumber\\
&=:&2J_3+2J_4.
\end{eqnarray}
As in \eqref{inte21} we easily get
\begin{eqnarray}
\label{inte31}
J_3\leq C(\tau_2-\tau_1)^{1-(\gamma+\rho)}.
\end{eqnarray}
Using triangle inequality, \lemref{propdisevol}, \eqref{semigroup1}, \cite[(3.60)]{Tammab} and the change of variable $\sigma=r(1-u)+\tau_2u$ we obtain
\begin{eqnarray}
\label{inte32}
J_4&=&\int_{\tau_1}^{\tau_2}\left\|\int_{r}^{\tau_2}(A_h(\tau))^{\frac{\gamma}2}S^h_{\sigma}(\tau_2-\sigma)R^h(\sigma,r)(A_h(r))^{\frac{\rho}2}d\sigma\right\|^2_{L(\mathcal{H})}dr\nonumber\\
&\leq&\int_{\tau_1}^{\tau_2}\left(\int_{r}^{\tau_2}\left\|(A_h(\tau))^{\frac{\rho}2}S^h_{\sigma}(\tau_2-\sigma)\right\|_{L(\mathcal{H})}\left\|R^h(\sigma,r)(A_h(r))^{-\frac{\eta}2}\right\|_{L(\mathcal{H})}d\sigma\right)^2 dr\nonumber\\
&\leq&\int_{\tau_1}^{\tau_2}\left(\int_{r}^{\tau_2}(\tau_2-\sigma)^{-\frac{\gamma}2}(\sigma-r)^{-\frac{\eta}2}d\sigma\right)^2 dr\nonumber\\
&\leq&C\int_{\tau_1}^{\tau_2}(\tau_2-r)^{2-(\gamma+\rho)}\left(\int_0^1(1-u)^{-\frac{\gamma}2}u^{-\frac{\rho}2}du\right)^2 dr\nonumber\\
&\leq&C\mathcal{B}^2\left(1-\frac{\gamma}2,1-\frac{\rho}2\right)\int_{\tau_1}^{\tau_2}(\tau_2-r)^{2-(\gamma+\rho)} dr\nonumber\\
&\leq&C(\tau_2-\tau_1)^{3-(\gamma+\rho)},
\end{eqnarray}
where $\mathcal{B}(x,y)$ is the Beta function. Hence putting \eqref{inte31} and \eqref{inte32} in \eqref{inte3'} yields
\begin{eqnarray*}
\int_{\tau_1}^{\tau_2}\left\|(A_h(\tau))^{\frac{\gamma}2}U_h(\tau_2,r)(A_h(r))^{\frac{\rho}2}\right\|^2_{L(\mathcal{H})}dr&\leq& C(\tau_2-\tau_1)^{1-(\gamma+\rho)}+C(\tau_2-\tau_1)^{3-(\gamma+\rho)}\nonumber\\
&\leq& C(\tau_2-\tau_1)^{1-(\gamma+\rho)}.\nonumber
\end{eqnarray*}
The proof of \eqref{inte3} is thus completed. Let us prove (ii), we begin by considering the case $0\leq \alpha+\kappa<2H$ and using \eqref{propdisevol5} we have
\begin{eqnarray*}
\label{inte41}
\int_{\tau_1}^{\tau_2}\left\|(A_h(\tau))^{\frac{\alpha}2}U_h(\tau_2,r)(A_h(r))^{\frac{\kappa}2}\right\|^{\frac 1H}_{L(\mathcal{H})}dr
&\leq& C\int_{\tau_1}^{\tau_2}(\tau_2-r)^{-\frac{\alpha+\kappa}{2H}}dr\nonumber\\
&\leq& C(\tau_2-\tau_1)^{1-\frac{\alpha+\kappa}{2H}}\nonumber\\
&=& C(\tau_2-\tau_1)^{\frac{2H-(\alpha+\kappa)}{2H}}.
\end{eqnarray*} 
Now we move on the critical case $\alpha+\kappa=2H$ with $H\in(\frac 12,1)$.
\begin{itemize}
\item For $H=\frac 12$, we have $\alpha+\kappa=1$ and using \eqref{inte3} with $\alpha=\gamma$ and $\kappa=\rho$, we obtain
\begin{eqnarray*}
\int_{\tau_1}^{\tau_2}\left\|(A_h(\tau))^{\frac{\alpha}2}U_h(\tau_2,r)(A_h(r))^{\frac{\kappa}2}\right\|^{\frac 1H}_{L(\mathcal{H})}dr
&=&\int_{\tau_1}^{\tau_2}\left\|(A_h(\tau))^{\frac{\alpha}2}U_h(\tau_2,r)(A_h(r))^{\frac{\kappa}2}\right\|^2_{L(\mathcal{H})}dr\nonumber\\
&\leq& C(\tau_2-\tau_1)^{1-(\alpha+\kappa)}\nonumber\\
&=& C(\tau_2-\tau_1)^{\frac{2H-(\alpha+\kappa)}{2H}}.\nonumber
\end{eqnarray*} 
\item 
For $H=1$, Using the definition of $U_h(\tau_2, r)$ and triangle inequality yields 
\begin{eqnarray*}
&&\int_{\tau_1}^{\tau_2}\left\|(A_h(\tau))^{\frac{\alpha}2}U_h(\tau_2,r)(A_h(r))^{\frac{\kappa}2}\right\|^{\frac 1H}_{L(\mathcal{H})}dr\nonumber\\
&=&\int_{\tau_1}^{\tau_2}\left\|(A_h(\tau))^{\frac{\alpha}2}U_h(\tau_2,r)(A_h(r))^{\frac{\kappa}2}\right\|_{L(\mathcal{H})}dr\nonumber\\
&\leq& \int_{\tau_1}^{\tau_2}\left\|(A_h(\tau))^{\frac{\alpha}2}S^h_r(\tau_2-r)(A_h(r))^{\frac{\kappa}2}\right\|^2_{L(\mathcal{H})}dr\nonumber\\
&+&\int_{\tau_1}^{\tau_2}\left\|\int_{r}^{\tau_2}(A_h(\tau))^{\frac{\alpha}2}S^h_{\sigma}(\tau_2-\sigma)R^h(\sigma,r)(A_h(r))^{\frac{\kappa}2}d\sigma\right\|_{L(\mathcal{H})}dr.\nonumber
\end{eqnarray*}
Note that $\alpha+\kappa=2$, using \lemref{propdisevol}, \eqref{semigroup1}, \cite[(94)]{Mukmaa} with $\rho=\frac{\alpha+\kappa}2$ and a similar reasoning as in $J_4$ yields
\begin{eqnarray*}
&&\int_{\tau_1}^{\tau_2}\left\|(A_h(\tau))^{\frac{\alpha}2}U_h(\tau_2,r)(A_h(r))^{\frac{\kappa}2}\right\|^{\frac 1H}_{L(\mathcal{H})}dr\nonumber\\
&\leq& C\int_{\tau_1}^{\tau_2}\left\|(A_h(r))^{\frac{\alpha}2}S^h_r(\tau_2-r)(A_h(r))^{\frac{\kappa}2}\right\|_{L(\mathcal{H})}dr\nonumber\\
&+&\int_{\tau_1}^{\tau_2}\left(\int_{r}^{\tau_2}\left\|(A_h(\tau))^{\frac{\alpha}2}S^h_{\sigma}(\tau_2-\sigma)\right\|_{L(\mathcal{H})}\left\|R^h(\sigma,r)(A_h(r))^{\frac{\kappa}2}\right\|_{L(\mathcal{H})}d\sigma\right) dr\nonumber\\
&\leq& C\int_{\tau_1}^{\tau_2}\left\|(A_h(r))^{\frac{\alpha+\kappa}2}S^h_r(\tau_2-r)\right\|_{L(\mathcal{H})}dr+C\int_{\tau_1}^{\tau_2}\left(\int_{r}^{\tau_2}(\tau_2-\sigma)^{-\frac{\alpha}2}(\sigma-r)^{-\frac{\kappa}2}d\sigma\right) dr\nonumber\\
&\leq& C(\tau_2-\tau_1)^{1-\frac{\alpha+\kappa}2}+C(\tau_2-\tau_1)^{2-\frac{\alpha+\kappa}2}\nonumber\\
&\leq& C(\tau_2-\tau_1)^{\frac{2-(\alpha+\kappa)}2}= C(\tau_2-\tau_1)^{\frac{2H-(\alpha+\kappa)}{2H}}.\nonumber
\end{eqnarray*}
\end{itemize}
Hence the proof of \eqref{inte4} is thus completed by interpolation theory. $\hfill\square$.

The following lemma proving  spatial and temporal regularity for the solution process \eqref{soldispb} will be useful in our convergence analysis.
\begin{lem}
\label{reg}
Suppose  that Assumptions \ref{linear}--\ref{initdata} hold, then the unique mild solution $X^h(t)$ of \eqref{disSPDE} satisfied
\begin{eqnarray} 
\label{bound1}
\left\|X^h(t)\right\|_{L^2(\Omega,\mathcal{H})}\leq C;\quad \left\|F(t,X^h(t))\right\|_{L^2(\Omega,\mathcal{H})}\leq C,\quad 0\leq t\leq T.
\end{eqnarray}
Moreover the following optimal regularity in space  holds
\begin{eqnarray}
\label{spareg}
\left\|(A_h(\tau))^{\frac{2H+\beta-1}2}X^h(t)\right\|_{L^2(\Omega,\mathcal{H})}\leq C\left(1+\left\|(A(0))^{\frac{2H+\beta-1}2}X_0\right\|_{L^2(\Omega,\mathcal{H})}\right),\quad 0\leq t,\tau\leq T.
\end{eqnarray}
The following optimal time regularity also holds
\begin{eqnarray}
\label{timreg}
\left\|X^h(t_2)-X^h(t_1)\right\|_{L^2(\Omega,\mathcal{H})}\leq C(t_2-t_1)^{\frac{\min(2H+\beta-1,1)}2}\left(1+\left\|(A(0))^{\frac{2H+\beta-1}2}X_0\right\|_{L^2(\Omega,\mathcal{H})}\right), 
\end{eqnarray}
for any $0\leq t_1\leq t_2\leq T$.
\end{lem}
\textit{Proof:} We first show that
\begin{eqnarray}
\label{bound}
\left\|X^h(t)\right\|_{L^2(\Omega,\mathcal{H})}\leq C,\quad 0\leq t\leq T.
\end{eqnarray}
By the definition of mild solution we have
\begin{eqnarray*}
X^h(t)&=&U_h(t,0)P_hX_0+\int_0^tU_h(t,s)P_hF(s,X^h(s))ds+\int_0^tU_h(t,s)P_h\phi(s)dB^H(s)\nonumber\\
&+&\int_0^t\int_{\chi}U_h(t,s)P_h\psi(z)\tilde{N}(dz,ds).\nonumber
\end{eqnarray*} 
Taking the norm in both sides and using triangle inequality yields
\begin{eqnarray*}
\label{bound*}
\left\|X^h(t)\right\|^2_{L^2(\Omega,\mathcal{H})}&\leq&4\left\|U_h(t,0)P_hX_0\right\|^2_{L^2(\Omega,\mathcal{H})}+4\left\|\int_0^tU_h(t,s)P_hF(s,X^h(s))ds\right\|^2_{L^2(\Omega,\mathcal{H})}\nonumber\\
&+&4\left\|\int_0^tU_h(t,s)P_h\phi(s)dB^H(s)\right\|^2_{L^2(\Omega,\mathcal{H})}\nonumber\\
&+&4\left\|\int_0^t\int_{\chi}U_h(t,s)P_h\psi(z)\tilde{N}(dz,ds)\right\|^2_{L^2(\Omega,\mathcal{H})}\nonumber\\
&=:& 4\sum_{i=1}^4 I_i.
\end{eqnarray*}
In what follows, we bound $I_i, i=1,2,3,4$ one by one. First, thanks to \eqref{propdisevol1} and the boundedness of $P_h$,  we easily have
\begin{eqnarray}
\label{bound*1}
I_1=\left\|U_h(t,0)P_hX_0\right\|^2_{L^2(\Omega,\mathcal{H})}\leq C\left\|X_0\right\|^2_{L^2(\Omega,\mathcal{H})}<C.
\end{eqnarray}
For $I_2$, using  H\"older inequality, \eqref{propdisevol1}, the boundedness of $P_h$ and \assref{nonlinear} yields
\begin{eqnarray}
%\label{bound*2}
I_2&=& \left\|\int_0^tU_h(t,s)P_hF(s,X^h(s))ds\right\|^2_{L^2(\Omega,\mathcal{H})}\nonumber\\
&\leq& C\int_0^t\left\|U_h(t,s)P_h\left(F(s,0)+F(s,X^h(s))-F(s,0)\right)\right\|^2_{L^2(\Omega,\mathcal{H})}ds\nonumber\\
&\leq& C\int_0^t\left\|F(s,0)\right\|^2_{L^2(\Omega,\mathcal{H})}+\left\|F(s,X^h(s))-F(s,0)\right\|^2_{L^2(\Omega,\mathcal{H})}ds\nonumber\\
&\leq& C\int_0^t \left(1+\left\|X^h(s)\right\|^2_{L^2(\Omega,\mathcal{H})}\right)ds\leq C+C\int_0^t\left\|X^h(s)\right\|^2_{L^2(\Omega,\mathcal{H})}ds.
\end{eqnarray}
For $I_3$, \eqref{fracint2}, inserting an appropriate power of $A_h(s)$, using \eqref{equi1} with $\alpha=\frac{1-\beta}2$, \eqref{equi3}, \assref{noiseandjump}(more precisely \eqref{noise2}) and \eqref{inte1} with $\rho=1-\beta$ yield
\begin{eqnarray}
%\label{bound*3}
I_3&=&\left\|\int_0^tU_h(t,s)P_h\phi(s)dB^H(s)\right\|^2_{L^2(\Omega,\mathcal{H})}\nonumber\\
&\leq& C\int_0^t\left\|U_h(t,s)P_h\phi(s)\right\|^2_{L^0_2}ds\nonumber\\
&\leq& C\int_0^t\left\|U_h(t,s)(A_h(s))^{\frac{1-\beta}2}\right\|^2_{L(\mathcal{H})}\left\|(A_h(s))^{\frac{\beta-1}2}P_h\phi(s)\right\|^2_{L^0_2}ds\nonumber\\
&\leq& C\int_0^t\left\|U_h(t,s)(A_h(s))^{\frac{1-\beta}2}\right\|^2_{L(\mathcal{H})}\left\|(A(0))^{\frac{\beta-1}2}\phi(s)\right\|^2_{L^0_2}ds\nonumber\\
&\leq& C\int_0^t\left\|U_h(t,s)(A_h(s))^{\frac{1-\beta}2}\right\|^2_{L(\mathcal{H})}ds\leq C t^{\beta}\leq C.
\end{eqnarray}
For $I_4$ we use It\^o isometry \eqref{jumpint}, insert an appropriate power of $A_h(s)$, using \eqref{equi2} and \eqref{propdisevol4} with $\alpha=0$ and $\gamma=\frac{2H+\beta-2}2$ (if $2H+\beta-2\geq 0$) or \eqref{equi1} and \eqref{propdisevol3} with $\alpha=0$ and $\gamma=-\frac{2H+\beta-2}2$ (if $2H+\beta-2<0$), \eqref{equi3} and \assref{noiseandjump} (more precisely \eqref{jump}), to get
\begin{eqnarray}
\label{bound*4}
I_4&=&\left\|\int_0^t\int_{\chi}U_h(t,s)P_h\psi(z)\tilde{N}(dz,ds)\right\|^2_{L^2(\Omega,\mathcal{H})}\nonumber\\
&=&  \int_0^t\int_{\chi}\left\|U_h(t,s)P_h\psi(z)\right\|^2\nu(dz)ds\nonumber\\
&\leq&  \int_0^t\int_{\chi}\left\|U_h(t,s)(A_h(s))^{-\frac{2H+\beta-2}2}\right\|^2_{L(\mathcal{H})}\left\|(A_h(s))^{\frac{2H+\beta-2}2}P_h\psi(z)\right\|^2\nu(dz)ds.\nonumber\\
\end{eqnarray}
Finally 
\begin{eqnarray}
\label{bound*44}
I_4&\leq&  C\left(\int_0^t(t-s)^{2H+\beta-2}ds\right)\left(\int_{\chi}\left\|(A(0))^{\frac{2H+\beta-2}2}\psi(z)\right\|^2\nu(dz)\right)\nonumber\\
&\leq&  C t^{2H+\beta-1}\leq C.
\end{eqnarray}
The combination of \eqref{bound*1}--\eqref{bound*4} and continuous Gronwall's lemma completes the proof of \eqref{bound}. The estimate $\left\|F(X^h(t))\right\|_{L^2(\Omega,\mathcal{H})}\leq C$ is just a consequence of \eqref{bound} using \assref{nonlinear}. \\
Let us now prove \eqref{spareg}, pre-multiplying \eqref{soldispb} by $(A_h(\tau))^{\frac{2H+\beta-1}2}$, taking the square of the both side and using triangle inequality yields
\begin{eqnarray*}
\left\|(A_h(\tau))^{\frac{2H+\beta-1}2}X^h(t)\right\|_{L^2(\Omega,\mathcal{H})}
&\leq&\left\|(A_h(\tau))^{\frac{2H+\beta-1}2}U_h(t,0)P_hX_0\right\|_{L^2(\Omega,\mathcal{H})}\nonumber\\
&+&\int_0^t\left\|U_h(t,s)P_hF(s,X^h(s))\right\|_{L^2(\Omega,\mathcal{H})}ds\nonumber\\
&+&\left\|\int_0^t(A_h(\tau))^{\frac{2H+\beta-1}2}U_h(t,s)P_h\phi(s)dB^H(s)\right\|_{L^2(\Omega,\mathcal{H})}\nonumber\\
&+&\left\|\int_0^t\int_{\chi}(A_h(\tau))^{\frac{2H+\beta-1}2}U_h(t,s)P_h\psi(z)\tilde{N}(dz,ds)\right\|_{L^2(\Omega,\mathcal{H})}\nonumber\\
&=:& \sum_{i=1}^4 II_i.\nonumber
\end{eqnarray*}
Following the same lines as that in \cite[(3.29)-(3.30)]{Tamma}, it holds that
\begin{eqnarray}
\label{spareg1}
II_1+II_2\leq C
\end{eqnarray}
For the estimate $II_3$, using triangle inequality, the estimate $(a+b)^2\leq 2a^2+2b^2$, we have
\begin{eqnarray}
\label{spareg2*}
II_3^2&=&\left\|\int_0^t(A_h(\tau))^{\frac{2H+\beta-1}2}U_h(t,s)P_h\phi(s)dB^H(s)\right\|_{L^2(\Omega,\mathcal{H})}\nonumber\\
&\leq& 4\left\|\int_0^t(A_h(\tau))^{\frac{2H+\beta-1}2}U_h(t,s)P_h(\phi(t)-\phi(s))dB^H(s)\right\|_{L^2(\Omega,\mathcal{H})}\nonumber\\
&+& 4\left\|\int_0^t(A_h(\tau))^{\frac{2H+\beta-1}2}U_h(t,s)P_h\phi(t)dB^H(s)\right\|_{L^2(\Omega,\mathcal{H})}\nonumber\\
&=:&4II_{31}^2+4II_{32}^2.
\end{eqnarray}
Firstly using \eqref{fracint2}, inserting an appropriate power of $A_h(s)$, \eqref{equi3}, \eqref{propdisevol5} and \assref{noiseandjump} \eqref{noise2} yields
\begin{eqnarray}
\label{spareg21}
II_{31}^2&=&\left\|\int_0^t(A_h(\tau))^{\frac{2H+\beta-1}2}U_h(t,s)P_h(\phi(t)-\phi(s))dB^H(s)\right\|^2_{L^2(\Omega,\mathcal{H})}\nonumber\\
&\leq& C\int_0^t\left\|(A_h(\tau))^{\frac{2H+\beta-1}2}U_h(t,s)(A_h(s))^{\frac{1-\beta}2}\right\|^2_{L(\mathcal{H})}\left\|(A_h(s))^{\frac{\beta-1}2}P_h(\phi(t)-\phi(s))\right\|^2_{L^0_2}ds\nonumber\\
&\leq& C\int_0^t(t-s)^{-2H}\left\|(A(0))^{\frac{\beta-1}2}(\phi(t)-\phi(s))\right\|^2_{L^0_2}ds\nonumber\\
&\leq& C\int_0^t(t-s)^{-2H}(t-s)^{2\delta}ds\nonumber\\
&\leq& C\int_0^t(t-s)^{2\delta-2H}ds\leq C t^{2\delta-2H+1}\leq C.
\end{eqnarray}
Secondly, using \eqref{fracint1}, inserting $(A_h(s))^{\frac{1-\beta}2}(A_h(s))^{\frac{\beta-1}2}$, \eqref{equi3}, \assref{noiseandjump} \eqref{noise1} and \eqref{inte4} with $\alpha=2H+\beta-1$, $\kappa=1-\beta$, we obtain
\begin{eqnarray}
\label{spareg22}
&&II_{32}^2\nonumber\\
&=&\left\|\int_0^t(A_h(\tau))^{\frac{2H+\beta-1}2}U_h(t,s)P_h\phi(t)dB^H(s)\right\|^2_{L^2(\Omega,\mathcal{H})}\nonumber\\
&\leq& C\sum_{i\in\N^d}\left(\int_0^t\left\|(A_h(\tau))^{\frac{2H+\beta-1}2}U_h(t,s)(A_h(s))^{\frac{1-\beta}2}\right\|^{\frac1H}_{L(\mathcal{H})}\left\|(A_h(s))^{\frac{\beta-1}2}P_h\phi(t)Q^{\frac 12}e_i\right\|^{\frac 1H}ds\right)^{2H}\nonumber\\
&\leq& C\left(\int_0^t\left\|(A_h(\tau))^{\frac{2H+\beta-1}2}U_h(t,s)(A_h(s))^{\frac{1-\beta}2}\right\|^{\frac1H}_{L(\mathcal{H})}ds\right)^{2H}\left(\sum_{i\in\N^d}\left\|(A_h(s))^{\frac{\beta-1}2}P_h\phi(t)Q^{\frac 12}e_i\right\|^2\right)\nonumber\\
&\leq& C\left\|(A(0))^{\frac{\beta-1}2}\phi(t)\right\|^2_{L^0_2}\leq C.
\end{eqnarray}
Putting \eqref{spareg21} and \eqref{spareg22} in \eqref{spareg2*} yields
\begin{eqnarray}
\label{spareg2}
II_3^2\leq C.
\end{eqnarray}
For the estimate of $II_4$, using It\^o isometry \eqref{jumpint}, inserting $(A_h(s))^{-\frac{2H+\beta-2}2}(A_h(s))^{\frac{2H+\beta-2}2}$, \eqref{equi2} and \eqref{inte2} with $\rho=\frac{2H+\beta-1}2$, $\gamma=\eta=\frac{2H+\beta-2}2$ (if $2H+\beta-2\geq0$) or \eqref{equi1} and \eqref{inte3} with $\rho=\frac{2H+\beta-1}2$, $\gamma=-\frac{2H+\beta-2}2$ if ($2H+\beta-2<0$), \eqref{equi3} and \assref{noiseandjump} \eqref{jump}, it holds that
\begin{eqnarray}
\label{spareg3}
&&II_4^2\nonumber\\
&=&\left\|\int_0^t\int_{\chi}(A_h(\tau))^{\frac{2H+\beta-1}2}U_h(t,s)P_h\psi(z)\tilde{N}(dz,ds)\right\|^2_{L^2(\Omega,\mathcal{H})}\nonumber\\
&=&\int_0^t\int_{\chi}\left\|(A_h(\tau))^{\frac{2H+\beta-1}2}U_h(t,s)(A_h(s))^{-\frac{2H+\beta-2}2}\right\|^2_{L(\mathcal{H})}\left\|(A_h(s))^{\frac{2H+\beta-2}2}P_h\psi(z)\right\|^2\nu(dz)ds,\nonumber\\
\end{eqnarray}
so
\begin{eqnarray}
&&II_4^2\nonumber\\
&\leq& C\left(\int_0^t\left\|(A_h(\tau))^{\frac{2H+\beta-1}2}U_h(t,s)(A_h(s))^{-\frac{2H+\beta-2}2}\right\|^2_{L(\mathcal{H})}ds\right)\left(\int_{\chi}\left\|(A(0))^{\frac{2H+\beta-2}2}\psi(z)\right\|^2\nu(dz)\right)\nonumber\\
&\leq& C\left(\int_{\chi}\left\|(A(0))^{\frac{2H+\beta-2}2}\psi(z)\right\|^2\nu(dz)\right)\leq C.
\end{eqnarray}
Combining \eqref{spareg1}, \eqref{spareg2} and \eqref{spareg3} completes the proof of \eqref{spareg}.\\
Now for the proof of \eqref{timreg}, using triangle inequality yields
\begin{eqnarray}
\label{timreg*}
\left\|X^h(t_2)-X^h(t_1)\right\|_{L^2(\Omega,\mathcal{H})}
&\leq&\left\|\left(U_h(t_2,t_1)-\mathbf{I}\right)X^h(t_1)\right\|_{L^2(\Omega,\mathcal{H})}\nonumber\\
&+&\left\|\int_{t_1}^{t_2}U_h(t_2,s)P_h\phi(s)dB^H(s)\right\|_{L^2(\Omega,\mathcal{H})}\nonumber\\
&+&\left\|\int_{t_1}^{t_2}\int_{\chi}U_h(t_2,s)P_h\psi(z)\tilde{N}(dz,ds)\right\|_{L^2(\Omega,\mathcal{H})}\nonumber\\
&=:& \sum_{i=1}^3 III_i.
\end{eqnarray}
Using \lemref{propdisevol} and \eqref{spareg} yields
\begin{eqnarray}
\label{timreg1}
 III_1&\leq& \left\|\left(U_h(t_2,t_1)-\mathbf{I}\right)(A(t_1))^{-\frac{2H+\beta-1}{2}}\right\|_{L(\mathcal{H})}\left\|(A(t_1))^{\frac{2H+\beta-1}{2}} X^h(t_1)\right\|_{L^2(\Omega,\mathcal{H})}\nonumber\\
 &\leq &C(t_2-t_1)^{\frac{2H+\beta-1}2}.
\end{eqnarray}
Triangle inequality, the estimate $(a+b)^2\leq 2a^2+2b^2$, \eqref{fracint2}, \eqref{fracint1}, inserting $(A_h(s))^{\frac{1-\beta}2}(A_h(s))^{\frac{\beta-1}2}$ yields
\begin{eqnarray*}
&&III_2^2\nonumber\\
&=&\left\|\int_{t_1}^{t_2}U_h(t_2,s)P_h\phi(s)dB^H(s)\right\|^2_{L^2(\Omega,\mathcal{H})}\nonumber\\
&\leq& 2\left\|\int_{t_1}^{t_2}U_h(t_2,s)P_h(\phi(t_2)-\phi(s))dB^H(s)\right\|^2_{L^2(\Omega,\mathcal{H})}+2\left\|\int_{t_1}^{t_2}U_h(t_2,s)P_h\phi(t_2)dB^H(s)\right\|^2_{L^2(\Omega,\mathcal{H})}\nonumber\\
&\leq& C\int_{t_1}^{t_2}\left\|U_h(t_2,s)(A_h(s))^{\frac{1-\beta}2}\right\|^{2}_{L(\mathcal{H})}\left\|(A_h(s))^{\frac{\beta-1}2}P_h(\phi(t_2)-\phi(s))\right\|^2_{L^0_2}ds\nonumber\\
&+&C\sum_{i\in\N^d}\left(\int_{t_1}^{t_2}\left\|U_h(t_1,s)(A_h(s))^{\frac{1-\beta}2}\right\|^{\frac 1H}_{L(\mathcal{H})}\left\|(A_h(s))^{\frac{\beta-1}2}P_h\phi(t_1)Q^{\frac 12}e_i\right\|^{\frac 1H}ds\right)^{2H}\nonumber
\end{eqnarray*}
Using \eqref{equi3}, \eqref{propdisevol3}, \eqref{inte4} with $\alpha=0$ and $\kappa=1-\beta$, \assref{noiseandjump} yields
\begin{eqnarray}
\label{timreg3}
III_2^2
&\leq& C\int_{t_1}^{t_2}(t_2-s)^{\beta-1}(t_2-s)^{2\delta}ds\nonumber\\
&+&C\left(\int_{t_1}^{t_2}\left\|U_h(t_1,s)(A_h(s))^{\frac{1-\beta}2}\right\|^{\frac 1H}_{L(\mathcal{H})}\right)^{2H}\left\|(A(0))^{\frac{\beta-1}2}\phi(t_2)\right\|^2\nonumber\\
&\leq& C(t_2-t_1)^{2\delta+\beta}+C(t_2-t_1)^{2H+\beta-1}\nonumber\\
&\leq& C(t_2-t_1)^{2H+\beta-1}.
\end{eqnarray}
%using the It\^o isometry \eqref{jumpint}, \lemref{equi}, inserting an appropriate power of  $A_h(t_1)$ and $A_h(s)$, \eqref{equi3} \eqref{propdisevol6}, \eqref{inte2} with $\rho=2H+\beta-1$ and $\eta=2H+\beta-2$ if $2H+\beta-2\geq0$ or \eqref{inte3} with $\rho=2H+\beta-1$ and $\gamma=-(2H+\beta-2)$ if $2H+\beta-2<0$, \eqref{equi3}, \assref{noiseandjump} \eqref{jump} and the $\sigma$-finiteness of the measure $\nu$, we obtain
%{\small
%\begin{eqnarray}
%\label{timreg4}
%&&III_6^2\nonumber\\
%&=&\left\|\int_0^{t_1}\int_{\chi}\left(U_h(t_2,s)-U_h(t_1,s)\right)P_hz_0\tilde{N}(dz,ds)\right\|^2_{L^2(\Omega,\mathcal{H})}\nonumber\\
%&=&\left\|\int_0^{t_1}\int_{\chi}\left(U_h(t_2,t_1)-I\right)U_h(t_1,s)P_hz_0\tilde{N}(dz,ds)\right\|^2_{L^2(\Omega,\mathcal{H})}\nonumber\\
%&=&\int_0^{t_1}\int_{\chi}\left\|\left(U_h(t_2,t_1)-I\right)(A_h(t_1))^{-\frac{2H+\beta-1}2}\right\|^2_{L(\mathcal{H})}\left\|(A_h(t_1))^{\frac{2H+\beta-1}2}U_h(t_1,s)(A_h(s))^{-\frac{2H+\beta-2}2}\right\|^2_{L(\mathcal{H})}\nonumber\\
%&&\left\|(A_h(s))^{\frac{2H+\beta-2}2}P_hz_0\right\|^2\nu(dz)ds\nonumber\\
%&\leq&C(t_2-t_1)^{2H+\beta-1}\left(\int_0^{t_1}\left\|(A_h(t_1))^{\frac{2H+\beta-1}2}U_h(t_1,s)(A_h(s))^{-\frac{2H+\beta-2}2}\right\|^2_{L(\mathcal{H})}ds\right)\nonumber\\
%&&\left(\int_{\chi}\left\|(A(0))^{\frac{2H+\beta-2}2}z_0\right\|^2\nu(dz)\right)\nonumber\\
%&\leq&C(t_2-t_1)^{2H+\beta-1}.
%\end{eqnarray}
%}
Using It\^o isometry \eqref{jumpint}, the insertion $(A_h(s))^{-\frac{2H+\beta-2}2}(A_h(s))^{\frac{2H+\beta-2}2}$, \eqref{propdisevol1} if $2H+\beta-2\geq 0$ or \eqref{propdisevol3} if $2H+\beta-2<0$, and \eqref{equi3}, \assref{noiseandjump} \eqref{jump} yields
\begin{eqnarray}
\label{timreg5}
III_3^2&=&\left\|\int_{t_1}^{t_2}\int_{\chi}U_h(t_2,s)P_h\psi(z)\tilde{N}(dz,ds)\right\|^2_{L^2(\Omega,\mathcal{H})}\nonumber\\
&=&\int_{t_1}^{t_2}\int_{\chi}\left\|U_h(t_2,s)P_h\psi(z)\right\|^2\nu(dz)ds\nonumber\\
&\leq&C\left(\int_{t_1}^{t_2}\left\|U_h(t_2,s)(A_h(s))^{-\frac{2H+\beta-2}2}\right\|^2_{L(\mathcal{H})}ds\right)\left(\int_{\chi}\left\|(A_h(s))^{\frac{2H+\beta-2}2}P_h\psi(z)\right\|^2\nu(dz)\right)\nonumber\\
&\leq&C\left(\int_{t_1}^{t_2}(t_2-s)^{\min(0,2H+\beta-2)}ds\right)\left(\int_{\chi}\left\|(A(0))^{\frac{2H+\beta-2}2}\psi(z)\right\|^2\nu(dz)\right)\nonumber\\
&\leq&C(t_2-t_1)^{\min(1,2H+\beta-1)}.
\end{eqnarray}
Putting \eqref{timreg1}-\eqref{timreg5} in \eqref{timreg*} completes the proof of \eqref{timreg} and therefore that of \lemref{reg}.$\hfill\square$

\subsection{Spatial error}
As in \cite{Tammab} we consider the following deterministic problem: find $u\in V$ such that
\begin{eqnarray*}
\label{detpb}
\frac{du}{dt}+A(t)u=0,\quad u(\tau)=v,\quad t\in(\tau,T],\quad \tau\geq 0.
\end{eqnarray*}
The corresponding semi-discrete problem in space is: find $u^h\in V_h$ such that 
\begin{eqnarray}
\label{detdispb}
\frac{du^h}{dt}+A_h(t)u^h=0,\quad u^h(\tau)=P_h v,\quad t\in(\tau,T].
\end{eqnarray}
Let us define the operator 
\begin{eqnarray*}
\label{semiop}
T_h(t,\tau):=U(t,\tau)-U_h(t,\tau)P_h.
\end{eqnarray*} 
Then we have the following lemma.
\begin{lem}(\cite[Lemma 3.7]{Tammab})
\label{semierror}
Let \assref{linear} be fulfilled. Then the following estimate holds for the semi-discrete approximation of \eqref{detdispb}. The exists a constant $C>0$ such that for $v\in \mathcal{\dot H}^{\gamma}$
\begin{eqnarray*}
\label{semierror1}
\|u(t)-u^h(t)\|=\|T_h(t,\tau)v\|\leq Ch^r(t-\tau)^{-(r-\gamma)/2}\|v\|_{\gamma},\quad t\in(0,T],\quad r\in[0,2],\quad\gamma\leq r.
\end{eqnarray*}
\end{lem}
The following lemma provides an estimate 
in mean square sense for the error between the solution of the SPDE \eqref{SPDE} 
and the spatially semi discrete approximation \eqref{soldispb}.
\begin{lem}
\label{spaerr}
Let $X$ and $X^h$ be the mild solutions of \eqref{SPDE} and \eqref{disSPDE}, respectively. Let Assumptions \ref{linear}--\ref{initdata} be fulfilled, then there exists a constant $C$ independent of $h$, such that
\begin{eqnarray}
\label{boundspaerr}
\|X(t)-X^h(t)\|_{L^2(\Omega,\mathcal{H})}\leq Ch^{2H+\beta-1-\epsilon},
\end{eqnarray} 
where $\epsilon$ is a positive constant small enough.
\end{lem}
\textit{Proof:} Define $e(t):=X(t)-X^h(t)$. By \eqref{solpb} and \eqref{soldispb}, taking the norm and using triangle inequality, we deduce
\begin{eqnarray*}
\|e(t)\|_{L^2(\Omega,\mathcal{H})}&\leq& \|U(t,0)X_0-U_h(t,0)P_hX_0\|_{L^2(\Omega,\mathcal{H})}\nonumber\\
&+&\left\|\int_0^t[U(t,s)F(s,X(s))-U_h(t,s)P_hF(s,X^h(s))]ds\right\|_{L^2(\Omega,\mathcal{H})}\nonumber\\
&+&\left\|\int_0^t[U(t,s)\phi(s)-U_h(t,s)P_h\phi(s)]dB^H(s)\right\|_{L^2(\Omega,\mathcal{H})}\nonumber\\
&+&\left\|\int_0^t\int_{\chi}[U(t,s)\psi(z)-U_h(t,s)P_h\psi(z)]\tilde{N}(dz,ds)\right\|_{L^2(\Omega,\mathcal{H})}\nonumber\\
&=:&\sum_{i=1}^4 IV_i. %\label{boundspaerr*}
\end{eqnarray*} 
We will analyse the above terms $IV_i$, $i=1,2,3,4$ one by one. For the first term $IV_1$, using \lemref{semierror} (i) with $r=\gamma=2H+\beta-1$ and \assref{initdata} yields
\begin{eqnarray}
\label{boundspaerr*1}
IV_1&=&\|U(t,0)X_0-U_h(t,0)P_hX_0\|_{L^2(\Omega,\mathcal{H})}= \|T_h(t,0)X_0\|_{L^2(\Omega,\mathcal{H})}\nonumber\\
&\leq& C h^{2H+\beta-1}\|X_0\|_{L^2(\Omega,\mathcal{\dot H}^{2H+\beta-1})} \leq C h^{2H+\beta-1}.
\end{eqnarray}
For the second term $IV_2$, following the same line as \cite[(3.47)]{Tammab}
\begin{eqnarray}
%\label{boundspaerr*2}
IV_2&=&\left\|\int_0^t[U(t,s)F(s,X(s))-U_h(t,s)P_hF(s,X^h(s))]ds\right\|_{L^2(\Omega,\mathcal{H})}\nonumber\\
&\leq& C h^{2H+\beta-1}+C\int_0^t\|X(s)-X^h(s)\|_{L^2(\Omega,\mathcal{H})}ds.
\end{eqnarray}
For the third term $IV_3$, by adding and subtracting a term, using the triangle inequality and the estimate $(a+b)^2\leq 2a^2+2b^2$ and  \eqref{fracint2},  we have
\begin{eqnarray*}
IV_3^2&=&\left\|\int_0^t[U(t,s)\phi(s)-U_h(t,s)P_h\phi(s)]dB^H(s)\right\|^2_{L^2(\Omega,\mathcal{H})}\nonumber\\
&=&\left\|\int_0^t[T_h(t,s)\phi(s)]dB^H(s)\right\|^2_{L^2(\Omega,\mathcal{H})}\nonumber\\
&\leq&4\int_0^t\left\|T_h(t,s)(\phi(t)-\phi(s))\right\|^2_{L^0_2}ds+4C\sum_{i\in\N^d}\left(\int_0^t\left\|T_h(t,s)\phi(t)Q^{\frac 12}e_i\right\|^{\frac 1H}ds\right)^{2H}.
\end{eqnarray*}
Using \lemref{semierror} with $r=2H+\beta-1$, $\gamma=\beta-1$, \assref{noiseandjump} (more precisely \eqref{noise2}) for the first term and  \eqref{fracint1}, \lemref{semierror} with $r=2H+\beta-1-\epsilon$, $\gamma=\beta-1$, \assref{noiseandjump} (more precisely \eqref{noise1}) for the second term, yields
\begin{eqnarray}
\label{boundspaerr*3}
IV_3^2&\leq&Ch^{2(2H+\beta-1)}\int_0^t(t-s)^{-2H}(t-s)^{2\delta}ds\nonumber\\
&+&C\sum_{i\in\N^d}\left(\int_0^th^{\frac{2H+\beta-1-\epsilon}H}(t-s)^{-1+\frac{\epsilon}{2H}}\left\|(A(0))^{\frac{\beta-1}2}\phi(t)Q^{\frac 12}e_i\right\|^{\frac 1H}ds\right)^{2H}\nonumber\\
&\leq&Ch^{2(2H+\beta-1)}t^{2\delta-2H+1}+Ch^{2(2H+\beta-1-\epsilon)}\left(\int_0^t(t-s)^{-1+\frac{\epsilon}{2H}}ds\right)^{2H}\left\|(A(0))^{\frac{\beta-1}2}\phi(t)\right\|^2_{L^0_2}\nonumber\\
&\leq&Ch^{2(2H+\beta-1-\epsilon)}.
\end{eqnarray}
For the fourth term $IV_4$, using the It\^{o} isometry \eqref{jumpint}, \lemref{semierror} with $r=2H+\beta-1-\epsilon$, $\gamma=2H+\beta-2$ and \assref{noiseandjump} \eqref{jump}, we derive
\begin{eqnarray}
\label{boundspaerr*4}
IV_4^2&=&\left\|\int_0^t\int_{\chi}[U(t,s)\psi(z)-U_h(t,s)P_h\psi(z)]\tilde{N}(dz,ds)\right\|^2_{L^2(\Omega,\mathcal{H})}\nonumber\\
&=&\int_0^t\int_{\chi}\left\|T_h(t,s)\psi(z)\right\|^2\nu(dz)ds\nonumber\\
&\leq&\int_{\chi}\int_0^t C h^{2(2H+\beta-1-\epsilon)}(t-s)^{-1+\epsilon}\left\|(A(0))^{\frac{2H+\beta-2}2}\psi(z)\right\|^2\nu(dz)ds\nonumber\\
&\leq&Ch^{2(2H+\beta-1-\epsilon)}t^{\epsilon}\int_{\chi}\left\|(A(0))^{\frac{2H+\beta-2}2}\psi(z)\right\|^2\nu(dz)\nonumber\\
&\leq&Ch^{2(2H+\beta-1-\epsilon)}.
\end{eqnarray}
A combination of the estimates \eqref{boundspaerr*1}-\eqref{boundspaerr*4} gives
\begin{eqnarray*}
\|e(t)\|_{L^2(\Omega,\mathcal{H})}\leq C h^{2H+\beta-1-\epsilon}+C\int_0^t\|e(s)\|_{L^2(\Omega,\mathcal{H})}ds.\nonumber
\end{eqnarray*}
Applying Gronwall's inequality finally completes the proof.$\hfill\square$

\section{Proof of \thmref{mainthm} for SMTI scheme}
\label{convSMTI}
Before moving to the proof of \thmref{mainthm}, we first present some preparatory results useful in our convergence proof. \subsection{Preliminaries results}
\begin{lem}
\label{prepaSMTI1}
Let \assref{linear} be fulfilled.
\begin{enumerate}
\item [(i)] For all $\alpha\in[-1,1]$, the following estimate holds
\begin{eqnarray*}
\left\|\left(U_h(t_j,t_{j-1})-e^{-\Delta t A_{h,j-1}}\right)\left(A_{h,j-1}\right)^{\alpha}\right\|_{L(\mathcal{H})}\leq C\Delta t^{1-\alpha}. \nonumber
\end{eqnarray*}
\item [(ii)] For any $0\leq s\leq t\leq T$, $r\in[0,T]$, $\gamma\in[0,1]$ and $\eta\in[-1,1)$ such that $0\leq \gamma-\eta\leq 1$, the following estimate holds
\begin{eqnarray*}
\left\|(A_{h}(r))^{-\gamma}\left(U_h(t,s)-I\right)(A_h(s))^{\eta}\right\|_{L(\mathcal{H})}\leq C\Delta t^{\gamma-\eta}. \nonumber
\end{eqnarray*}
\end{enumerate}
\end{lem} 
\textit{Proof:} See \cite[Lemma 3.10,(i),(ii)]{Tammab} for the proof of (i). Now concerning the proof of (ii), from \eqref{disevol1}, triangle inequality, \lemref{equi} (i), smoothing properties of discrete semigroup $S^h_{\tau}(t)$ and \cite[Lemma 3.10 (i)]{Tammab}, it holds that
\begin{eqnarray*}
&&\left\|(A_{h}(r))^{-\gamma}\left(U_h(t,s)-I\right)(A_h(s))^{\eta}\right\|_{L(\mathcal{H})}\nonumber\\
&\leq& \left\|(A_{h}(r))^{-\gamma}\left(e^{-(t-s)A_h(s)}-I\right)(A_h(s))^{\eta}\right\|_{L(\mathcal{H})}\nonumber\\
&+&\int_s^t\left\|(A_{h}(r))^{-\gamma}S^h_{\tau}(t-\tau)R^h(\tau,s)(A_h(s))^{\eta}\right\|_{L(\mathcal{H})}d\tau\nonumber\\
&\leq& C\left\|(A_h(s))^{-\gamma+\eta}\left(e^{-(t-s)A_h(s)}-I\right)\right\|_{L(\mathcal{H})}+C\int_s^t\left\|R^h(\tau,s)(A_h(s))^{\eta}\right\|_{L(\mathcal{H})}d\tau\nonumber\\
&\leq& C(t-s)^{\gamma-\eta}+C\int_s^t(\tau-s)^{\min(-\eta,0)}d\tau\nonumber\\
&\leq& C(t-s)^{\gamma-\eta}+C(t-s)^{\min(1-\eta,1)}\leq C(t-s)^{\gamma-\eta}.\nonumber\hfill\square 
\end{eqnarray*}
As in \cite{Mukmaa,Mukmab,Tamma,Tammab}, for non commutative operators $H_j$ in a Banach space, we introduce the following notation.
\begin{eqnarray*}
\prod_{j=l}^kH_j:=\left\{
  \begin{array}{ll}
    H_kH_{k-1}\cdot\cdot\cdot H_l,\hspace{1cm}\text{if}\hspace{0.5cm}k\geq l \\
    I,\hspace{3.5cm}\text{if}\hspace{0.5cm} k< l .
  \end{array}
\right.
\end{eqnarray*}
\begin{lem}(\cite[Lemma 3.9]{Tammab})
\label{prepaSMTI3}
Let \assref{linear} be fulfilled. Then the following estimate holds
\begin{eqnarray*}
\label{propsemi}
\left\|\left(\prod_{j=l}^m e^{-\Delta t A_{h,j}}\right)(A_{h,l})^{\gamma}\right\|_{L(\mathcal{H})}\leq C t_{m+1-l}^{-\gamma},\quad 0\leq l<m,\quad 0\leq \gamma<1.
\end{eqnarray*}  
\end{lem}
\begin{lem}(\cite[Lemma 2.15]{Mukmaa})
\label{prepaSMTI2}
\begin{enumerate}
\item [(i)] The following estimate holds for $1\leq i\leq m$
\begin{eqnarray*}
\left\|\left(\prod_{j=i}^m U_h(t_j,t_{j-1})\right)-\left(\prod_{j=i-1}^{m-1} e^{-\Delta t A_{h,j}}\right)\right\|_{L(\mathcal{H})}\leq C \Delta t^{1-\epsilon}.\nonumber
\end{eqnarray*}  
\item [(ii)] The following estimate holds for $1\leq i\leq m$
\begin{eqnarray*}
\left\|\left(\prod_{j=i}^m U_h(t_j,t_{j-1})\right)-\left(\prod_{j=i-1}^{m-1} e^{-\Delta t A_{h,j}}\right)(A_h(0))^{-\epsilon}\right\|_{L(\mathcal{H})}\leq C \Delta t.\nonumber
\end{eqnarray*}  
\item [(iii)] For all $1\leq i\leq m\leq M$, for any $\alpha\in[0,1]$, the following estimate holds
\begin{eqnarray*}
\left\|\left(\prod_{j=i}^m U_h(t_j,t_{j-1})\right)-\left(\prod_{j=i-1}^{m-1} e^{-\Delta t A_{h,j}}\right)(A_{h,i-1})^{\alpha}\right\|_{L(\mathcal{H})}\leq C \Delta t^{1-\alpha-\epsilon}t_{m-i+1}^{-\alpha+\epsilon}.\nonumber
\end{eqnarray*}  
\end{enumerate}
for arbitrary $\epsilon>0$.
\end{lem}
With these three lemmas, we are now in position to prove our  main theorem for the SMTI scheme. 
\subsection{Proof of \thmref{mainthm} for Magnus-type integrator}
 Using triangle inequality yields
\begin{eqnarray*}
\label{err}
\|X(t_m)-X^h_m\|_{L^2(\Omega,\mathcal{H})}\leq\|X(t_m)-X^h(t_m)\|_{L^2(\Omega,\mathcal{H})}+\|X^h(t_m)-X^h_m\|_{L^2(\Omega,\mathcal{H})}.
\end{eqnarray*}
The space error is estimated in \lemref{spaerr}. It remains to estimate the temporal error. Let us recall that the mild solution \eqref{soldispb} can be written as follows
\begin{eqnarray}
\label{mildsol}
X^h(t_m)&=&U_h(t_m,t_{m-1})P_hX_{m-1}^h+\int_{t_{m-1}}^{t_m}U_h(t_m,s)P_hF(s,X^h(s))ds\nonumber\\
&+&\int_{t_{m-1}}^{t_m}U_h(t_m,s)P_h\phi(s)dB^H(s)+\int_{t_{m-1}}^{t_m}\int_{\chi}U_h(t_m,s)P_h\psi(z)\tilde{N}(dz,ds).
\end{eqnarray}
Iterating the mild solution \eqref{mildsol} yields
\begin{eqnarray}
\label{mildsol*}
X^h(t_m)&=&\left(\prod_{j=1}^m U_h(t_j,t_{j-1})\right)P_hX_0+\int_{t_{m-1}}^{t_m}U_h(t_m,s)P_hF(s,X^h(s))ds\nonumber\\
&+&\int_{t_{m-1}}^{t_m}U_h(t_m,s)P_h\phi(s)dB^H(s)+\int_{t_{m-1}}^{t_m}\int_{\chi}U_h(t_m,s)P_h\psi(z)\tilde{N}(dz,ds)\nonumber\\
&+&\sum_{k=1}^{m-1}\int_{t_{m-k-1}}^{t_{m-k}}\left(\prod_{j=m-k+1}^m U_h(t_j,t_{j-1})\right)U_h(t_{m-k},s)P_hF(s,X^h(s))ds\nonumber\\
&+&\sum_{k=1}^{m-1}\int_{t_{m-k-1}}^{t_{m-k}}\left(\prod_{j=m-k+1}^m U_h(t_j,t_{j-1})\right)U_h(t_{m-k},s)P_h\phi(s)dB^H(s)\nonumber\\
&+&\sum_{k=1}^{m-1}\int_{t_{m-k-1}}^{t_{m-k}}\int_{\chi}\left(\prod_{j=m-k+1}^m U_h(t_j,t_{j-1})\right)U_h(t_{m-k},s)P_h\psi(z)\tilde{N}(dz,ds).
\end{eqnarray}
Iterating the numerical solution \eqref{SMTI} as in \cite{Mukmaa} yields
\begin{eqnarray}
\label{SMTI*}
X^h_m&=&\left(\prod_{j=0}^{m-1} e^{-\Delta t A_{h,j}}\right)X^h_0+\int_{t_{m-1}}^{t_m}e^{-(t_m-s)A_{h,m-1}}P_hF(t_{m-1},X^h_{m-1})ds\nonumber\\
&+&\int_{t_{m-1}}^{t_m}e^{-(t_m-s)A_{h,m-1}}P_h\phi(t_{m-1})dB^H(s)+\int_{t_{m-1}}^{t_m}\int_{\chi}e^{-(t_m-s)A_{h,m-1}}P_h\psi(z)\tilde{N}(dz,ds)\nonumber\\
&+&\sum_{k=1}^{m-1}\int_{t_{m-k-1}}^{t_{m-k}}\left(\prod_{j=m-k}^{m-1} e^{-\Delta t A_{h,j}}\right)e^{-\Delta t A_{h,m-k-1}}P_hF(t_{m-k-1},X^h_{m-k-1})ds\nonumber\\
&+&\sum_{k=1}^{m-1}\int_{t_{m-k-1}}^{t_{m-k}}\left(\prod_{j=m-k}^{m-1} e^{-\Delta t A_{h,j}}\right)e^{-\Delta t A_{h,m-k-1}}P_h\phi(t_{m-k-1})dB^H(s)\nonumber\\
&+&\sum_{k=1}^{m-1}\int_{t_{m-k-1}}^{t_{m-k}}\int_{\chi}\left(\prod_{j=m-k}^{m-1} e^{-\Delta t A_{h,j}}\right)e^{-\Delta t A_{h,m-k-1}}P_h\psi(z)\tilde{N}(dz,ds).
\end{eqnarray}
Subtracting \eqref{SMTI*} from \eqref{mildsol*}, taking the $L^2$-norm and triangle inequality yields
\begin{eqnarray*}
\label{timerr*}
\|X^h(t_m)-X^h_m\|_{L^2(\Omega,\mathcal{H})}\leq \sum_{i=1}^7 \|V_i\|_{L^2(\Omega,\mathcal{H})}.
\end{eqnarray*}
Next we will estimate these seven terms separately.
\subsubsection{Estimate of $V_1$, $V_2$ and $V_3$}
 For the first term, inserting $(A_h(0))^{-\frac{2H+\beta-1}2}$ $(A_h(0))^{\frac{2H+\beta-1}2}$, applying \lemref{prepaSMTI2} (i), \lemref{equi} (iii) with $\alpha=\frac{2H+\beta-1}2$ and \assref{initdata}, we get
\begin{eqnarray}
\label{timerr1}
\|V_1\|_{L^2(\Omega,\mathcal{H})}
&:=&\left\|\left(\prod_{j=1}^m U_h(t_j,t_{j-1})\right)P_hX_0-\left(\prod_{j=0}^{m-1} e^{-\Delta t A_{h,j}}\right)P_hX_0\right\|_{L^2(\Omega,\mathcal{H})}\nonumber\\
&\leq&\left\|\left[\left(\prod_{j=1}^m U_h(t_j,t_{j-1})\right)-\left(\prod_{j=0}^{m-1} e^{-\Delta t A_{h,j}}\right)\right](A_h(0))^{-\frac{2H+\beta-1}2}\right\|_{L(\mathcal{H})}\nonumber\\
&\times&\left\|(A_h(0))^{\frac{2H+\beta-1}2}P_hX_0\right\|_{L^2(\Omega,\mathcal{H})}\nonumber\\
&\leq& C\Delta t\left\|(A(0))^{\frac{2H+\beta-1}2}X_0\right\|_{L^2(\Omega,\mathcal{H})}\leq C\Delta t.
\end{eqnarray}
As in \cite[(3.79)]{Tammab} we have
\begin{eqnarray}
\label{timerr2}
&&\|V_2\|_{L^2(\Omega,\mathcal{H})}\nonumber\\
&:=&\left\| \int_{t_{m-1}}^{t_m}U_h(t_m,s)P_hF(s,X^h(s))ds-\int_{t_{m-1}}^{t_m}e^{-(t_m-s)A_{h,m-1}}P_hF(t_{m-1},X^h_{m-1})ds\right\|_{L^2(\Omega,\mathcal{H})} \nonumber\\
&\leq& C\Delta t+C\Delta t\|X^h(t_{m-1})-X^h_{m-1}\|_{L^2(\Omega,\mathcal{H})}.
\end{eqnarray}
To estimate  $\|V_3\|^2_{L^2(\Omega,\mathcal{H})}$, we use  triangle inequality and the estimate $(a+b+c)^2\leq 3a^2+3b^2+3c^2$ to split it into three terms as follows.
\begin{eqnarray*}
\label{timerr*3}
\|V_3\|^2_{L^2(\Omega,\mathcal{H})}&:=&\left\| \int_{t_{m-1}}^{t_m}U_h(t_m,s)P_h\phi(s)dB^H(s)-\int_{t_{m-1}}^{t_m}e^{-\Delta t A_{h,m-1}}P_h\phi(t_{m-1})dB^H(s)\right\|^2_{L^2(\Omega,\mathcal{H})}\nonumber\\
&\leq& 9\left\| \int_{t_{m-1}}^{t_m}U_h(t_m,s)P_h(\phi(s)-\phi(t_{m-1}))dB^H(s)\right\|^2_{L^2(\Omega,\mathcal{H})}\nonumber\\
&+&9\left\|\int_{t_{m-1}}^{t_m}\left(U_h(t_m,s)-U_h(t_m,t_{m-1})\right)P_h\phi(t_{m-1})dB^H(s)\right\|^2_{L^2(\Omega,\mathcal{H})}\nonumber\\
&+&9\left\|\int_{t_{m-1}}^{t_m}\left(U_h(t_m,t_{m-1})-e^{-\Delta t A_{h,m-1}}\right)P_h\phi(t_{m-1})dB^H(s)\right\|^2_{L^2(\Omega,\mathcal{H})}\nonumber\\
&=:&9\|V_{31}\|^2_{L^2(\Omega,\mathcal{H})}+9\|V_{32}\|^2_{L^2(\Omega,\mathcal{H})}+9\|V_{33}\|^2_{L^2(\Omega,\mathcal{H})}.
\end{eqnarray*}
Using \eqref{fracint2}, inserting $(A_h(s))^{\frac{1-\beta}2}(A_h(s))^{\frac{\beta-1}2}$, \eqref{propdisevol3}, \eqref{equi3} and \assref{noiseandjump} \eqref{noise2}, it holds that
\begin{eqnarray}
\label{timerr*31}
\|V_{31}\|^2_{L^2(\Omega,\mathcal{H})}&=&\left\| \int_{t_{m-1}}^{t_m}U_h(t_m,s)P_h(\phi(s)-\phi(t_{m-1}))dB^H(s)\right\|^2_{L^2(\Omega,\mathcal{H})}\nonumber\\
&\leq& C \int_{t_{m-1}}^{t_m}\left\|U_h(t_m,s)(A_h(s))^{\frac{1-\beta}2}\right\|^2_{L(\mathcal{H})}\left\|(A_h(s))^{\frac{\beta-1}2}P_h(\phi(s)-\phi(t_{m-1}))\right\|^2_{L^0_2}ds\nonumber\\
&\leq& C \int_{t_{m-1}}^{t_m}(t_m-s)^{\beta-1}\left\|(A(0))^{\frac{\beta-1}2}(\phi(s)-\phi(t_{m-1}))\right\|^2_{L^0_2}ds\nonumber\\
&\leq& C \int_{t_{m-1}}^{t_m}(t_m-s)^{\beta-1}(s-t_{m-1})^{2\delta}ds\nonumber\\
&\leq& C\Delta t^{2H+\beta-1} \left(\int_{t_{m-1}}^{t_m}(t_m-s)^{\beta-1}ds\right)\nonumber\\
&\leq& C\Delta t^{2H+\beta-1}.
\end{eqnarray}
Applying \eqref{fracint1}, \lemref{disevol} (iii), inserting an appropriate power of $A_h(s)$ and $A_{h,m-1}$, \eqref{propdisevol5} with $\gamma=H-\frac{\epsilon}2$ and $\eta=\frac{1-\beta}2$, \lemref{equi} (iii) with $\alpha=\frac{\beta-1}2$, \assref{noiseandjump} \eqref{noise1} yields 
\begin{eqnarray*}
%\label{timerr*32}
\|V_{32}\|^2_{L^2(\Omega,\mathcal{H})}
&=&\left\|\int_{t_{m-1}}^{t_m}\left(U_h(t_m,s)-U_h(t_m,t_{m-1})\right)P_h\phi(t_{m-1})dB^H(s)\right\|^2_{L^2(\Omega,\mathcal{H})}\nonumber\\
&\leq& C\sum_{i\in\N^d}\left(\int_{t_{m-1}}^{t_m}\left\|U_h(t_m,s)\left(U_h(s,t_{m-1})-I\right)(A_{h,m-1})^{\frac{1-\beta}2}\right\|^{\frac 1H}_{L(\mathcal{H})}\right.\nonumber\\
&&\left.\left\|(A_{h_m-1})^{\frac{\beta-1}2}P_h\phi(t_{m-1})Q^{\frac 12}e_i\right\|^{\frac 1H}ds\right)^{2H}\nonumber\\
&\leq& C\left(\sum_{i\in\N^d}\left\|(A_h(0))^{\frac{\beta-1}2}P_h\phi(t_{m-1})Q^{\frac 12}e_i\right\|^2\right)\left(\int_{t_{m-1}}^{t_m}\left\|U_h(t_m,s)(A_h(s))^{H-\frac{\epsilon}2}\right\|^{\frac 1H}_{L(\mathcal{H})}\right.\nonumber\\
&&\left.\left\|(A_h(s))^{-H+\frac{\epsilon}2}\left(U_h(s,t_{m-1})-I\right)(A_{h,m-1})^{\frac{1-\beta}2}\right\|^{\frac 1H}_{L(\mathcal{H})}ds\right)^{2H}\nonumber\\
&\leq& C\left\|(A(0))^{\frac{\beta-1}2}\phi(t_{m-1})\right\|^2_{L^0_2}\left(\int_{t_{m-1}}^{t_m}(t_m-s)^{-1+\frac{\epsilon}{2H}}(s-t_{m-1})^{\frac{2H+\beta-1-\epsilon}{2H}}ds\right)^{2H}\nonumber\\
&\leq& C\Delta t^{2H+\beta-1-\epsilon}\left(\int_{t_{m-1}}^{t_m}(t_m-s)^{-1+\frac{\epsilon}{2H}}ds\right)^{2H}\nonumber\\
&\leq& C\Delta t^{2H+\beta-1-\epsilon}\Delta t^{\epsilon}=C\Delta t^{2H+\beta-1}.
\end{eqnarray*}
Thanks to \eqref{fracint2}, \lemref{prepaSMTI1} (i) with $\alpha=\frac{1-\beta}2$, \eqref{equi3} and \assref{noiseandjump} \eqref{noise1}, we get
\begin{eqnarray}
\label{timerr*33}
&&\|V_{33}\|^2_{L^2(\Omega,\mathcal{H})}\nonumber\\
&=&\left\|\int_{t_{m-1}}^{t_m}\left(U_h(t_m,t_{m-1})-e^{-\Delta t A_{h,m-1}}\right)P_h\phi(t_{m-1})dB^H(s)\right\|^2_{L^2(\Omega,\mathcal{H})}\nonumber\\
&\leq& C\int_{t_{m-1}}^{t_m}\left\|\left(U_h(t_m,t_{m-1})-e^{-\Delta t A_{h,m-1}}\right)(A_{h,m-1})^{\frac{1-\beta}2}\right\|^2_{L(\mathcal{H})}\left\|(A_{h,m-1})^{\frac{\beta-1}2}P_h\phi(t_{m-1})\right\|^2_{L^0_2}ds\nonumber\\
&\leq& C\int_{t_{m-1}}^{t_m}\Delta t^{1+\beta}\left\|(A(0))^{\frac{\beta-1}2}\phi(t_{m-1})\right\|^2_{L^0_2}ds\nonumber\\
&\leq& C \Delta t^{2+\beta}\leq C \Delta t^{2H+\beta-1}.
\end{eqnarray}
Combining \eqref{timerr*31}-\eqref{timerr*33} and taking the square-root yields
\begin{eqnarray}
\label{timerr3}
\|V_3\|_{L^2(\Omega,\mathcal{H})}\leq C \Delta t^{\frac{2H+\beta-1}2}.
\end{eqnarray}
\subsubsection{Estimate of $V_4$}
 By adding and subtracting a term, using triangle inequality and the estimate $(a+b)^2\leq 2a^2+2b^2$, we get
\begin{eqnarray}
\label{timerr*4}
&&\|V_4\|^2_{L^2(\Omega,\mathcal{H})}\nonumber\\
&:=&\left\|\int_{t_{m-1}}^{t_m}\int_{\chi}U_h(t_m,s)P_h\psi(z)\tilde{N}(dz,ds)-\int_{t_{m-1}}^{t_m}\int_{\chi}e^{-\Delta t A_{h,m-1}}P_h\psi(z)\tilde{N}(dz,ds)\right\|^2_{L^2(\Omega,\mathcal{H})}\nonumber\\
&\leq& 4\left\|\int_{t_{m-1}}^{t_m}\int_{\chi}\left(U_h(t_m,s)-U_h(t_m,t_{m-1})\right)P_h\psi(z)\tilde{N}(dz,ds)\right\|^2_{L^2(\Omega,\mathcal{H})}\nonumber\\
&+&4\left\|\int_{t_{m-1}}^{t_m}\int_{\chi}\left(U_h(t_m,t_{m-1})-e^{-\Delta t A_{h,m-1}}\right)P_h\psi(z)\tilde{N}(dz,ds)\right\|^2_{L^2(\Omega,\mathcal{H})}\nonumber\\
&=:&4\|V_{41}\|^2_{L^2(\Omega,\mathcal{H})}+4\|V_{42}\|^2_{L^2(\Omega,\mathcal{H})}.
\end{eqnarray}
From the It\^{o} isometry \eqref{jumpint}, \lemref{disevol} (ii), \eqref{propdisevol5} with $\gamma=\frac{1-\epsilon}2$, \eqref{equi3}, \lemref{prepaSMTI1} (ii) with $\gamma=\frac{-1+\epsilon}2$, $\eta=-\frac{2H+\beta-2}2$ and \assref{noiseandjump} \eqref{jump},
\begin{eqnarray}
\label{timerr*41}
&&\|V_{41}\|^2_{L^2(\Omega,\mathcal{H})}\nonumber\\
%&=&\left\|\int_{t_{m-1}}^{t_m}\int_{\chi}\left(U_h(t_m,s)-U_h(t_m,t_{m-1})\right)P_h\psi(z)\tilde{N}(dz,ds)\right\|^2_{L^2(\Omega,\mathcal{H})}\nonumber\\
&=&\int_{t_{m-1}}^{t_m}\int_{\chi}\left\|\left(U_h(t_m,s)-U_h(t_m,t_{m-1})\right)P_h\psi(z)\right\|^2\nu(dz)ds\nonumber\\
&\leq&\int_{t_{m-1}}^{t_m}\int_{\chi}\left\|U_h(t_m,s)\left(I-U_h(s,t_{m-1})\right)(A_{h,m-1})^{-\frac{2H+\beta-2}2}\right\|^2_{L(\mathcal{H})}\nonumber\\
&&\times\left\|(A_{h,m-1})^{\frac{2H+\beta-2}2}P_h\psi(z)\right\|^2\nu(dz)ds\nonumber\\
&\leq&\int_{t_{m-1}}^{t_m}\int_{\chi}\left\|U_h(t_m,s)(A_h(s))^{\frac{1-\epsilon}2}\right\|^{\frac 1H}_{L(\mathcal{H})}\nonumber\\
&\times&\left\|(A_h(s))^{\frac{-1+\epsilon}2}\left(I-U_h(s,t_{m-1})\right)(A_{h,m-1})^{-\frac{2H+\beta-2}2}\right\|^2_{L(\mathcal{H})}\left\|(A(0))^{\frac{2H+\beta-2}2}\psi(z)\right\|^2\nu(dz)ds\nonumber\\
&\leq&C\left(\int_{t_{m-1}}^{t_m}(t_m-s)^{-1+\epsilon}(s-t_{m-1})^{2H+\beta-1-\epsilon}ds\right)\left(\int_{\chi}\left\|(A(0))^{\frac{2H+\beta-2}2}\psi(z)\right\|^2\nu(dz)\right)\nonumber\\
&\leq&C\Delta t^{2H+\beta-1-\epsilon}\left(\int_{t_{m-1}}^{t_m}(t_m-s)^{-1+\epsilon}ds\right)\leq C\Delta t^{2H+\beta-1-\epsilon}\Delta t^{\epsilon}=C\Delta t^{2H+\beta-1}.
\end{eqnarray}
Applying also the It\^{o} isometry property \eqref{jumpint}, inserting $(A_{h,m-1})^{-\frac{2H+\beta-2}2}(A_{h,m-1})^{\frac{2H+\beta-2}2}$, \lemref{prepaSMTI1} (i) with $\alpha=-\frac{2H+\beta-2}2$, \eqref{equi3} and \assref{noiseandjump} \eqref{jump} yields
\begin{eqnarray}
\label{timerr*42}
&&\|V_{42}\|^2_{L^2(\Omega,\mathcal{H})}\nonumber\\
&=&\left\|\int_{t_{m-1}}^{t_m}\int_{\chi}\left(U_h(t_m,t_{m-1})-e^{-\Delta t A_{h,m-1}}\right)P_h\psi(z)\tilde{N}(dz,ds)\right\|^2_{L^2(\Omega,\mathcal{H})}\nonumber\\
&\leq&\int_{t_{m-1}}^{t_m}\int_{\chi}\left\|\left(U_h(t_m,t_{m-1})-e^{-\Delta t A_{h,m-1}}\right)(A_{h,m-1})^{-\frac{2H+\beta-2}2}\right\|^2_{L(\mathcal{H})}\nonumber\\
&&\times\left\|(A_{h,m-1})^{\frac{2H+\beta-2}2}P_h\psi(z)\right\|^2\nu(dz)ds\nonumber\\
&\leq&C\left(\int_{t_{m-1}}^{t_m}\Delta t^{4-2H-\beta}ds\right)\left(\int_{\chi}\left\|(A(0))^{\frac{2H+\beta-2}2}\psi(z)\right\|^2\nu(dz)\right)\nonumber\\
&\leq&C\Delta t^{5-2H-\beta}\leq C\Delta t^{2H+\beta-1}.
\end{eqnarray} 
Putting \eqref{timerr*41} and \eqref{timerr*42} in \eqref{timerr*4} and taking the square-root gives
\begin{eqnarray}
\label{timerr4}
\|V_{42}\|_{L^2(\Omega,\mathcal{H})}\leq C\Delta t^{\frac{2H+\beta-1}2}.
\end{eqnarray}
\subsubsection{Estimate of $V_5$} 
As in \cite{Mukmaa}, we split it in five terms 
\begin{eqnarray}
\label{timerr*5}
V_5
&=&\sum_{k=1}^{m-1}\int_{t_{m-k-1}}^{t_{m-k}}\left(\prod_{j=m-k+1}^m U_h(t_j,t_{j-1})\right)U_h(t_{m-k},s)P_hF(s,X^h(s))ds\nonumber\\
&-&\sum_{k=1}^{m-1}\int_{t_{m-k-1}}^{t_{m-k}}\left(\prod_{j=m-k}^{m-1} e^{-\Delta t A_{h,j}}\right)e^{-\Delta t A_{h,m-k-1}}P_hF(t_{m-k-1},X^h_{m-k-1})ds\nonumber\\
&=&\sum_{k=1}^{m-1}\int_{t_{m-k-1}}^{t_{m-k}}\left(\prod_{j=m-k+1}^m U_h(t_j,t_{j-1})\right)U_h(t_{m-k},s)\nonumber\\
&&P_h\left[F(s,X^h(s))-F(t_{m-k-1},X^h(t_{m-k-1}))\right]ds\nonumber\\
&+&\sum_{k=1}^{m-1}\int_{t_{m-k-1}}^{t_{m-k}}\left(\prod_{j=m-k+1}^m U_h(t_j,t_{j-1})\right)\left[U_h(t_{m-k},s)-U_h(t_{m-k},t_{m-k-1})\right]\nonumber\\
&&P_hF(t_{m-k-1},X^h(t_{m-k-1}))ds\nonumber\\
&+&\sum_{k=1}^{m-1}\int_{t_{m-k-1}}^{t_{m-k}}\left[\left(\prod_{j=m-k}^m U_h(t_j,t_{j-1})\right)-\left(\prod_{j=m-k-1}^{m-1} e^{-\Delta t A_{h,j}}\right)\right]\nonumber\\
&&\left[U_h(t_{m-k},s)-U_h(t_{m-k},t_{m-k-1})\right]P_hF(t_{m-k-1},X^h(t_{m-k-1}))ds\nonumber\\
&+&\sum_{k=1}^{m-1}\int_{t_{m-k-1}}^{t_{m-k}}\left(\prod_{j=m-k}^{m-1}-e^{-\Delta t A_{h,j}}\right)\left(e^{-\Delta t A_{h,m-k-1}}-e^{-(t_{m-k}-s)A_{h,m-k-1}}\right)\nonumber\\
&&P_hF(t_{m-k-1},X^h(t_{m-k-1}))ds\nonumber\\
&+&\sum_{k=1}^{m-1}\int_{t_{m-k-1}}^{t_{m-k}}\left(\prod_{j=m-k}^{m-1}e^{-\Delta t A_{h,j}}\right)e^{-(t_{m-k}-s)A_{h,m-k-1}}\nonumber\\
&&P_h[F(t_{m-k-1},X^h(t_{m-k-1}))-F(t_{m-k-1},X^h_{m-k-1})]ds\nonumber\\
&=:&\sum_{i=1}^5V_{5i}.
\end{eqnarray}
Taking the $L^2$- norm and applying triangle inequality yields
\begin{eqnarray*}
\label{timerr**5}
\|V_5\|_{L^2(\Omega,\mathcal{H})}\leq \sum_{i=1}^5\|V_{5i}\|_{L^2(\Omega,\mathcal{H})}.
\end{eqnarray*}
As in \cite[(129),(130)]{Mukmaa}, we have 
\begin{eqnarray}
\label{timerr**52}
\sum_{i=2}^5\|V_{5i}\|_{L^2(\Omega,\mathcal{H})}\leq C\Delta t+C\Delta t\sum_{j=0}^{m-2}\|X^h(t_j)-X^h_j\|_{L^2(\Omega,\mathcal{H})}.
\end{eqnarray}
To estimate  $\|V_{51}\|_{L^2(\Omega,\mathcal{H})}$, by adding and subtracting a term gives
\begin{eqnarray}
\label{timerr*51}
V_{51}&=&\sum_{k=1}^{m-1}\int_{t_{m-k-1}}^{t_{m-k}}\left(\prod_{j=m-k+1}^m U_h(t_j,t_{j-1})\right)U_h(t_{m-k},s)\nonumber\\
&&P_h\left[F(s,X^h(s))-F(t_{m-k-1},X^h(t_{m-k-1}))\right]ds\nonumber\\
&=&\sum_{k=1}^{m-1}\int_{t_{m-k-1}}^{t_{m-k}}U_h(t_m,s)P_h\left[F(s,X^h(s))-F(t_{m-k-1},X^h(s))\right]ds\nonumber\\
&+&\sum_{k=1}^{m-1}\int_{t_{m-k-1}}^{t_{m-k}}U_h(t_m,s)P_h\left[F(t_{m-k-1},X^h(s))-F(t_{m-k-1},X^h(t_{m-k-1}))\right]ds\nonumber\\
&=:&V_{511}+V_{512}.
\end{eqnarray}
Using triangle inequality, boundedness of $P_h$ and $U_h(t,s)$, \eqref{propdisevol1}, \assref{nonlinear}, it holds that
\begin{eqnarray}
\label{timerr*511}
\|V_{511}\|_{L^2(\Omega,\mathcal{H})}\leq C\Delta t^{\frac{2H+\beta-1}2}.
\end{eqnarray}
As in \cite{manto,Tammab}, to achieve higher order in $V_{512}$, we apply Taylor expansion in the Banach space to $F$, we obtain for all $m=2,\cdot\cdot\cdot,M$ and $k=1,\cdot\cdot\cdot,m-1$
\begin{eqnarray}
\label{taylor}
F(t_{m-k-1},X^h(s))-F(t_{m-k-1},X^h(s))=I^h_{m-k-1,s}(X^h(s)-X^h(t_{m-k-1})),
\end{eqnarray} 
where $I^h_{m-k-1,s}$ is defined for $t_{m-k-1}\leq s\leq t_{m-k}$ as follows
\begin{eqnarray}
\label{not}
I^h_{m-k-1,s}:=\int_0^1F'\left(t_{m-k-1},X^h(t_{m-k-1})+\lambda(X^h(s)-X^h(t_{m-k-1}))\right)d\lambda.
\end{eqnarray}
Thanks to \assref{derivate} and \eqref{not}, we easily have
\begin{eqnarray}
\label{boundtayl}
\left\|I^h_{m-k-1,s}\right\|_{L(\mathcal{H})}\leq C,\hspace{1cm}m=2,\cdot\cdot\cdot,M,\hspace{0.5cm}k=1,\cdot\cdot\cdot,m-1,\hspace{0.5cm}t_{m-k-1}\leq s\leq t_{m-k}.
\end{eqnarray}
Note that $X^h(s)$, $t_{m-k-1}\leq s\leq t_{m-k}$ can be written as follows
\begin{eqnarray}
\label{mildsol3}
X^h(s)&=&U_h(s,t_{m-k-1})X^h(t_{m-k-1})+\int_{t_{m-k-1}}^{s}U_h(s,r)P_hF(r,X^h(r))dr\nonumber\\
&+&\int_{t_{m-k-1}}^sU_h(s,r)P_h\phi(r)dB^H(r)+\int_{t_{m-k-1}}^s\int_{\chi}U_h(s,r)P_h\psi(z)\tilde{N}(dz,dr).
\end{eqnarray}
Putting \eqref{mildsol3} in \eqref{taylor} yields
\begin{eqnarray*}
\label{diff}
&&F(t_{m-k-1},X^h(s))-F(t_{m-k-1},X^h(s))\\
&=& I^h_{m-k-1,s}\left(U_h(s,t_{m-k-1})-I\right)X^h(t_{m-k-1})+I^h_{m-k-1,s}\int_{t_{m-k-1}}^{s}U_h(s,r)P_hF(r,X^h(r))dr\nonumber\\
&+&I^h_{m-k-1,s}\int_{t_{m-k-1}}^sU_h(s,r)P_h\phi(r)dB^H(r)+I^h_{m-k-1,s}\int_{t_{m-k-1}}^s\int_{\chi}U_h(s,r)P_h\psi(z)\tilde{N}(dz,dr).\nonumber
\end{eqnarray*}
Hence
\begin{eqnarray}
\label{timerr*512}
V_{512}&=&\sum_{k=1}^{m-1}\int_{t_{m-k-1}}^{t_{m-k}}U_h(t_m,s)P_h\left[F(t_{m-k-1},X^h(s))-F(t_{m-k-1},X^h(s))\right]ds\nonumber\\
&=&\sum_{k=1}^{m-1}\int_{t_{m-k-1}}^{t_{m-k}}U_h(t_m,s)P_hI^h_{m-k-1,s}\left(U_h(s,t_{m-k-1})-I\right)X^h(t_{m-k-1})ds\nonumber\\
&+&\sum_{k=1}^{m-1}\int_{t_{m-k-1}}^{t_{m-k}}U_h(t_m,s)P_hI^h_{m-k-1,s}\int_{t_{m-k-1}}^{s}U_h(s,r)P_hF(r,X^h(r))drds\nonumber\\
&+&\sum_{k=1}^{m-1}\int_{t_{m-k-1}}^{t_{m-k}}U_h(t_m,s)P_hI^h_{m-k-1,s}\int_{t_{m-k-1}}^sU_h(s,r)P_h\phi(r)dB^H(r)ds\nonumber\\
&+&\sum_{k=1}^{m-1}\int_{t_{m-k-1}}^{t_{m-k}}U_h(t_m,s)P_hI^h_{m-k-1,s}\int_{t_{m-k-1}}^s\int_{\chi}U_h(s,r)P_h\psi(z)\tilde{N}(dz,dr)ds\nonumber\\
&=& \sum_{i=1}^{4}V_{512}^{(i)}.
\end{eqnarray}
Inserting an appropriate power of $A_{h,m-k-1}$, using triangle inequality, \eqref{propdisevol1}, boundedness of $P_h$, \eqref{boundtayl}, \eqref{propdisevol6} and \lemref{reg} (more precisely \eqref{spareg}) yields
\begin{eqnarray}
\label{timerr*5121}
&&\left\|V_{512}^{(1)}\right\|_{L^2(\Omega,\mathcal{H})}\nonumber\\
&\leq& C\sum_{k=1}^{m-1}\int_{t_{m-k-1}}^{t_{m-k}}\left\|U_h(t_m,s)P_hI^h_{m-k-1,s}\right\|_{L(\mathcal{H})}\left\|\left(U_h(s,t_{m-k-1})-I\right)(A_{h,m-k-1})^{-\frac{2H+\beta-1}2}\right\|_{L(\mathcal{H})}\nonumber\\
&&\times\left\|(A_{h,m-k-1})^{\frac{2H+\beta-1}2}X^h(t_{m-k-1})\right\|_{L^2(\Omega,\mathcal{H})}ds\nonumber\\
&\leq& C\sum_{k=1}^{m-1}\int_{t_{m-k-1}}^{t_{m-k}}(s-t_{m-k-1})^{\frac{2H+\beta-1}2}ds\leq C\Delta t^{\frac{2H+\beta-1}2}.
\end{eqnarray}
Using triangle inequality and thanks to \eqref{propdisevol1}, boundedness of $P_h$, \eqref{boundtayl}, \lemref{reg} \eqref{bound1}, it holds that
\begin{eqnarray*}
%\label{timerr*5122}
&&\left\|V_{512}^{(2)}\right\|_{L^2(\Omega,\mathcal{H})}\nonumber\\
&\leq& \sum_{k=1}^{m-1}\int_{t_{m-k-1}}^{t_{m-k}}\left\|U_h(t_m,s)P_hI^h_{m-k-1,s}\right\|_{L(\mathcal{H})}\int_{t_{m-k-1}}^{s}\left\|U_h(s,r)P_hF(r,X^h(r))\right\|_{L^2(\Omega,\mathcal{H})}drds\nonumber\\
&\leq& C\sum_{k=1}^{m-1}\int_{t_{m-k-1}}^{t_{m-k}}(s-t_{m-k-1})ds\leq C\Delta t.
\end{eqnarray*}
Using also triangle inequality, \eqref{propdisevol1}, boundedness of $P_h$, \eqref{boundtayl} and \eqref{timreg3} with $t_1=t_{m-k-1}$ and $t_2=s$ we obtain
\begin{eqnarray*}
\label{timerr*5123}
&&\left\|V_{512}^{(3)}\right\|_{L^2(\Omega,\mathcal{H})}\nonumber\\
&\leq& \sum_{k=1}^{m-1}\int_{t_{m-k-1}}^{t_{m-k}}\left\|U_h(t_m,s)P_hI^h_{m-k-1,s}\right\|_{L(\mathcal{H})}\left\|\int_{t_{m-k-1}}^{s}U_h(s,r)P_h\phi(r)dB^H(r)\right\|_{L^2(\Omega,\mathcal{H})}ds\nonumber\\
&\leq& C\sum_{k=1}^{m-1}\int_{t_{m-k-1}}^{t_{m-k}}(s-t_{m-k-1})^{\frac{2H+\beta-1}2}ds\leq C\Delta t^{\frac{2H+\beta-1}2}.
\end{eqnarray*}
We use triangle inequality to split $V_{512}^{(4)}$ in two terms as follows
\begin{eqnarray}
\label{timerr*5124}
&&V_{512}^{(4)}\nonumber\\
&=&\sum_{k=1}^{m-1}\int_{t_{m-k-1}}^{t_{m-k}}U_h(t_m,s)P_hI^h_{m-k-1,t_{m-k-1}}\int_{t_{m-k-1}}^s\int_{\chi}U_h(s,r)P_h\psi(z)\tilde{N}(dz,dr)ds\nonumber\\
&+&\sum_{k=1}^{m-1}\int_{t_{m-k-1}}^{t_{m-k}}U_h(t_m,s)P_h\left(I^h_{m-k-1,s}-I^h_{m-k-1, t_{m-k-1}}\right)\int_{t_{m-k-1}}^s\int_{\chi}U_h(s,r)P_h\psi(z)\tilde{N}(dz,dr)ds\nonumber\\
&=:&V_{512}^{(41)}+V_{512}^{(42)}.
\end{eqnarray}
Using the martingale property of stochastic integral, H\"older inequality, \eqref{propdisevol1}, boundedness of $P_h$, \eqref{boundtayl} and \eqref{timreg5} with $t_1=t_{m-k-1}$ and $t_2=s$ yields
\begin{eqnarray}
\label{timerr*51241}
&&\left\|V_{512}^{(41)}\right\|^2_{L^2(\Omega,\mathcal{H})}\nonumber\\
&=& \mathbb{E}\left\|\sum_{k=1}^{m-1}\int_{t_{m-k-1}}^{t_{m-k}}U_h(t_m,s)P_hI^h_{m-k-1,t_{m-k-1}}\int_{t_{m-k-1}}^s\int_{\chi}U_h(s,r)P_h\psi(z)\tilde{N}(dz,dr)ds\right\|^2\nonumber\\
&=& \sum_{k=1}^{m-1}\mathbb{E}\left\|\int_{t_{m-k-1}}^{t_{m-k}}U_h(t_m,s)P_hI^h_{m-k-1,t_{m-k-1}}\int_{t_{m-k-1}}^s\int_{\chi}U_h(s,r)P_h\psi(z)\tilde{N}(dz,dr)ds\right\|^2\nonumber\\
&\leq& \Delta t \sum_{k=1}^{m-1}\mathbb{E}\int_{t_{m-k-1}}^{t_{m-k}}\left\|U_h(t_m,s)P_hI^h_{m-k-1,t_{m-k-1}}\right\|^2_{L(\mathcal{H})}\nonumber\\
&&\times\left\|\int_{t_{m-k-1}}^s\int_{\chi}U_h(s,r)P_h\psi(z)\tilde{N}(dz,dr)\right\|^2_{L^2(\Omega,\mathcal{H})}ds\nonumber\\
&\leq& C\Delta t \sum_{k=1}^{m-1}\int_{t_{m-k-1}}^{t_{m-k}}(s-t_{m-k-1})^{\min(1,2H+\beta-1)}ds\nonumber\\
&\leq& C\Delta t \Delta t^{\min(1,2H+\beta-1)}\leq C\Delta t^{2H+\beta-1}.
\end{eqnarray}
Using triangle inequality and H\"{o}lder inequality and It\^{o} isometry yields
\begin{eqnarray}
\label{timerr*51242}
\Vert V_{512}^{(42)}\Vert_{L^2(\Omega, H)}&\leq& \sum_{k=1}^{m-1}\left\Vert\int_{t_{m-k-1}}^{t_{m-k}}U_h(t_m,s)P_h\left(I^h_{m-k-1,s}-I^h_{m-k-1, t_{m-k-1}}\right)\right.\nonumber\\
&&\left.\int_{t_{m-k-1}}^s\int_{\chi}U_h(s,r)P_h\psi(z)\tilde{N}(dz,dr)ds\right\Vert_{L^2(\Omega, H)}\nonumber\\
&\leq& \Delta t^{\frac{1}{2}}\sum_{k=1}^{m-1}\left[\int_{t_{m-k-1}}^{t_{m-k}}\mathbb{E}\left\Vert \int_{t_{m-k-1}}^s U_h(t_m,s)P_h\left(I^h_{m-k-1,s}-I^h_{m-k-1, t_{m-k-1}}\right)\right.\right.\nonumber\\
&&\left.\left.\int_{\chi}U_h(s,r)P_h\psi(z)\tilde{N}(dz,dr)\right\Vert^2ds\right]^{\frac{1}{2}}\nonumber\\
&\leq& \Delta t^{\frac{1}{2}}\sum_{k=1}^{m-1}\left[\int_{t_{m-k-1}}^{t_{m-k}} \int_{t_{m-k-1}}^s\int_{\chi}\mathbb{E}\left\Vert U_h(t_m,s)P_h\left(I^h_{m-k-1,s}-I^h_{m-k-1, t_{m-k-1}}\right)\right.\right.\nonumber\\
&&\left.\left.U_h(s,r)P_h\psi(z)\right\Vert^2\nu(dz)drds\right]^{\frac{1}{2}}
\end{eqnarray}
Using  \assref{noiseandjump} yields
\begin{eqnarray}
\label{dan1}
&&\int_{\chi}\mathbb{E}\left\Vert U_h(t_m,s)P_h\left(I^h_{m-k-1,s}-I^h_{m-k-1, t_{m-k-1}}\right)U_h(s,r)P_h\psi(z)\right\Vert^2\nu(dz)\nonumber\\
&\leq& \int_{\chi}\mathbb{E}\left\Vert U_h(t_m,t_{m-k})(A_h(0))^{\eta}\right\Vert^2_{L(\mathcal{H})}\left\Vert (A_h(0))^{-\eta}U_h(t_{m-k},s)(A_h(0))^{\eta}\right\Vert^2_{L(\mathcal{H})}\nonumber\\
&&\times\left\Vert(A_h(0))^{-\eta} P_h\left(I^h_{m-k-1,s}-I^h_{m-k-1, t_{m-k-1}}\right)U_h(s,r)P_h\psi(z)\right\Vert^2\nu(dz)\nonumber\\
&\leq& Ct_k^{-2\eta}\left\Vert(A_h(0))^{-\eta} P_h\left(I^h_{m-k-1,s}-I^h_{m-k-1, t_{m-k-1}}\right)\right\Vert^2_{L(\mathcal{H})}\nonumber\\
&&\times\int_{\chi}\left\Vert U_h(s,r)(A_h(0))^{-\frac{2H+\beta-2}{2}}\right\Vert^2_{L(\mathcal{H}}\left\Vert (A_h(0))^{\frac{2H+\beta-2}{2}} P_h\psi(z)\right\Vert^2\nu(dz)\nonumber\\
&\leq& Ct_k^{-2\eta}(s-r)^{\min(0, 2H+\beta-2)}\left\Vert(A_h(0))^{-\eta} P_h\left(I^h_{m-k-1,s}-I^h_{m-k-1, t_{m-k-1}}\right)\right\Vert^2_{L(\mathcal{H})}.
\end{eqnarray}
 From the definition of $I^h_{m-k-1,s}$, by using \assref{derivate} and \lemref{reg} we arrive at 
 \begin{eqnarray}
 \label{dan2}
 &&\Vert (-A_h(s))^{-\eta}\left(I^h_{m-k-1,s}-I^h_{m-k1,t_{m-k-1}}\right)\Vert_{L(\mathcal{H})}\nonumber\\
&\leq& \int_0^1\left\Vert (-A_h(s))^{-\eta}P_h\left( F'\left(t_{m-k-1},X^h(t_{m-k-1})+\lambda\left(X^h(s)-X^h(t_{m-k-1})\right)\right)\right.\right.\nonumber\\
&&\left.\left.-F'\left(t_{m-k-1}, X^h(t_{m-k-1})\right)\right)\right\Vert_{L(\mathcal{H})}d\lambda\nonumber\\
&\leq& \int_0^1\left\Vert (-A(s))^{-\frac{\eta}{2}}\left( F'\left(t_{m-k-1},X^h(t_{m-k-1})+\lambda\left(X^h(s)-X^h(t_{m-k-1})\right)\right)\right.\right.\nonumber\\
&&\left.\left.-F'\left(t_{m-k-1}, X^h(t_{m-k-1})\right)\right)\right\Vert_{L(\mathcal{H})}d\lambda\nonumber\\
&\leq& C\int_0^1\lambda\Vert X^h(s)-X^h(t_{m-k-1})\Vert d\lambda\nonumber\\
&\leq & C\Vert X^h(s)-X^h(t_{m-k-1})\Vert\leq C(s-t_{m-k-1})^{\min(1, 2H+\beta-1)/2}.
 \end{eqnarray}
Substituting \eqref{dan2} in \eqref{dan1} yields
\begin{eqnarray}
\label{dan3}
&&\int_{\chi}\mathbb{E}\left\Vert U_h(t_m,s)P_h\left(I^h_{m-k-1,s}-I^h_{m-k-1, t_{m-k-1}}\right)U_h(s,r)P_h\psi(z)\right\Vert^2\nu(dz)\nonumber\\
&\leq& Ct_k^{-2\eta}(s-r)^{\min(0, 2H+\beta-2)}(s-t_{m-k-1})^{\min(1, 2H+\beta-1)}.
\end{eqnarray}
Substituting \eqref{dan3} in \eqref{timerr*51242} yields
\begin{eqnarray}
\label{timerr**51242}
&&\Vert V_{512}^{(42)}\Vert_{L^2(\Omega, \mathcal{H})}\nonumber\\
&\leq& C\Delta t^{\frac{1}{2}}\sum_{k=1}^{m-1}\left[t_k^{-2\eta}\int_{t_{m-k-1}}^{t_{m-k}}\int_{t_{m-k-1}}^s(s-r)^{\min(0, 2H+\beta-2)}(s-t_{m-k-1})^{\min(1, 2H+\beta-1)}drds\right]^{\frac{1}{2}}\nonumber\\
&\leq& C\Delta t^{\frac{1}{2}}\sum_{k=1}^{m-1}t_k^{-\eta}\left[\int_{t_{m-k-1}}^{t_{m-k}}\int_{t_{m-k-1}}^s(s-t_{m-k-1})^{\min(1, 4H+2\beta-3)}drds\right]^{\frac{1}{2}}\nonumber\\
&\leq& C\Delta t^{\frac{2H+\beta-1}{2}}. 
\end{eqnarray}
Substituting  \eqref{timerr**51242}--\eqref{timerr*512}  and \eqref{timerr*511} in \eqref{timerr*51}, then combining with \eqref{timerr**52} yields
\begin{eqnarray}
\label{timerr5}
\|V_{5}\|_{L^2(\Omega,\mathcal{H})}&\leq& C\Delta t^{\frac{2H+\beta-1}2}+C\Delta t+C\Delta t\sum_{j=0}^{m-2}\|X^h(t_j)-X^h_j\|_{L^2(\Omega,\mathcal{H})}\nonumber\\
&\leq& C\Delta t^{\frac{2H+\beta-1}2}+C\Delta t\sum_{j=0}^{m-2}\|X^h(t_j)-X^h_j\|_{L^2(\Omega,\mathcal{H})}.
\end{eqnarray}
\subsubsection{Estimate of $V_6$}
The  $V_{6}$ can we rewrite as follows
\begin{eqnarray}
\label{timerr*6}
V_6&=&\int_0^{t_{m-1}}\left(\prod_{j=\llcorner s\lrcorner^m+2}^m U_h(t_j,t_{j-1})\right)U_h(\llcorner s\lrcorner+\Delta t,s)P_h\phi(s)dB^H(s)\nonumber\\
&-&\int_0^{t_{m-1}}\left(\prod_{j=\llcorner s\lrcorner^m+1}^{m-1} e^{-\Delta t A_{h,j}}\right)e^{-\Delta t A_{h,\llcorner s\lrcorner^m}}P_h\phi(\llcorner s\lrcorner)dB^H(s),
\end{eqnarray}
where the notation $\llcorner s\lrcorner$ and $\llcorner s\lrcorner^m$ are defined by
\begin{eqnarray}
\label{not2}
\llcorner s\lrcorner:=\left[\frac{s}{\Delta t}\right]\Delta t \quad\text{and}\quad\llcorner s\lrcorner^m:=\left[\frac{s}{\Delta t}\right],
\end{eqnarray}
and we split \eqref{timerr*6} into four terms as follows
{\small
\begin{eqnarray}
\label{timerr**6}
&&V_6\nonumber\\
&=&\int_0^{t_{m-1}}\left(\prod_{j=\llcorner s\lrcorner^m+2}^m U_h(t_j,t_{j-1})\right)\left(U_h(\llcorner s\lrcorner+\Delta t,s)-U_h(\llcorner s\lrcorner+\Delta t,\llcorner s\lrcorner)\right)\nonumber\\
&&P_h(\phi(s)-\phi(t_{m-1}))dB^H(s)\nonumber\\
&+&\int_0^{t_{m-1}}\left(\prod_{j=\llcorner s\lrcorner^m+2}^m U_h(t_j,t_{j-1})\right)\left(U_h(\llcorner s\lrcorner+\Delta t,s)-U_h(\llcorner s\lrcorner+\Delta t,\llcorner s\lrcorner)\right)P_h\phi(t_{m-1})dB^H(s)\nonumber\\
&+&\int_0^{t_{m-1}}\left[\left(\prod_{j=\llcorner s\lrcorner^m+1}^m U_h(t_j,t_{j-1})\right)-\left(\prod_{j=\llcorner s\lrcorner^m}^{m-1} e^{-\Delta t A_{h,j}}\right)\right]P_h\phi(s)dB^H(s)\nonumber\\
&+&\int_0^{t_{m-1}}\left(\prod_{j=\llcorner s\lrcorner^m}^{m-1} e^{-\Delta t A_{h,j}}\right)P_h(\phi(s)-\phi(\llcorner s\lrcorner))dB^H(s)\nonumber\\
&=:&\sum_{i=1}^4V_{6i}.
\end{eqnarray}
}
Taking the $L^2$-norm, squaring both sides and  using the estimate $\left(\sum_{i=1}^na_i\right)^2\leq n^2\sum_{i=1}^na_i^2$, we have
\begin{eqnarray}
\label{timerr_6}
\left\|V_6\right\|^2_{L^2(\Omega,\mathcal{H})}\leq 16\sum_{i=1}^4\left\|V_{6i}\right\|^2_{L^2(\Omega,\mathcal{H})}.
\end{eqnarray} 
Using \eqref{fracint2}, \lemref{disevol} (iii), inserting an appropriate power of $A_h(s)$ and $A_{h,\llcorner s\lrcorner^m}$, \eqref{propdisevol5} with $\gamma=H$, \lemref{prepaSMTI1} (ii) with $\gamma=H$ and $\eta=\frac{1-\beta}2$, \lemref{equi} \eqref{equi3} with $\alpha=\frac{\beta-1}2$ and \assref{noiseandjump} \eqref{noise2} yields
\begin{eqnarray}
\label{timerr_61}
&&\left\|V_{61}\right\|^2_{L^2(\Omega,\mathcal{H})}\nonumber\\
&=&\|\int_0^{t_{m-1}}\left(\prod_{j=\llcorner s\lrcorner^m+2}^m U_h(t_j,t_{j-1})\right)\left(U_h(\llcorner s\lrcorner+\Delta t,s)-U_h(\llcorner s\lrcorner+\Delta t,\llcorner s\lrcorner)\right)\nonumber\\
&&P_h(\phi(s)-\phi(t_{m-1}))dB^H(s)\|^2_{L^2(\Omega,\mathcal{H})}\nonumber\\
&\leq& C\int_0^{t_{m-1}}\left|U_h(t_m,s)(A_h(s))^H\right|^2_{L(\mathcal(H))}\left\|(A_h(s))^{-H}\left(I-U_h(s,\llcorner s\lrcorner)(A_{h,\llcorner s\lrcorner^m})^{\frac{1-\beta}2}\right)\right\|^2_{L(\mathcal{H})} \nonumber\\
&&\left\| (A_{h,\llcorner s\lrcorner^m})^{\frac{\beta-1}2}P_h(\phi(s)-\phi(t_{m-1}))\right\|^2_{L^0_2}ds\nonumber\\
&\leq& C\int_0^{t_{m-1}}(t_m-s)^{-2H}(s-\llcorner s\lrcorner)^{2H+\beta-1}\left\| (A(0))^{\frac{\beta-1}2}(\phi(t_{m-1})-\phi(s))\right\|^2_{L^0_2}ds\nonumber\\
&\leq& C\Delta t^{2H+\beta-1}\int_0^{t_{m-1}}(t_m-s)^{-2H}(t_{m-1}-s)^{2\delta}ds\nonumber\\
&\leq& C\Delta t^{2H+\beta-1}\int_0^{t_{m-1}}(t_{m-1}-s)^{-2H+2\delta}ds\leq C\Delta t^{2H+\beta-1}t_{m-1}^{-2H+2\delta+1}\leq C\Delta t^{2H+\beta-1}.
\end{eqnarray}
By \eqref{fracint1}, \lemref{disevol} (iii), inserting an appropriate power of $A_h(s)$ and $A_{h,\llcorner s\lrcorner^m}$, \eqref{propdisevol1}, \lemref{prepaSMTI1} (ii) with $\gamma=H$ and $\eta=\frac{1-\beta}2$, \lemref{equi} \eqref{equi3} with $\alpha=\frac{\beta-1}2$, \eqref{inte4} with $\rho=0$, $\gamma=2H$ and \assref{noiseandjump} \eqref{noise1}, it holds that.
\begin{eqnarray*}
%\label{timerr_62}
&&\left\|V_{62}\right\|^2_{L^2(\Omega,\mathcal{H})}\nonumber\\
&=&\left\|\int_0^{t_{m-1}}\left(\prod_{j=\llcorner s\lrcorner^m+2}^m U_h(t_j,t_{j-1})\right)\right.\nonumber\\
&&\left.\left(U_h(\llcorner s\lrcorner+\Delta t,s)-U_h(\llcorner s\lrcorner+\Delta t,\llcorner s\lrcorner)\right)P_h\phi(t_{m-1})dB^H(s)\right\|^2_{L^2(\Omega,\mathcal{H})}\nonumber\\
&\leq& C\sum_{i\in\N^d}(\int_0^{t_{m-1}}\left|U_h(t_m,s)(A_h(s))^H\right|^{\frac 1H}_{L(\mathcal(H))}\left\|(A_h(s))^{-H}\left(I-U_h(s,\llcorner s\lrcorner)(A_{h,\llcorner s\lrcorner^m})^{\frac{1-\beta}2}\right)\right\|^{\frac 1H}_{L(\mathcal{H})} \nonumber\\
&&\left\|(A_{h,\llcorner s\lrcorner^m})^{\frac{\beta-1}2}P_h\phi(t_{m-1})Q^{\frac 12}e_i\right\|^{\frac 1H}ds)^{2H}\nonumber\\
&\leq& C\left(\int_0^{t_{m-1}}\left\|U_h(t_m,t_{m-1})\right\|^{\frac 1H}_{L(\mathcal{H})}\left\|U_h(t_{m-1},s)(A_h(s))^H\right\|^{\frac 1H}_{L(\mathcal{H})}(s-\llcorner s\lrcorner)^{\frac{2H+\beta-1}{2H}}ds\right)^{2H} \nonumber\\
&&\left\|(A(0))^{\frac{\beta-1}2}\phi(t_{m-1})\right\|^2_{L^0_2}\nonumber\\
&\leq& C\Delta t^{2H+\beta-1}\left(\int_0^{t_{m-1}}\left\|U_h(t_{m-1},s)(A_h(s))^H\right\|^{\frac 1H}_{L(\mathcal{H})}ds\right)^{2H} \leq C\Delta t^{2H+\beta-1}.
\end{eqnarray*}
Using \eqref{fracint2}, inserting an appropriate power of $A_{h,\llcorner s\lrcorner^m}$, \lemref{prepaSMTI2} (iii) with $\alpha=\frac{1-\beta}2$, \eqref{equi3} with $\alpha=\frac{\beta-1}2$, \eqref{noise1}, the variable change $j=m-k$ and \cite[(169)]{Muka} yields
\begin{eqnarray*}
%\label{timerr_63}
&&\left\|V_{63}\right\|^2_{L^2(\Omega,\mathcal{H})}\nonumber\\
&=&\left\|\int_0^{t_{m-1}}\left[\left(\prod_{j=\llcorner s\lrcorner^m+1}^m U_h(t_j,t_{j-1})\right)-\left(\prod_{j=\llcorner s\lrcorner^m}^{m-1} e^{-\Delta t A_{h,j}}\right)\right]P_h\phi(s)dB^H(s)\right\|^2_{L^2(\Omega,\mathcal{H})}\nonumber\\
&\leq& C\int_0^{t_{m-1}}\left\|\left[\left(\prod_{j=\llcorner s\lrcorner^m+1}^m U_h(t_j,t_{j-1})\right)-\left(\prod_{j=\llcorner s\lrcorner^m}^{m-1} e^{-\Delta t A_{h,j}}\right)\right](A_{h,\llcorner s\lrcorner^m})^{\frac{1-\beta}2}\right\|^2_{L(\mathcal{H})} \nonumber\\
&&\left\| (A_{h,\llcorner s\lrcorner^m})^{\frac{\beta-1}2}P_h\phi(s)\right\|^2_{L^0_2}ds\nonumber\\
&\leq& C\int_0^{t_{m-1}}\Delta t^{1+\beta-\epsilon} t_{m-\llcorner s\lrcorner^m}^{-1+\beta+\epsilon}\left\| (A(0))^{\frac{\beta-1}2}\phi(s)\right\|^2_{L^0_2}ds\nonumber\\
&\leq& C\Delta t^{1+\beta-\epsilon}\sum_{k=0}^{m-2}\int_{t_k}^{t_{k+1}} t_{m-k}^{-1+\beta+\epsilon}ds\nonumber\\
&\leq& C\Delta t^{1+\beta-\epsilon}\left(\sum_{j=2}^{m} t_j^{-1+\beta+\epsilon}\Delta t\right)\leq C\Delta t^{1+\beta-\epsilon}.
\end{eqnarray*}
At last $\eqref{fracint2}$, inserting  , inserting an appropriate power of $A_{h,\llcorner s\lrcorner^m}$, \lemref{prepaSMTI3} with $\gamma=\frac{1-\beta}2$, \eqref{equi3} with $\alpha=\frac{\beta-1}2$, \eqref{noise2}, the variable change $j=m-k$ and \cite[(169)]{Muka} yields
\begin{eqnarray*}
&&\left\|V_{64}\right\|^2_{L^2(\Omega,\mathcal{H})}\nonumber\\
&=&\left\|\int_0^{t_{m-1}}\left(\prod_{j=\llcorner s\lrcorner^m}^{m-1} e^{-\Delta t A_{h,j}}\right)P_h(\phi(s)-\phi(\llcorner s\lrcorner))dB^H(s)\right\|^2_{L^2(\Omega,\mathcal{H})}\nonumber\\
&\leq& C\int_0^{t_{m-1}}\left\|\left(\prod_{j=\llcorner s\lrcorner^m}^{m-1} e^{-\Delta t A_{h,j}}\right)(A_{h,\llcorner s\lrcorner^m})^{\frac{1-\beta}2}\right\|^2_{L(\mathcal{H})}\left\| (A_{h,\llcorner s\lrcorner^m})^{\frac{\beta-1}2}P_h(\phi(s)-\phi(\llcorner s\lrcorner))\right\|^2_{L^0_2}ds\nonumber\\
&\leq& C\int_0^{t_{m-1}} t_{m-\llcorner s\lrcorner^m}^{-1+\beta}\left\| (A(0))^{\frac{\beta-1}2}(\phi(s)-\phi(\llcorner s\lrcorner))\right\|^2_{L^0_2}ds\nonumber\\
&\leq& C\int_0^{t_{m-1}} t_{m-\llcorner s\lrcorner^m}^{-1+\beta}(s-\llcorner s\lrcorner)^{2\delta}ds.
\end{eqnarray*}
Classical estimates yields
\begin{eqnarray}
\label{timerr_64}
\left\|V_{64}\right\|^2_{L^2(\Omega,\mathcal{H})}
&\leq& C\Delta t^{2H+\beta-1}\int_0^{t_{m-1}} t_{m-\llcorner s\lrcorner^m}^{-1+\beta}ds\nonumber\\
&\leq& C\Delta t^{2H+\beta-1}\sum_{k=0}^{m-2}\int_{t_k}^{t_{k+1}} t_{m-k}^{-1+\beta}ds\nonumber\\
&\leq& C\Delta t^{2H+\beta-1}\left(\sum_{j=2}^{m} t_j^{-1+\beta}\Delta t\right)\leq C\Delta t^{2H+\beta-1}.
\end{eqnarray}
Combining \eqref{timerr_61}-\eqref{timerr_64} related to \eqref{timerr_6} and taking the square-root yields
\begin{eqnarray}
\label{timerr6}
\left\|V_6\right\|_{L^2(\Omega,\mathcal{H})}\leq C\Delta t^{\frac{2H+\beta-1}2}+C\Delta t^{\frac{1+\beta-\epsilon}2}\leq C\Delta t^{\frac{2H+\beta-1}2}.
\end{eqnarray}
\subsubsection{Estimate of $V_7$}
 By adding and subtracting a term in $V_7$, we split it as
\begin{eqnarray}
\label{timerr*7}
V_7&=&\sum_{k=1}^{m-1}\int_{t_{m-k-1}}^{t_{m-k}}\int_{\chi}\left(\prod_{j=m-k+1}^m U_h(t_j,t_{j-1})\right)U_h(t_{m-k},s)P_hz_0\tilde{N}(dz,ds)\nonumber\\
&-&\sum_{k=1}^{m-1}\int_{t_{m-k-1}}^{t_{m-k}}\int_{\chi}\left(\prod_{j=m-k}^{m-1} e^{-\Delta t A_{h,j}}\right)e^{-\Delta t A_{h,m-k-1}}P_hz_0\tilde{N}(dz,ds)\nonumber\\
&=&\sum_{k=1}^{m-1}\int_{t_{m-k-1}}^{t_{m-k}}\int_{\chi}\left(\prod_{j=m-k+1}^m U_h(t_j,t_{j-1})\right)\left[U_h(t_{m-k},s)-U_h(t_{m-k},t_{m-k-1})\right]P_hz_0\tilde{N}(dz,ds)\nonumber\\
&+&\sum_{k=1}^{m-1}\int_{t_{m-k-1}}^{t_{m-k}}\int_{\chi}\left[\left(\prod_{j=m-k}^m U_h(t_j,t_{j-1})\right)-\left(\prod_{j=m-k-1}^{m-1} e^{-\Delta t A_{h,j}}\right)\right] P_h\psi(z)\tilde{N}(dz,ds)\nonumber\\
&=&V_{71}+V_{72}.
\end{eqnarray}
Using the martingale property of stochastic integral, the It\^{o} isometry \eqref{jumpint}, inserting an appropriate power of $A_h(s)$ and $A_{h,m-k-1}$ yields
\begin{eqnarray*}
&&\left\|V_{71}\right\|^2_{L^2(\Omega;\mathcal{H})}\nonumber\\
&=&\mathbb{E}\left\|\sum_{k=1}^{m-1}\int_{t_{m-k-1}}^{t_{m-k}}\int_{\chi}\left(\prod_{j=m-k+1}^m U_h(t_j,t_{j-1})\right)\right.\nonumber\\
&&\left.\left[U_h(t_{m-k},s)-U_h(t_{m-k},t_{m-k-1})\right]P_h\psi(z)\tilde{N}(dz,ds)\right\|^2\nonumber\\
&=&\sum_{k=1}^{m-1}\int_{t_{m-k-1}}^{t_{m-k}}\int_{\chi}\left\|\left(\prod_{j=m-k+1}^m U_h(t_j,t_{j-1})\right)U_h(t_{m-k},s)(I-U_h(s,t_{m-k-1}))P_h\psi(z)\right\|^2\nu(dz)ds\nonumber\\
&\leq&\sum_{k=1}^{m-1}\int_{t_{m-k-1}}^{t_{m-k}}\left\|U_h(t_m,t_{m-1})\right\|^2_{L(\mathcal{H})}\left\|U_h(t_{m-1},s)(A_h(s))^{\frac 12}\right\|^2_{L(\mathcal{H})}\nonumber\\
&&\times\left\|(A_h(s))^{-\frac 12}(I-U_h(s,t_{m-k-1}))(A_{h,m-k-1})^{\frac {2-2H-\beta}2}\right\|^2_{L(\mathcal{H})}ds\nonumber\\
&&\times\left(\int_{\chi}\left\| (A_{h,m-k-1})^{\frac {2H+\beta-2}2}P_h\psi(z)\right\|^2\nu(dz)\right).
\end{eqnarray*}
Using \eqref{propdisevol1}, \eqref{propdisevol5} with $\gamma=\frac 12$, \lemref{prepaSMTI1}  with $\gamma=\frac 12$ and $\eta=-\frac{2H+\beta-2}2$, \eqref{equi3} with $\alpha=\frac{2H+\beta-2}2$, \assref{noiseandjump} \eqref{jump} and \eqref{inte1} with $\rho=1$, it holds that
\begin{eqnarray}
\label{timerr*71}
&&\left\|V_{71}\right\|^2_{L^2(\Omega;\mathcal{H})}\nonumber\\
&\leq&C\sum_{k=1}^{m-1}\int_{t_{m-k-1}}^{t_{m-k}}\left\|U_h(t_m,t_{m-1})\right\|^2_{L(\mathcal{H})}\left\|U_h(t_{m-1},s)(A_h(s))^{\frac 12}\right\|^2_{L(\mathcal{H})}(s-t_{m-k-1})^{\frac {2H+\beta-1}2}ds\nonumber\\
&&\times\left(\int_{\chi}\left\| (A(0))^{\frac {2H+\beta-2}2}\psi(z)\right\|^2\nu(dz)\right)\nonumber\\
%&\leq&C\sum_{k=1}^{m-1}\int_{t_{m-k-1}}^{t_{m-k}}\left\|U_h(t_{m-1},s)(A_h(s))^{\frac 12}\right\|^2_{L(\mathcal{H})}(s-t_{m-k-1})^{\frac {2H+\beta-1}2}ds\nonumber\\
&\leq&C\Delta t^{\frac {2H+\beta-1}2}\int_0^{t_{m-1}}\left\|U_h(t_{m-1},s)(A_h(s))^{\frac 12}\right\|^2_{L(\mathcal{H})}ds\nonumber\\
&\leq&C\Delta t^{\frac {2H+\beta-1}2}.
\end{eqnarray}
Using again the martingale property of stochastic integral, It\^{o} isometry \eqref{jumpint}, inserting an appropriate power of $A_{h,m-k-1}$ yields
\begin{eqnarray*}
&&\left\|V_{72}\right\|^2_{L^2(\Omega;\mathcal{H})}\nonumber\\
&=&\mathbb{E}\left\|\sum_{k=1}^{m-1}\int_{t_{m-k-1}}^{t_{m-k}}\int_{\chi}\left[\left(\prod_{j=m-k}^m U_h(t_j,t_{j-1})\right)-\left(\prod_{j=m-k-1}^{m-1} e^{-\Delta t A_{h,j}}\right)\right] P_h\psi(z)\tilde{N}(dz,ds)\right\|^2\nonumber\\
&=&\sum_{k=1}^{m-1}\int_{t_{m-k-1}}^{t_{m-k}}\int_{\chi}\left\|\left[\left(\prod_{j=m-k}^m U_h(t_j,t_{j-1})\right)-\left(\prod_{j=m-k-1}^{m-1} e^{-\Delta t A_{h,j}}\right)\right](A_{h,m-k-1})^{-\frac{2H+\beta-2}2}\right\|^2_{L(\mathcal{H})}\nonumber\\
&&\left\|(A_{h,m-k-1})^{\frac{2H+\beta-2}2} P_h\psi(z)\right\|^2\nu(dz)ds.
\end{eqnarray*}
Using \lemref{prepaSMTI2} ((i) if $2H+\beta=2$, (ii) if $2H+\beta>2$ and (iii) if $2H+\beta<2$), \eqref{equi3} with $\alpha=\frac{2H+\beta-2}2$, \assref{noiseandjump} \eqref{jump} and \cite[(169)]{Muka} yields
\begin{eqnarray}
\label{timerr*72}
\left\|V_{72}\right\|^2_{L^2(\Omega;\mathcal{H})}
&\leq& C\sum_{k=1}^{m-1}\int_{t_{m-k-1}}^{t_{m-k}}\int_{\chi}\Delta t^{2H+\beta-\epsilon}t_{k+1}^{2-2H-\beta+\epsilon}\left\|(A(0))^{\frac{2H+\beta-2}2}\psi(z)\right\|^2\nu(dz)ds\nonumber\\
&\leq& C\Delta t^{2H+\beta-\epsilon}\left(\sum_{k=1}^{m-1}\int_{t_{m-k-1}}^{t_{m-k}}t_{k+1}^{2-2H-\beta+\epsilon}ds\right)\left(\int_{\chi}\left\|(A(0))^{\frac{2H+\beta-2}2}\psi(z)\right\|^2\nu(dz)\right)\nonumber\\
&\leq& C\Delta t^{2H+\beta-\epsilon}\left(\sum_{k=1}^{m-1}\Delta t \hspace{0.1cm}t_{k+1}^{-1+\epsilon}t_{k+1}^{3-2H-\beta}\right)\nonumber\\
&\leq& C\Delta t^{2H+\beta-\epsilon}\left(\sum_{k=1}^{m-1} t_{k+1}^{-1+\epsilon}\Delta t\right)\nonumber\\
&\leq& C\Delta t^{2H+\beta-\epsilon}\leq C\Delta t^{2H+\beta-1}.
\end{eqnarray}
Combining \eqref{timerr*71}, \eqref{timerr*72} and taking the square-root yields
\begin{eqnarray}
\label{timerr7}
\left\|V_7\right\|_{L^2(\Omega;\mathcal{H})}\leq C\Delta t^{\frac {2H+\beta-1}2}.
\end{eqnarray}
Adding the estimates \eqref{timerr1}, \eqref{timerr2}, \eqref{timerr3}, \eqref{timerr4}, \eqref{timerr5}, \eqref{timerr6} and \eqref{timerr7} yields
\begin{equation*}
%\label{Timerr}
\left\|X^h(t_m)-X^h_m\right\|_{L^2(\Omega;\mathcal{H})}\leq C\Delta t^{\frac {2H+\beta-1}2}+C\Delta t\sum_{j=0}^{m-1}\left\|X^h(t_j)-X^h_j\right\|_{L^2(\Omega;\mathcal{H})}.
\end{equation*}
The desired result follows from the discrete Gronwall's inequality added to \eqref{boundspaerr}. $\hfill\square$

\section{Proof of \thmref{mainthm} for Magnus-implicit scheme}
\label{convimpl}
It is important to mention that the estimates made in this section are inspired by the results in \cite{Mukmab} when the non-autonomous problem is only driven by a standard Brownian motion. But in our case as we know, it's driven by a cylindrical fractional Brownian motion and a random Poisson measure. To reach our goal, we present some preparatory results.
\subsection{Preparatory results}
\begin{lem}
\label{prepaimpl1}
For any $\alpha_1,\alpha_2\in[0,1)$ and $0\leq j,k\leq M$. The following estimates hold
\begin{eqnarray}
&&\left\|(A_{h,k})^{-\alpha_1}\left(e^{-\Delta t A_{h,j}}-S^j_{h,\Delta t}\right)(A_{h,j})^{-\alpha_2}\right\|_{L(\mathcal{H})}\leq C\Delta t^{\alpha_1+\alpha_2},\label{prepaimpl11}\\
&&\left\|(A_{h,k})^{\alpha_1}\left(e^{-\Delta t A_{h,j}}-S^j_{h,\Delta t}\right)(A_{h,j})^{-\alpha_2}\right\|_{L(\mathcal{H})}\leq C\Delta t^{-\alpha_1+\alpha_2},\label{prepaimpl12}\\
&&\left\|(A_{h,k})^{-\alpha_1}\left(e^{-\Delta t A_{h,j}}-S^j_{h,\Delta t}\right)(A_{h,j})^{\alpha_2}\right\|_{L(\mathcal{H})}\leq C\Delta t^{\alpha_1-\alpha_2},\label{prepaimpl13}
\end{eqnarray}
\end{lem}
\textit{Proof:} See \cite[Lemma 3.9]{Mukmab} for the proof of \eqref{prepaimpl11} and \eqref{prepaimpl12}. As in \cite[(69)]{Mukmab} we set
\begin{eqnarray*}
%\label{prepaimpl131}
K^j_{h,\Delta t}:=e^{-\Delta t A_{h,j}}-S^j_{h,\Delta t}.
\end{eqnarray*} 
So using \cite[(72)]{Mukmab}, we get 
\begin{eqnarray}
\label{prepaimpl132}
-K^j_{h,\Delta t}&=&\int_0^{\Delta t}sA^2_{h,j}(I+sA_{h,j})^{-2}e^{(\Delta t-s)A_{h,j}}ds\\
&=&\int_0^{\Delta t}sA^{2-\epsilon}_{h,j}(I+sA_{h,j})^{-2}A^{\epsilon}_{h,j}e^{(\Delta t-s)A_{h,j}}ds.
\end{eqnarray}
Using \lemref{equi}, we have
\begin{eqnarray}
\label{prepaimpl133}
\left\|(A_{h,k})^{-\alpha_1}K^j_{h,\Delta t}(A_{h,j})^{\alpha_2}\right\|_{L(\mathcal{H})}\leq C\left\|(A_{h,j})^{-\alpha_1}K^j_{h,\Delta t}(A_{h,j})^{\alpha_2}\right\|_{L(\mathcal{H})}.
\end{eqnarray}
Using \eqref{prepaimpl132},\eqref{prepaimpl133}, \cite[Lemma 3.7]{Mukmab} and \eqref{semigroup3}, it holds that
\begin{eqnarray*}
%\label{prepaimpl134}
\left\|(A_{h,k})^{-\alpha_1}K^j_{h,\Delta t}(A_{h,j})^{\alpha_2}\right\|_{L(\mathcal{H})}&\leq&
\int_0^{\Delta t}s\left\|A^{2-\alpha_1-\epsilon}_{h,j}(I+sA_{h,j})^{-2}\right\|_{L(\mathcal{H})}\left\| A^{\alpha_2+\epsilon}_{h,j}e^{(\Delta t-s)A_{h,j}}\right\|_{L(\mathcal{H})}ds\nonumber\\
&\leq& C\int_0^{\Delta t}ss^{-2+\alpha_1+\epsilon}(\Delta t-s)^{-\alpha_2-\epsilon}ds\nonumber\\
&\leq& C\int_0^{\Delta t}s^{-1+\alpha_1+\epsilon}(\Delta t-s)^{-\alpha_2-\epsilon}ds\nonumber\\
&\leq& C\Delta t^{\alpha_1-\alpha_2}\mathcal{B}(\alpha_1+\epsilon,1-\alpha_2-\epsilon)\nonumber\\
&\leq& C\Delta t^{\alpha_1-\alpha_2}.
\end{eqnarray*}
The proof of \lemref{prepaimpl1} is thus completed. $\hfill\square$
\begin{lem}(\cite[Lemma 3.8]{Mukmab})
\label{prepaimpl2}
Let \assref{linear} be fulfilled.
\begin{enumerate}
\item [(i)] For any $\alpha\in[0,1)$, it holds that
\begin{eqnarray}
\label{prepaimpl21}
\left\|(A_{h,k})^{\alpha}\left(\prod_{j=i}^mS^j_{h,\Delta t}\right)\right\|_{L(\mathcal{H})}\leq C t_{m-i+1}^{-\alpha} ,\quad 0\leq i\leq m\leq M,\quad 0\leq k\leq M.
\end{eqnarray}
\item [(ii)] For any $\alpha_1,\alpha_2\in[0,1)$, it holds that
\begin{eqnarray*}
%\label{prepaimpl22}
\left\|(A_{h,k})^{\alpha_1}\left(\prod_{j=i}^mS^j_{h,\Delta t}\right)(A_{h,i})^{-\alpha_2}\right\|_{L(\mathcal{H})}\leq C t_{m-i+1}^{-\alpha_1+\alpha_2},\quad 0\leq i\leq m\leq M,\quad 0\leq k\leq M.
\end{eqnarray*}
\end{enumerate}
\end{lem}
\begin{lem}
\label{prepaimpl3}
For all $\alpha\in[0,1)$, there exists $C_{\alpha}\geq 0$ such that
\begin{eqnarray*}
%\label{prepaimpl31}
\Delta t\sum_{j=1}^m t_{m-j+1}^{-\alpha}t_j^{-\alpha}\leq C_{\alpha} t_m^{1-2\alpha}.
\end{eqnarray*} 
\end{lem}
\textit{Proof:} This estimate is an application of the second estimate in \cite[(76)]{Mukmab} with $\alpha_2=1-\alpha$.
\begin{lem}
\label{prepaimpl4}
Let $0\leq \alpha<2$ and let \assref{linear} be fulfilled.
\begin{enumerate}
\item [(i)] If $v\in\mathcal{D}\left((A(0))^{\frac{\alpha}2}\right)$, then the following estimate holds
\begin{eqnarray*}
\left\|\left(\prod_{j=i}^me^{-\Delta tA_{h,j}}\right)P_h v-\left(\prod_{j=i}^mS^j_{h,\Delta t}\right) P_h v\right\|_{L(\mathcal{H})} \leq C\Delta t^{\frac{\alpha}2}\left\|v\right\|_{\alpha},\quad 1\leq i\leq m\leq M.\nonumber
\end{eqnarray*}
\item [(ii)] Moreover, for non smooth data, ie for $v\in \mathcal{H}$, it holds that 
\begin{eqnarray*}
\left\|\left(\prod_{j=i}^me^{-\Delta tA_{h,j}}\right)P_h v-\left(\prod_{j=i}^mS^j_{h,\Delta t}\right) P_h v\right\|_{L(\mathcal{H})} \leq C\Delta t^{\frac{\alpha}2}t_{m-i}^{-\frac{\alpha}2}\left\|v\right\|,\quad 1\leq i< m\leq M.\nonumber
\end{eqnarray*}
\item [(iii)] For any $\alpha_1,\alpha_2\in[0,1)$, such that $\alpha_1\leq \alpha_2$, it holds that 
\begin{eqnarray*}
\left\|\left[\left(\prod_{j=i}^me^{-\Delta tA_{h,j}}\right) -\left(\prod_{j=i}^mS^j_{h,\Delta t}\right)\right](A_{h,i})^{\alpha_1-\alpha_2}\right\|_{L(\mathcal{H})} \leq C\Delta t^{\alpha_2}t_{m-i}^{-\alpha_1},\quad 1\leq i< m\leq M.\nonumber
\end{eqnarray*}
\item [(iv)] For any $\gamma\in[0,1)$, it holds that 
\begin{eqnarray*}
\left\|\left[\left(\prod_{j=i}^me^{-\Delta tA_{h,j}}\right) -\left(\prod_{j=i}^mS^j_{h,\Delta t}\right)\right](A_{h,i})^{\frac{\gamma}2}\right\|_{L(\mathcal{H})} \leq C\Delta t^{\frac{1-\gamma-\epsilon}2}t_{m-i}^{\frac{-1+\epsilon}2},\quad 1\leq i< m\leq M.\nonumber
\end{eqnarray*}
\item [(v)] For any $\gamma\in[0,1)$ and $H\in(\frac 12, 1)$, it holds that 
\begin{eqnarray*}
\left\|\left[\left(\prod_{j=i}^me^{-\Delta tA_{h,j}}\right) -\left(\prod_{j=i}^mS^j_{h,\Delta t}\right)\right](A_{h,i})^{\frac{\gamma}2}\right\|_{L(\mathcal{H})} \leq C\Delta t^{\frac{2H-\gamma}2}t_{m-i}^{-H},\quad 1\leq i< m\leq M.\nonumber
\end{eqnarray*}
\end{enumerate}
\end{lem}
\textit{Proof:} See \cite[Lemma 3.11]{Mukmab} for the proof of (i)--(iv). Concerning the proof of (v), using the telescopic sum, inserting an appropriate power of $A_{h,m}$
, $A_{h,i+1}$ and $A_{h,i+k}$, taking the norm and using triangle inequality yields
\begin{eqnarray*}
%\label{prepaimpl41}
&&\left\|\left[\left(\prod_{j=i}^me^{-\Delta tA_{h,j}}\right) -\left(\prod_{j=i}^mS^j_{h,\Delta t}\right)\right](A_{h,i})^{\frac{\gamma}2}\right\|_{L(\mathcal{H})}\nonumber\\
&\leq& \left\|\left( e^{-\Delta tA_{h,m}}-S^m_{h,\Delta t}\right)(A_{h,m})^{-\frac{2H-\gamma}2}(A_{h,m})^{\frac{2H-\gamma}2}\left(\prod_{j=i}^{m-1}S^j_{h,\Delta t}\right)(A_{h,i})^{\frac{\gamma}2}\right\|_{L(\mathcal{H})}\nonumber\\
&+&\left\|\left(\prod_{j=i+1}^me^{-\Delta tA_{h,j}}\right)(A_{h,i+1})^{H}(A_{h,i+1})^{-H}\left( e^{-\Delta tA_{h,i}}-S^i_{h,\Delta t}\right)(A_{h,i})^{\frac{\gamma}2}\right\|_{L(\mathcal{H})}\nonumber\\
&+& \sum_{k=2}^{m-i}\left\|\left(\prod_{j=i+k}^me^{-\Delta tA_{h,j}}\right)(A_{h,i+k})^{H}(A_{h,i+k})^{-H}\left( e^{-\Delta tA_{h,i+k-1}}-S^{i+k-1}_{h,\Delta t}\right)(A_{h,i+k-1})^{-\frac{2H-\gamma}2}\right.\nonumber\\
&&\left.(A_{h,i+k-1})^{\frac{2H-\gamma}2}\left(\prod_{j=i}^{i+k-2}S^j_{h,\Delta t}\right)(A_{h,i})^{\frac{\gamma}2}\right\|_{L(\mathcal{H})}\nonumber\\
&=:&K_1+K_2+K_3.
\end{eqnarray*}
Using \eqref{prepaimpl11} with $\alpha_1=0$, $\alpha_2=\frac{2H-\gamma}{2}$ and \eqref{prepaimpl21} with $\alpha=H$, we have
\begin{eqnarray}
\label{prepaimpl42}
K_1&\leq& \left\|\left( e^{-\Delta tA_{h,m}}-S^m_{h,\Delta t}\right)(A_{h,m})^{-\frac{2H-\gamma}2}\right\|_{L(\mathcal{H})} \left\|(A_{h,m})^{\frac{2H-\gamma}2}\left(\prod_{j=i}^{m-1}S^j_{h,\Delta t}\right)(A_{h,i})^{\frac{\gamma}2}\right\|_{L(\mathcal{H})}\nonumber\\
&\leq& C\Delta t^{\frac{2H-\gamma}2}t_{m-i}^{-H}.
\end{eqnarray}
Using \lemref{prepaSMTI3}, with $\gamma=H$ and \lemref{prepaimpl1} \eqref{prepaimpl13} with $\alpha_1=H$, $\alpha_2=\frac{\gamma}2$ yields
\begin{eqnarray*}
%\label{prepaimpl43}
K_2&\leq& \left\|\left(\prod_{j=i+1}^me^{-\Delta tA_{h,j}}\right)(A_{h,i+1})^{H}\right\|_{L(\mathcal{H})} \left\|(A_{h,i+1})^{-H}\left( e^{-\Delta tA_{h,i}}-S^i_{h,\Delta t}\right)(A_{h,i})^{\frac{\gamma}2}\right\|_{L(\mathcal{H})}\nonumber\\
&\leq& C\Delta t^{\frac{2H-\gamma}2}t_{m-i}^{-H}.
\end{eqnarray*}
Applying \lemref{prepaSMTI3} with $\gamma=\frac{1+H}2$, $\alpha_2=\frac{1+H-\gamma}2$, \eqref{prepaimpl21} with $\alpha=\frac{1+H}2$ and \lemref{prepaimpl3} with $\alpha=\frac{1+H}2$, it holds that
\begin{eqnarray}
\label{prepaimpl44}
K_3&\leq& \sum_{k=2}^{m-i}\left\|\left(\prod_{j=i+k}^me^{-\Delta tA_{h,j}}\right)(A_{h,i+k})^{H}\right\|_{L(\mathcal{H})}\nonumber\\
&\times&\left\|(A_{h,i+k})^{-H}\left( e^{-\Delta tA_{h,i+k-1}}-S^{i+k-1}_{h,\Delta t}\right)(A_{h,i+k-1})^{-\frac{2H-\gamma}2}\right\|_{L(\mathcal{H})} \nonumber\\
&\times&\left\|(A_{h,i+k-1})^{\frac{2H-\gamma}2}\left(\prod_{j=i}^{i+k-2}S^j_{h,\Delta t}\right)(A_{h,i})^{\frac{\gamma}2}\right\|_{L(\mathcal{H})}\nonumber\\
&\leq& C\sum_{k=2}^{m-i}t_{m+1-i-k}^{-\frac{1+H}2}\Delta t^{1+\frac{2H-\gamma}2}t_{k-1}^{-\frac{1+H}2}\nonumber\\
&\leq& C\Delta t^{\frac{2H-\gamma}2}\left(\Delta t\sum_{k=2}^{m-i}t_{m+1-i-k}^{-\frac{1+H}2}t_{k-1}^{-\frac{1+H}2}\right)\leq C\Delta t^{\frac{2H-\gamma}2}t_{m-i}^{-H}.
\end{eqnarray}
Adding \eqref{prepaimpl42}-\eqref{prepaimpl44} completes the proof of $(v)$.$\hfill\square$

With theses lemmas, we are now ready to prove \thmref{mainthm} for the Magnus-implicit scheme.
\subsection{Main proof of \thmref{mainthm} for semi-implicit scheme}
 Using triangle inequality yields
\begin{eqnarray*}
%\label{proofimpl*}
\left\|X(t_m)-Y^h_m\right\|_{L^2(\Omega;\mathcal{H})}\leq \left\|X(t_m)-X^h(t_m)\right\|_{L^2(\Omega;\mathcal{H})}+\left\|X^h(t_m)-Y^h_m\right\|_{L^2(\Omega;\mathcal{H})}.
\end{eqnarray*} 
The space error is estimated in \lemref{spaerr}. It remains to estimate the time error. We recall that the mild solution \eqref{soldispb} can be rewritten as follows
\begin{eqnarray}
\label{itmildsol}
X^h(t_m)&=&\left(\prod_{j=1}^m U_h(t_j,t_{j-1})\right)P_hX_0+\int_{t_{m-1}}^{t_m}U_h(t_m,s)P_hF(s,X^h(s))ds\nonumber\\
&+&\int_{t_{m-1}}^{t_m}U_h(t_m,s)P_h\phi(s)dB^H(s)+\int_{t_{m-1}}^{t_m}\int_{\chi}U_h(t_m,s)P_h\psi(z)\tilde{N}(dz,ds)\nonumber\\
&+&\sum_{k=1}^{m-1}\int_{t_{m-k-1}}^{t_{m-k}}\left(\prod_{j=m-k+1}^m U_h(t_j,t_{j-1})\right)U_h(t_{m-k},s)P_hF(s,X^h(s))ds\nonumber\\
&+&\sum_{k=1}^{m-1}\int_{t_{m-k-1}}^{t_{m-k}}\left(\prod_{j=m-k+1}^m U_h(t_j,t_{j-1})\right)U_h(t_{m-k},s)P_h\phi(s)dB^H(s)\nonumber\\
&+&\sum_{k=1}^{m-1}\int_{t_{m-k-1}}^{t_{m-k}}\int_{\chi}\left(\prod_{j=m-k+1}^m U_h(t_j,t_{j-1})\right)U_h(t_{m-k},s)P_h\psi(z)\tilde{N}(dz,ds).
\end{eqnarray}
Iterating the numerical solution \eqref{impl} as in \cite{Mukmab} yields
\begin{eqnarray}
\label{itimpl}
Y^h_m&=&\left(\prod_{j=0}^{m-1} S^j_{h,\Delta t}\right)P_h X_0+\int_{t_{m-1}}^{t_m}S^{m-1}_{h,\Delta t}P_hF(t_{m-1},Y^h_{m-1})ds+\int_{t_{m-1}}^{t_m}S^{m-1}_{h,\Delta t}P_h\phi(t_{m-1})dB^H(s)\nonumber\\
&+&\int_{t_{m-1}}^{t_m}\int_{\chi}S^{m-1}_{h,\Delta t}P_hz_0\tilde{N}(dz,ds)+\sum_{i=2}^{m}\int_{t_{m-i}}^{t_{m-i+1}}\left(\prod_{j=m-i}^{m-1} S^j_{h,\Delta t}\right)P_hF(t_{m-i},Y^h_{m-i})ds\nonumber\\
&+&\sum_{i=2}^{m}\int_{t_{m-i}}^{t_{m-i+1}}\left(\prod_{j=m-i}^{m-1} S^j_{h,\Delta t}\right)P_h\phi(t_{m-i})dB^H(s)\nonumber\\
&+&\sum_{i=2}^{m}\int_{t_{m-i}}^{t_{m-i+1}}\int_{\chi}\left(\prod_{j=m-i}^{m-1}S^j_{h,\Delta t}\right)P_h\psi(z)\tilde{N}(dz,ds).
\end{eqnarray}
Subtracting \eqref{itmildsol} to \eqref{itimpl}, taking the $L^2$-norm and using triangle inequality, we deduce
%\begin{eqnarray}
%X^h(t_m)-Y^h_m&=&\left(\prod_{j=1}^m U_h(t_j,t_{j-1})\right)P_hX_0-\left(\prod_{j=0}^{m-1} S^j_{h,\Delta t}\right)P_hX_0\nonumber\\
%&+&\int_{t_{m-1}}^{t_m}\left[U_h(t_m,s)P_hF(s,X^h(s))-S^{m-1}_{h,\Delta t}P_hF(t_{m-1},Y^h_{m-1})\right]ds\nonumber\\
%&+&\int_{t_{m-1}}^{t_m}\left[U_h(t_m,s)P_h\phi(s)-S^{m-1}_{h,\Delta t}P_h\phi(t_{m-1})\right]dB^H(s)\nonumber\\
%&+&\int_{t_{m-1}}^{t_m}\int_{\chi}\left[U_h(t_m,s)-S^{m-1}_{h,\Delta t}\right]P_hz_0\tilde{N}(dz,ds)\nonumber\\
%&+&\sum_{i=2}^{m}\int_{t_{m-i}}^{t_{m-i+1}}\left(\prod_{j=m-i+2}^m U_h(t_j,t_{j-1})\right)U_h(t_{m-i+1},s)P_hF(s,X^h(s))ds\nonumber\\
%&-&\sum_{i=2}^{m}\int_{t_{m-i}}^{t_{m-i+1}}\left(\prod_{j=m-i}^{m-1} S^j_{h,\Delta t}\right)P_hF(t_{m-i},Y^h_{m-i})ds\nonumber\\
%&+&\sum_{i=2}^{m}\int_{t_{m-i}}^{t_{m-i+1}}\left(\prod_{j=m-i+2}^m U_h(t_j,t_{j-1})\right)U_h(t_{m-i+1},s)P_h\phi(s)dB^H(s)\nonumber\\
%&-&\sum_{i=2}^{m}\int_{t_{m-i}}^{t_{m-i+1}}\left(\prod_{j=m-i}^{m-1} S^j_{h,\Delta t}\right)P_h\phi(t_{m-i})dB^H(s)\nonumber\\
%&+&\sum_{i=2}^{m}\int_{t_{m-i}}^{t_{m-i+1}}\int_{\chi}\left(\prod_{j=m-i+2}^m U_h(t_j,t_{j-1})\right)U_h(t_{m-i+1},s)P_hz_0\tilde{N}(dz,ds)\nonumber\\
%&-&\sum_{i=2}^{m}\int_{t_{m-i}}^{t_{m-i+1}}\int_{\chi}\left(\prod_{j=m-i}^{m-1}S^j_{h,\Delta t}\right)P_hz_0\tilde{N}(dz,ds)\nonumber\\
%&=:&\sum_{i=1}^7 VI_i.
%\end{eqnarray}
%Using triangle inequality yields 
\begin{eqnarray*}
%\label{proofimpl}
\|X^h(t_m)-Y^h_m\|_{L^2(\Omega,\mathcal{H})}\leq \sum_{i=1}^7 \|VI_i\|_{L^2(\Omega,\mathcal{H})}.
\end{eqnarray*}
We estimate $\|VI_i\|_{L^2(\Omega,\mathcal{H})}$, $i\in\{	1,2,\cdot\cdot\cdot,7\}$ separately. 
\subsubsection{Estimate of $V_1$, $V_2$, $V_3$, $V_4$ and $V_5$}
By adding and subtracting a term, applying triangle inequality, \lemref{prepaimpl4} (i) with $\alpha=\frac{2H+\beta-1}2$ and \assref{initdata}, it holds that
\begin{eqnarray}
\label{proofimpl1}
\|VI_1\|_{L^2(\Omega,\mathcal{H})}&:=&\left\|\left(\prod_{j=1}^m U_h(t_j,t_{j-1})\right)P_hX_0-\left(\prod_{j=0}^{m-1} S^j_{h,\Delta t}\right)P_hX_0\right\|_{L^2(\Omega,\mathcal{H})}\nonumber\\
&\leq& \left\|\left(\prod_{j=1}^m U_h(t_j,t_{j-1})\right)P_hX_0-\left(\prod_{j=0}^{m-1} e^{-\Delta t A_{h,j}}\right)P_hX_0\right\|_{L^2(\Omega,\mathcal{H})}\nonumber\\
&+&\left\|\left(\prod_{j=0}^{m-1} e^{-\Delta t A_{h,j}}\right)P_hX_0-\left(\prod_{j=0}^{m-1} S^j_{h,\Delta t}\right)P_hX_0\right\|_{L^2(\Omega,\mathcal{H})}\nonumber\\
&\leq& C\Delta t+C\Delta t^{\frac{2H+\beta-1}2}\|(A(0))^{\frac{2H+\beta-1}2}X_0\|_{L^2(\Omega,\mathcal{H})}\leq C\Delta t^{\frac{2H+\beta-1}2}.
\end{eqnarray}
By a similar reasoning as \eqref{timerr2}, \eqref{timerr*5}--\eqref{timerr5} in the previous section and \cite[(104), (135)--(153)]{Mukmab}, we obtain
\begin{eqnarray}
\label{proofimpl2}
&&\|VI_2\|_{L^2(\Omega,\mathcal{H})}+\|VI_5\|_{L^2(\Omega,\mathcal{H})}\nonumber\\
&:=&\int_{t_{m-1}}^{t_m}\left[U_h(t_m,s)P_hF(s,X^h(s))-S^{m-1}_{h,\Delta t}P_hF(t_{m-1},Y^h_{m-1})\right]ds\nonumber\\
&+&\sum_{i=2}^{m}\int_{t_{m-i}}^{t_{m-i+1}}\left(\prod_{j=m-i+2}^m U_h(t_j,t_{j-1})\right)U_h(t_{m-i+1},s)P_hF(s,X^h(s))ds\nonumber\\
&-&\sum_{i=2}^{m}\int_{t_{m-i}}^{t_{m-i+1}}\left(\prod_{j=m-i}^{m-1} S^j_{h,\Delta t}\right)P_hF(t_{m-i},Y^h_{m-i})ds\nonumber\\
&\leq& C\Delta t^{\frac{2H+\beta-1}2}+C\Delta t\sum_{i=1}^{m-1}\|X^h(t_i)-Y^h_i\|_{L^2(\Omega,\mathcal{H})}.
\end{eqnarray}
For the estimate of $\|VI_3\|_{L^2(\Omega,\mathcal{H})}$, by adding and subtracting a term, using triangle inequality and the estimate $(a+b)^2\leq 2a^2+2b^2$, we recast it in two terms as follows
\begin{eqnarray}
\label{proofimpl3*}
\|VI_3\|^2_{L^2(\Omega,\mathcal{H})}&:=&\left\|\int_{t_{m-1}}^{t_m}\left[U_h(t_m,s)P_h\phi(s)-S^{m-1}_{h,\Delta t}P_h\phi(t_{m-1})\right]dB^H(s)\right\|^2_{L^2(\Omega,\mathcal{H})}\nonumber\\
&\leq& 4\left\|\int_{t_{m-1}}^{t_m}\left[U_h(t_m,s)P_h\phi(s)-e^{-\Delta t A_{h,m-1}}P_h\phi(t_{m-1})\right]dB^H(s)\right\|^2_{L^2(\Omega,\mathcal{H})}\nonumber\\
&+& 4\left\|\int_{t_{m-1}}^{t_m}\left(e^{-\Delta t A_{h,m-1}}-S^{m-1}_{h,\Delta t}\right)P_h\phi(t_{m-1})dB^H(s)\right\|^2_{L^2(\Omega,\mathcal{H})}\nonumber\\
&=:&4\|VI_{31}\|^2_{L^2(\Omega,\mathcal{H})}+4\|VI_{32}\|^2_{L^2(\Omega,\mathcal{H})}.
\end{eqnarray} 
The first term $VI_{31}$ is the same that $V_3$, hence thanks \eqref{timerr3} we have
\begin{eqnarray}
\label{proofimpl3*1}
\|VI_{31}\|^2_{L^2(\Omega,\mathcal{H})}\leq C\Delta t^{2H+\beta-1}.
\end{eqnarray}
Using \eqref{fracint1}, inserting an appropriate power of $A_{h,m-1}$, \eqref{prepaimpl13} with $\alpha_1=0$
, $\alpha_2=\frac{1-\beta}2$, \eqref{equi3} with $\alpha=\frac{\beta-1}2$ and \assref{noiseandjump} \eqref{noise1}, it holds that
\begin{eqnarray}
\label{proofimpl3*2}
&&\|VI_{32}\|^2_{L^2(\Omega,\mathcal{H})}\nonumber\\
&=&\left\|\int_{t_{m-1}}^{t_m}\left(e^{-\Delta t A_{h,m-1}}-S^{m-1}_{h,\Delta t}\right)P_h\phi(t_{m-1})dB^H(s)\right\|^2_{L^2(\Omega,\mathcal{H})}\nonumber\\
&\leq& C\sum_{i\in\N^d}\left(\int_{t_{m-1}}^{t_m}\left\|\left(e^{-\Delta t A_{h,m-1}}-S^{m-1}_{h,\Delta t}\right)P_h\phi(t_{m-1})Q^{\frac 12}e_i\right\|^{\frac 1H}ds\right)^{2H}\nonumber\\
&\leq& C\left(\int_{t_{m-1}}^{t_m}\left\|\left(e^{-\Delta t A_{h,m-1}}-S^{m-1}_{h,\Delta t}\right)(A_{h,m-1})^{\frac{1-\beta}2}\right\|^{\frac 1H}ds\right)^{2H}\left\| (A_{h,m-1})^{\frac{\beta-1}2}P_h\phi(t_{m-1})\right\|^2_{L^0_2}\nonumber\\
&\leq& C\left(\int_{t_{m-1}}^{t_m}\Delta t^{\frac{-1+\beta}{2H}}ds\right)^{2H}\left\| (A(0))^{\frac{\beta-1}2}\phi(t_{m-1})\right\|^2_{L^0_2}\nonumber\\
&\leq& C\Delta t^{2H+\beta-1}.
\end{eqnarray}
Putting \eqref{proofimpl3*1} and \eqref{proofimpl3*2} in \eqref{proofimpl3*} and taking the square-root yields
\begin{eqnarray}
\label{proofimpl3}
\|VI_3\|_{L^2(\Omega,\mathcal{H})}\leq C\Delta t^{\frac{2H+\beta-1}2}.
\end{eqnarray}
By adding and subtracting a term, we split $VI_4$ in two following terms
\begin{eqnarray*}
%\label{proofimpl4*}
VI_4&:=&\int_{t_{m-1}}^{t_m}\int_{\chi}\left[U_h(t_m,s)-S^{m-1}_{h,\Delta t}\right]P_h\psi(z)\tilde{N}(dz,ds)\nonumber\\
&=&\int_{t_{m-1}}^{t_m}\int_{\chi}\left[U_h(t_m,s)-e^{-\Delta t A_{h,m-1}}\right]P_h\psi(z)\tilde{N}(dz,ds)\nonumber\\
&+&\int_{t_{m-1}}^{t_m}\int_{\chi}\left[e^{-\Delta t A_{h,m-1}}-S^{m-1}_{h,\Delta t}\right]P_h\psi(z)\tilde{N}(dz,ds)\nonumber\\
&=:&VI_{41}+VI_{42}.
\end{eqnarray*}
By the estimate done in \eqref{timerr*4}--\eqref{timerr4}, we have
\begin{eqnarray}
\label{proofimpl4*1}
\left\|VI_{41}\right\|_{L^2(\Omega;\mathcal{H})}\leq C\Delta t^{\frac{2H+\beta-1}2}.
\end{eqnarray}
Using It\^{o} isometry \eqref{jumpint}, inserting an appropriate power of $A_{h,m-1}$, \lemref{prepaimpl1} \eqref{prepaimpl11} with $\alpha_1=0$, $\alpha_2=\frac{2H+\beta-2}2$ if $2H+\beta>2$ and \eqref{prepaimpl13} with $\alpha_1=0$, $\alpha_2=-\frac{2H+\beta-2}2$ if $2H+\beta\leq 2$, \eqref{equi3} with $\alpha=\frac{2H+\beta-2}2$, \assref{noiseandjump} and \eqref{jump}, it holds that
\begin{eqnarray}
\label{proofimpl4*2}
&&\left\|VI_{42}\right\|^2_{L^2(\Omega;\mathcal{H})}\nonumber\\
&=&\left\| \int_{t_{m-1}}^{t_m}\int_{\chi}\left[e^{-\Delta t A_{h,m-1}}-S^{m-1}_{h,\Delta t}\right]P_h\psi(z)\tilde{N}(dz,ds)\right\|^2_{L^2(\Omega;\mathcal{H})}\nonumber\\
&\leq&  \int_{t_{m-1}}^{t_m}\int_{\chi}\left\|\left[e^{-\Delta t A_{h,m-1}}-S^{m-1}_{h,\Delta t}\right](A_{h,m-1})^{-\frac{2H+\beta-2}2}\right\|^2\left\| (A_{h,m-1})^{\frac{2H+\beta-2}2}P_h\psi(z)\right\|^2\nu(dz)ds\nonumber\\
&\leq&  C\int_{t_{m-1}}^{t_m}\int_{\chi}\Delta t^{2H+\beta-2}\left\| (A(0))^{\frac{2H+\beta-2}2}\psi(z)\right\|^2\nu(dz)ds\nonumber\\
&\leq&  C\Delta t^{2H+\beta-2}\left(\int_{t_{m-1}}^{t_m}ds\right)\left(\int_{\chi} \left\| (A(0))^{\frac{2H+\beta-2}2}\psi(z)\right\|^2\nu(dz)\right)\nonumber\\
&\leq& C\Delta t^{2H+\beta-1}.
\end{eqnarray}
Adding \eqref{proofimpl4*1} and \eqref{proofimpl4*2} gets 
 \begin{eqnarray}
 \label{proofimpl4}
 \left\|VI_4\right\|_{L^2(\Omega;\mathcal{H})}\leq C\Delta t^{\frac{2H+\beta-1}2}.
 \end{eqnarray}
 \subsubsection{Estimate of $V_6$}
By adding and subtracting a term in $V_6$, it holds that
\begin{eqnarray*}
%\label{proofimpl6*}
VI_6&:=&\sum_{i=2}^{m}\int_{t_{m-i}}^{t_{m-i+1}}\left(\prod_{j=m-i+2}^m U_h(t_j,t_{j-1})\right)U_h(t_{m-i+1},s)P_h\phi(s)dB^H(s)\nonumber\\
&-&\sum_{i=2}^{m}\int_{t_{m-i}}^{t_{m-i+1}}\left(\prod_{j=m-i}^{m-1} S^j_{h,\Delta t}\right)P_h\phi(t_{m-i})dB^H(s)\nonumber\\
&=&\sum_{i=2}^{m}\int_{t_{m-i}}^{t_{m-i+1}}\left(\prod_{j=m-i+2}^m U_h(t_j,t_{j-1})\right)U_h(t_{m-i+1},s)P_h\phi(s)dB^H(s)\nonumber\\
&-&\sum_{i=2}^{m}\int_{t_{m-i}}^{t_{m-i+1}}\left(\prod_{j=m-i}^{m-1} e^{-\Delta t A_{h,j}}\right)P_h\phi(t_{m-i})dB^H(s)\nonumber\\
&+&\sum_{i=2}^{m}\int_{t_{m-i}}^{t_{m-i+1}}\left[\left(\prod_{j=m-i}^{m-1} e^{-\Delta t A_{h,j}}\right)-\left(\prod_{j=m-i}^{m-1} S^j_{h,\Delta t}\right)\right]P_h\phi(t_{m-i})dB^H(s)\nonumber\\
&=:&VI_{61}+VI_{62}
\end{eqnarray*} 
As in the estimate $\left\|V_6\right\|_{L^2(\Omega;\mathcal{H})}$, \eqref{timerr**6}--\eqref{timerr6}, we have
\begin{eqnarray}
\label{proofimpl6*1}
\left\|VI_{61}\right\|_{L^2(\Omega;\mathcal{H})}\leq C\Delta t^{\frac{2H+\beta-1}2}.
\end{eqnarray}
Let us estimate $\left\|VI_{62}\right\|_{L^2(\Omega;\mathcal{H})}$, using notation \eqref{not2}, adding and subtracting a term, we rewrite $VI_{62}$ as follows
\begin{eqnarray*}
%\label{proofimpl6**2}
VI_{62}&=&\int_0^{t_{m-1}}\left[\left(\prod_{j=\llcorner s\lrcorner^m}^{m-1} e^{-\Delta t A_{h,j}}\right)-\left(\prod_{j=\llcorner s\lrcorner^m}^{m-1} S^j_{h,\Delta t}\right)\right]P_h\phi(\llcorner s\lrcorner)dB^H(s)\nonumber\\
&=& \int_0^{t_{m-1}}\left[\left(\prod_{j=\llcorner s\lrcorner^m}^{m-1} e^{-\Delta t A_{h,j}}\right)-\left(\prod_{j=\llcorner s\lrcorner^m}^{m-1} S^j_{h,\Delta t}\right)\right]P_h(\phi(\llcorner s\lrcorner)-\phi(t_{m-1}))dB^H(s)\nonumber\\
&+&\int_0^{t_{m-1}}\left[\left(\prod_{j=\llcorner s\lrcorner^m}^{m-1} e^{-\Delta t A_{h,j}}\right)-\left(\prod_{j=\llcorner s\lrcorner^m}^{m-1} S^j_{h,\Delta t}\right)\right]P_h\phi(t_{m-1})dB^H(s)\nonumber\\
&=:& V_{621}+V_{622}.
\end{eqnarray*}
Taking the square of the $L^2$-norm and the estimate $(a+b)^2\leq 2a^2+2b^2$ yields
\begin{eqnarray}
\label{proofimpl6_2}
\left\|VI_{62}\right\|^2_{L^2(\Omega;\mathcal{H})}\leq 2\left\|VI_{621}\right\|^2_{L^2(\Omega;\mathcal{H})}+2\left\|VI_{621}\right\|^2_{L^2(\Omega;\mathcal{H})}.
\end{eqnarray}
Using \eqref{fracint2} and inserting an appropriate power of 
$A_{h,\llcorner s\lrcorner^m}$ yields
\begin{eqnarray*}
&&\left\|VI_{621}\right\|^2_{L^2(\Omega;\mathcal{H})}\nonumber\\
&=&\left\|\int_0^{t_{m-1}}\left[\left(\prod_{j=\llcorner s\lrcorner^m}^{m-1} e^{-\Delta t A_{h,j}}\right)-\left(\prod_{j=\llcorner s\lrcorner^m}^{m-1} S^j_{h,\Delta t}\right)\right]P_h(\phi(\llcorner s\lrcorner)-\phi(t_{m-1}))dB^H(s)\right\|^2_{L^2(\Omega;\mathcal{H})}\nonumber\\
&=& \int_0^{t_{m-1}}\left\|\left[\left(\prod_{j=\llcorner s\lrcorner^m}^{m-1} e^{-\Delta t A_{h,j}}\right)-\left(\prod_{j=\llcorner s\lrcorner^m}^{m-1} S^j_{h,\Delta t}\right)\right](A_{h,\llcorner s\lrcorner^m})^{\frac{1-\beta}2}\right\|^2_{L(\mathcal{H})}\nonumber\\
&\times&\left\|(A_{h,\llcorner s\lrcorner^m})^{\frac{\beta-1}2}P_h(\phi(\llcorner s\lrcorner)-\phi(t_{m-1}))\right\|^2_{L^0_2}ds.
\end{eqnarray*}
Using \lemref{prepaimpl4} (v) with $\gamma=1-\beta$, \eqref{equi3}, \assref{noiseandjump} \eqref{noise2}, the variable change $j=m-k-1$ and \cite[(169)]{Muka}, we get
\begin{eqnarray}
\label{proofimpl6_21}
\left\|VI_{621}\right\|^2_{L^2(\Omega;\mathcal{H})}
&\leq& C\Delta t^{2H+\beta-1}\int_0^{t_{m-1}}t_{m-\llcorner s\lrcorner^m-1}^{-2H}\left\|(A(0))^{\frac{\beta-1}2}(\phi(\llcorner s\lrcorner)-\phi(t_{m-1}))\right\|^2_{L^0_2}ds\nonumber\\
&\leq& C\Delta t^{2H+\beta-1}\int_0^{t_{m-1}}t_{m-\llcorner s\lrcorner^m-1}^{-2H+2\delta}ds\nonumber\\
&\leq& C\Delta t^{2H+\beta-1}\int_0^{t_{m-1}}t_{m-\llcorner s\lrcorner^m-1}^{-1+\beta}ds\nonumber\\
&\leq& C\Delta t^{2H+\beta-1}\sum_{k=0}^{m-2}\int_{t_k}^{t_{k+1}}t_{m-k-1}^{-1+\beta}ds\nonumber\\
&\leq& C\Delta t^{2H+\beta-1}\sum_{j=1}^{m-1}t_j^{-1+\beta}\Delta t\leq C\Delta t^{2H+\beta-1}.
\end{eqnarray}
Applying \eqref{fracint1}, inserting an appropriate power of 
$(A_{h,\llcorner s\lrcorner^m})$, \lemref{prepaimpl4} (v) with $\gamma=1-\beta$, \eqref{equi3}, it holds that
\begin{eqnarray*}
\left\|VI_{622}\right\|^2_{L^2(\Omega;\mathcal{H})}&=&\left\|\int_0^{t_{m-1}}\left[\left(\prod_{j=\llcorner s\lrcorner^m}^{m-1} e^{-\Delta t A_{h,j}}\right)-\left(\prod_{j=\llcorner s\lrcorner^m}^{m-1} S^j_{h,\Delta t}\right)\right]P_h\phi(t_{m-1})dB^H(s)\right\|^2_{L^2(\Omega;\mathcal{H})}\nonumber\\
&\leq& C\sum_{i\in\N^d}\left(\int_0^{t_{m-1}}\left\|\left[\left(\prod_{j=\llcorner s\lrcorner^m}^{m-1} e^{-\Delta t A_{h,j}}\right)-\left(\prod_{j=\llcorner s\lrcorner^m}^{m-1} S^j_{h,\Delta t}\right)\right](A_{h,\llcorner s\lrcorner^m})^{\frac{1-\beta}2}\right\|^{\frac 1H}_{L(\mathcal{H})}\right.\nonumber\\
&\times&\left.\left\|(A_{h,\llcorner s\lrcorner^m})^{\frac{\beta-1}2}P_h\phi(t_{m-1})Q^{\frac 12}e_i\right\|^{\frac 1H}ds\right)^{2H}\nonumber\\
&\leq& C\left(\int_0^{t_{m-1}}\Delta t^{\frac{2H+\beta-1}{2H}}t_{m-\llcorner s\lrcorner^m-1}^{-1}ds\right)^{2H}\left\|(A(0)^{\frac{\beta-1}2}P_h\phi(t_{m-1})\right\|^2_{L^0_2},
\end{eqnarray*}
furthermore, \assref{noiseandjump} \eqref{noise1}, the variable change $j=m-k-1$ and \cite[(169)]{Muka} yield
\begin{eqnarray}
\label{proofimpl6_22}
\left\|VI_{622}\right\|^2_{L^2(\Omega;\mathcal{H})}&\leq& C\Delta t^{2H+\beta-1-\epsilon}\left(\sum_{k=0}^{m-2}\int_{t_k}^{t_{k+1}}\Delta t^{\frac{\epsilon}{2H}}t_{m-k-1}^{-1}ds\right)^{2H}\nonumber\\
&\leq& C\Delta t^{2H+\beta-1-\epsilon}\left(\sum_{j=1}^{m-1}\Delta t  \hspace{0.1cm}t_1^{\frac{\epsilon}{2H}}t_j^{-1}ds\right)^{2H}\nonumber\\
&\leq& C\Delta t^{2H+\beta-1-\epsilon}\left(\sum_{j=1}^{m-1}t_j^{-1+\frac{\epsilon}{2H}}ds\right)^{2H}\leq C\Delta t^{2H+\beta-1-\epsilon}.
\end{eqnarray}
Substituting \eqref{proofimpl6_21} and \eqref{proofimpl6_22} in \eqref{proofimpl6_2} and taking the square-root hence yields
\begin{eqnarray}
\label{proofimpl6*2}
\left\|VI_{62}\right\|_{L^2(\Omega;\mathcal{H})}\leq C\Delta t^{\frac{2H+\beta-1}2-\epsilon}.
\end{eqnarray}
Adding \eqref{proofimpl6*1} and \eqref{proofimpl6*2} gives
\begin{eqnarray}
\label{proofimpl6}
\left\|VI_6\right\|_{L^2(\Omega;\mathcal{H})}\leq C\Delta t^{\frac{2H+\beta-1}2-\epsilon}.
\end{eqnarray}
\subsubsection{Estimate of $V_7$}
By adding and subtracting a term in $VI_7$, we have
{\small
\begin{eqnarray*}
&&VI_7\nonumber\\
&:=&\sum_{i=2}^{m}\int_{t_{m-i}}^{t_{m-i+1}}\int_{\chi}\left[\left(\prod_{j=m-i+2}^m U_h(t_j,t_{j-1})\right)U_h(t_{m-i+1},s)-\left(\prod_{j=m-i}^{m-1}S^j_{h,\Delta t}\right)\right]P_h\psi(z)\tilde{N}(dz,ds)\nonumber\\
&=&\sum_{i=2}^{m}\int_{t_{m-i}}^{t_{m-i+1}}\int_{\chi}\left[\left(\prod_{j=m-i+2}^m U_h(t_j,t_{j-1})\right)U_h(t_{m-i+1},s)-\left(\prod_{j=m-i}^{m-1}e^{-\Delta t A_{h,j}}\right)\right]P_h\psi(z)\tilde{N}(dz,ds)\nonumber\\
&+&\sum_{i=2}^{m}\int_{t_{m-i}}^{t_{m-i+1}}\int_{\chi}\left[\left(\prod_{j=m-i}^{m-1}e^{-\Delta t A_{h,j}}\right)-\left(\prod_{j=m-i}^{m-1}S^j_{h,\Delta t}\right)\right]P_h\psi(z)\tilde{N}(dz,ds)\nonumber\\
&=:&VI_{71}+VI_{72}.
\end{eqnarray*}
}
Thanks to \eqref{timerr*7}--\eqref{timerr7}, it holds that
 \begin{eqnarray}
 \label{proofimpl7*1}
 \left\|VI_{71}\right\|_{L^2(\Omega;\mathcal{H})}\leq C\Delta t^{\frac{2H+\beta-1}2}.
 \end{eqnarray}
Using the martingale property of stochastic integral, the It\^{o} isometry \eqref{jumpint}, inserting an appropriate power of $A_{h,m-i}$, \eqref{equi3} with $\alpha=\frac{2H+\beta-2}2$, \assref{noiseandjump} and \eqref{jump}, we have
\begin{eqnarray*}
&&\left\|VI_{72}\right\|^2_{L^2(\Omega;\mathcal{H})}\nonumber\\
 &=& \left\|\sum_{i=2}^{m}\int_{t_{m-i}}^{t_{m-i+1}}\int_{\chi}\left[\left(\prod_{j=m-i}^{m-1}e^{-\Delta t A_{h,j}}\right)-\left(\prod_{j=m-i}^{m-1}S^j_{h,\Delta t}\right)\right]P_h\psi(z)\tilde{N}(dz,ds)\right\|^2_{L^2(\Omega;\mathcal{H})}\nonumber\\
 &\leq&\sum_{i=2}^{m}\int_{t_{m-i}}^{t_{m-i+1}}\int_{\chi}\left\|\left[\left(\prod_{j=m-i}^{m-1}e^{-\Delta t A_{h,j}}\right)-\left(\prod_{j=m-i}^{m-1}S^j_{h,\Delta t}\right)\right](A_{h,m-i})^{-\frac{2H+\beta-2}2}\right\|^{2}_{L(\mathcal{H})}\nonumber\\
 &&\left\|(A_{h,m-i})^{\frac{2H+\beta-2}2}P_h\psi(z)\right\|^2\nu(dz)ds\nonumber\\
 &\leq&C\sum_{i=2}^{m}\left(\int_{t_{m-i}}^{t_{m-i+1}}\left\|\left[\left(\prod_{j=m-i}^{m-1}e^{-\Delta t A_{h,j}}\right)-\left(\prod_{j=m-i}^{m-1}S^j_{h,\Delta t}\right)\right](A_{h,m-i})^{-\frac{2H+\beta-2}2}\right\|^{2}_{L(\mathcal{H})}ds\right)\nonumber\\
  &&\int_{\chi}\left\|(A(0))^{\frac{2H+\beta-2}2}\psi(z)\right\|^2\nu(dz)\nonumber\\
   &\leq&C\sum_{i=2}^{m}\int_{t_{m-i}}^{t_{m-i+1}}\left\|\left[\left(\prod_{j=m-i}^{m-1}e^{-\Delta t A_{h,j}}\right)-\left(\prod_{j=m-i}^{m-1}S^j_{h,\Delta t}\right)\right](A_{h,m-i})^{-\frac{2H+\beta-2}2}\right\|^{2}_{L(\mathcal{H})}ds.\nonumber
\end{eqnarray*} 
As in \cite{Mukmab}, if $2H+\beta<2$, then applying \lemref{prepaimpl4} (iv) with $\gamma=2-2H-\beta$, it follows that
\begin{eqnarray*}
 \left\|VI_{72}\right\|^2_{L^2(\Omega;\mathcal{H})}
   &\leq&C\sum_{i=2}^{m}\int_{t_{m-i}}^{t_{m-i+1}}\Delta t^{2H+\beta-1-\epsilon}t_{i-1}^{-1+\epsilon}ds\nonumber\\
   &\leq&C\Delta t^{2H+\beta-1-\epsilon}\sum_{i=2}^{m}t_{i-1}^{-1+\epsilon}\Delta t \leq C\Delta t^{2H+\beta-1-\epsilon}.\nonumber
\end{eqnarray*} 
And if $2H+\beta<2$, then applying \lemref{prepaimpl4} (iii) with $\alpha_1=\frac{1-\epsilon}2$, $\alpha_2=\frac{2H+\beta-1-\epsilon}2$ yields
\begin{eqnarray*}
 \left\|VI_{72}\right\|^2_{L^2(\Omega;\mathcal{H})}
   &\leq&C\sum_{i=2}^{m}\int_{t_{m-i}}^{t_{m-i+1}}\Delta t^{2H+\beta-1-\epsilon}t_{i-1}^{-1+\epsilon}\leq C\Delta t^{2H+\beta-1-\epsilon}.\nonumber
\end{eqnarray*}
Therefore for all $\beta\in(0,1]$ and $H\in(\frac 12,1)$, we have
\begin{eqnarray}
\label{proofimpl7*2}
 \left\|VI_{72}\right\|_{L^2(\Omega;\mathcal{H})}
   \leq C\Delta t^{\frac{2H+\beta-1}2-\epsilon}.
\end{eqnarray}
Adding the estimates \eqref{proofimpl7*1}, \eqref{proofimpl7*2} gives
\begin{eqnarray}
\label{proofimpl7}
 \left\|VI_7\right\|_{L^2(\Omega;\mathcal{H})}
   \leq C\Delta t^{\frac{2H+\beta-1}2-\epsilon}.
\end{eqnarray}
Combining the estimates \eqref{proofimpl1}, \eqref{proofimpl2}, \eqref{proofimpl3}, \eqref{proofimpl4}, \eqref{proofimpl6} and \eqref{proofimpl7} yields
\begin{eqnarray}
\label{implerr}
\|X(t_m)-Y^h_m\|_{L^2(\Omega,\mathcal{H})}\leq C\Delta t^{\frac{2H+\beta-1}2-\epsilon}+C\Delta t\sum_{j=0}^{m-1}\|X^h(t_j)-Y^h_j\|_{L^2(\Omega,\mathcal{H})}.
\end{eqnarray}
Applying discrete Gronwall's inequality to \eqref{implerr} and adding to \eqref{boundspaerr} completes the proof.$\hfill\square$.

\end{document}